\documentclass[floatfix, nofootinbib, letterpaper, onecolumn, showpacs]{revtex4}

\usepackage[total={6.9in,9.3in}, top=0.8in, left=0.8in, includefoot]{geometry}

\usepackage{amsmath,amsthm, amssymb}
\usepackage{amsfonts,enumerate,amsbsy}
\usepackage{subcaption}
\usepackage{mathrsfs}
\usepackage{mathtools}
\usepackage{tikz,tikzscale}
\usepackage{graphicx,caption}
\usetikzlibrary{decorations.markings}
\usetikzlibrary{graphs,graphs.standard,arrows.meta}
\usetikzlibrary{arrows,shapes,cd}
\usetikzlibrary{positioning}
\usepackage[nodisplayskipstretch]{setspace}
\usepackage{environ}
\usepackage{pgfplots,pgfplotstable,pgf}
\usepgfplotslibrary{groupplots}
\pgfplotsset{compat=newest, ticks=none}
\usepackage[final]{pdfpages}
\usepackage{import}
\usepackage{xifthen}
\usepackage{transparent}
\usepackage{enumerate}
\usepackage{xparse}
\usepackage{url} 
\usepackage[pdftex,bookmarks, colorlinks, citecolor =blue, breaklinks]{hyperref} 

\usepackage[capitalise]{cleveref}
\crefname{section}{\S}{\S}

\makeatletter
\newsavebox{\measure@tikzpicture}
\NewEnviron{scaletikzpicturetowidth}[1]{%
	\def\tikz@width{#1}%
	\begin{lrbox}{\measure@tikzpicture}%
		\BODY
	\end{lrbox}%
	\pgfmathparse{#1/\wd\measure@tikzpicture}%
	\BODY
}
\makeatother
\tikzset{->-/.style={decoration={
			markings,
			mark=at position #1 with {\arrow{>}}},postaction={decorate}}}
\tikzset{cross/.style={cross out, draw=black, minimum size=2*(#1-\pgflinewidth), inner sep=0pt, outer sep=0pt},
	cross/.default={3pt}}

\NewDocumentCommand\set{mg}{%
	\ensuremath{\bigl\{ #1 \IfNoValueTF{#2}{}{\bigm| #2} \bigr\}}%
}

\newcommand{\grad}[1]{\nabla #1} 
\renewcommand{\div}[1]{\nabla \cdot #1} 
\newcommand{\curl}[1]{\nabla \times #1} 

\newcommand{\me}{\mathrm{e}}
\newcommand{\R}{\mathbb{R}}

\newcommand{\N}{\mathbb{N}}

\newcommand{\T}{\mathbb{T}}

\newcommand{\cC}{\mathcal{C}}
\newcommand{\cF}{\mathcal{F}}
\newcommand{\cL}{\mathcal{L}}
\newcommand{\cM}{\mathcal{M}}
\newcommand{\cP}{\mathcal{P}}

\newcommand{\g}{\gamma}

\DeclareMathOperator{\Span}{Span}

\theoremstyle{plain}
\newtheorem{thm}{Theorem}[section]

\newtheorem{lemma}[thm]{Lemma}

\theoremstyle{definition}
\newtheorem*{definition}{Definition}

\begin{document}
\title [Invariant Tori Transverse to Foliations]
      {Nonexistence of Invariant Tori Transverse to Foliations:\\ 
        An Application of Converse KAM Theory}

\author{Nathan Duignan and James~D. Meiss}\thanks 
      {
        The authors acknowledge support from the Simons Foundation grant \#601972 
        ``Hidden Symmetries and Fusion Energy."
        Useful conversations with Robert MacKay and Josh Burby are gratefully acknowledged. 
      }

\affiliation{Department of Applied Mathematics, University of Colorado, Boulder, CO, USA}

\date{\today}

\begin{abstract}
\vspace*{1ex}
\noindent
	Invariant manifolds are of fundamental importance to the qualitative understanding of dynamical systems.
	In this work, we explore and extend MacKay's converse KAM condition to obtain a sufficient 
	condition for the nonexistence of invariant surfaces that are transverse to a chosen 1D foliation. 
	We show how useful foliations can be constructed from approximate integrals of the system. 
	This theory is implemented numerically for two models,
	a particle in a two-wave potential and a Beltrami flow studied by Zaslavsky (Q-flows). 
	These are both 3D volume-preserving flows, and they exemplify
	the dynamics seen in time-dependent Hamiltonian systems and incompressible fluids, respectively. Through both numerical and theoretical considerations, it is revealed how to choose foliations that
	capture the nonexistence of invariant tori with varying homologies. 
\end{abstract}
	
\maketitle
\section*{}
\textbf{
Kolmogorov-Arnold-Moser (KAM) theory provides conditions under which a nearly-integrable Hamiltonian system or symplectic map can be guaranteed to have invariant tori on which the dynamics is conjugate to rigid rotations. Dynamical behavior on these invariant sets is just like that of an integrable system, for which the phase space is almost everywhere foliated by such tori. By contrast, \textit{converse KAM} theory, as initiated by Mather in 1984, seeks to obtain conditions for the nonexistence of invariant tori. Mather's 
idea is based on the fact that if a two-dimensional map has an invariant circle, then the iterates of an infinitesimal vector that is attached to a point on the circle cannot rotate through the tangent space of the circle. A converse KAM condition then becomes: no such circle can exist if a vector rotates ``too much".   This condition was also extended to flows by MacKay, who gave a condition for the nonexistence of 
``rotational" invariant tori on three-dimensional energy surfaces of a two degree-of freedom Hamiltonian system.
\\
Here we use MacKay's ideas to explore the definition of rotation \textit{relative to a choice of foliation}---a 
family of curves that cover the manifold. A codimension-one
torus that is transverse to the foliation is ruled out when a positive tangent vector to the foliation flows past a negative one--exceeding a half rotation. Selection of the foliation is crucial in detecting the nonexistence of tori 
that might have different topological structures. 
We explore choosing foliations based on gradients of approximate
invariants for several examples of flows in three dimensions.
}
\section{Introduction}

Mather's converse KAM theory \cite{Mather84, Mather86} is based on a theorem due to Birkhoff: for an
area and orientation-preserving twist (or more generally, tilt) map, 
$f: M \to M$, on $M = \T \times \R$, every rotational
invariant circle is a Lipschitz graph over the angle variable \cite{Mather84,Meiss92}. Recall that a set
$\cC \subset M$ is invariant when $f(\cC) = \cC$ and it is a \textit{rotational} circle if
it is homotopic to the horizontal circle $\{(x,0)\,|\, x \in \T\}$. The
implication is that, if $v \in TM_z$ is a tangent vector at a point $z = (x,y)
\in M$ to a rotational invariant circle $\cC$, then it is contained in a
Lipschitz cone.

Constraining the tangent vectors to a Lipschitz cone can be turned into a nonexistence condition in the following way. Consider
a point $z_0 \in M$ and an arbitrary vector $\xi_0 \in T_{z_0}M$. If there is a
rotational invariant circle $\cC$ that contains $z_0$, then under iteration by linearized map
$Df_{z_t}: T_{z_t}M \to T_{z_{t+1}} M$, this vector $\xi_0$ cannot cross the circle. Indeed, since
the tangent to the circle lies in the Lipschitz cone, the image at time $t$, $\xi_t := Df_{z_{t-1}}\xi_{t-1}$, cannot cross this cone.
The simplest manifestation of this is that an initially vertical vector $\xi_0 = (0,1)^T$ can
never rotate past the negative vertical $(0,-1)^T$. For example, when the twist is
positive, vertical vectors tilt to the right, so if
there is a $t$ for which $\xi_t$ is in the third quadrant, there is no invariant
circle through $z_0$. This gives a sufficient condition for nonexistence of a
rotational circle. 

Mather used this nonexistence condition, together with explicit Lipschitz
bounds, to show that there are no rotational invariant circles for Chirikov's
standard map when $k > 4/3$; this amazing result is not too far from the 
numerically obtained value
$k_{cr} \approx 0.971$, for the last such circle \cite{Greene79}. Using interval
arithmetic, MacKay and Percival \cite{MacKay85} showed nonexistence for $k > 63/64$. 
Subsequent calculations have improved this bound \cite{Jungreis91}.

This converse KAM condition has been generalized to higher-dimensional, symplectic
maps with nondegenerate twist \cite{mackay_converse_1989} and monotone-positive, exact-symplectic maps \cite{Haro99},\footnote
{If $Df$ is represented in block form by four $n\times n$ matrices $A,B,C,D$, then the map has nondegenerate twist if $B + B^T$ is positive definite, and is monotone positive when $B^{-1}A$ and $D B^{-1}$ are positive definite}
as well as to Hamiltonian systems on a cotangent bundle that satisfy a Legendre condition
\cite{mackay_criterion_1989}. In these cases, the criteria can rule out the
existence of tori that are Lagrangian graphs. Note, however, that for the higher-dimensional
case there is no generalization of Birkhoff's theorem: it is not clear that all rotational tori
are indeed Lagrangian graphs, even when an analogue of the twist condition is
satisfied.

Other applications include a study of the tilting of a rigid body in an elliptical 
Kepler orbit for both conservative and dissipative cases \cite{Celletti07}. 
This system is effectively three-dimensional, and the authors use the conjugate point technique of
\cite{mackay_criterion_1989} to compute regions of phase space that cannot contain rotational tori.
Ideas similar to the converse KAM theorem show that in a ``natural Hamiltonian system'' 
(of kinetic plus potential energy form),
any compact energy surface cannot be completely foliated by invariant $n$-tori \cite{Knauf90}.
Another result for these natural systems gives a criterion for the nonexistence of so-called 
``viscosity solutions'' to the Hamilton-Jacobi equation \cite{Gomes08}.

White introduces an analogous condition in a study of guiding center dynamics \cite{whiteModificationParticleDistributions2011, whiteModificationParticleDistributions2012, whiteDeterminationBrokenKAM2015}. This two degree-of-freedom Hamiltonian system can be reduced to an area-preserving map on an energy surface by taking a Poincar\'{e} section. In this work, instead of utilizing the tangent map idea of Mather, White computes a nearby trajectory $z^*_t$ to $z_t$ and studies the finite-size vector $\Delta z_t^* = z^*_t - z_t$. A similar converse KAM theorem still holds; if $\Delta z_t^*$ rotates by more than $\pi$, then $z_t$ does not lie on a rotational invariant torus. White calls his method ``phase vector rotation''. Analyzing this finite-size vector $\Delta z_t^*$ has a deficiency that Mather's tangent map does not. Consider the dynamics of a pendulum and suppose that $z_t$ is on a librational torus, so it oscillates about the stable equilibrium, but that $z_t^*$ is instead chosen to be on the other side of the separatrix, on a rotational torus. Then, in such a situation, $\Delta z_t^*$ will not necessarily rotate by more than $\pi$ --- it could remain more-or-less vertical
throughout the evolution. Thus, $z_0$ would not satisfy the phase vector rotation condition and hence it could not be concluded to not lie on a rotational torus. The implication is that one must be careful in interpreting the results of this method. We will use below a useful diagnostic that White introduces: a plot of the effective rotation rate as a function of initial condition.

A reformulation of the converse KAM condition was given by MacKay \cite{mackay_finding_2018} for 
three-dimensional flows, the context that we will be working on in this paper. 
In particular, we will numerically explore a key result of MacKay
\cite[thm.~1]{mackay_finding_2018}. Given a flow $\phi$ generated by a
$C^1$ vector field $v$ on some $n$-dimensional manifold $M$, MacKay finds
conditions under which a given point $x\in M$ does not lie on an invariant,
codimension-one, orientable surface $S$. The class of surfaces considered are those
that are everywhere transverse to the leaves of a given one-dimensional,
orientable foliation $\cF$ of $M$. Recall that a $k$-dimensional foliation is
an equivalence relation on the manifold $M$, where the equivalence classes, called
leaves, are injectively immersed submanifolds of dimension $k$. Essentially
$\cF$ is a collection of disjoint $k$-dimensional submanifolds whose
union is the manifold $M$.

Given a choice of foliation $\cF$, we use the theorem of MacKay to construct a
numerical scheme that detects the nonexistence of invariant surfaces
transverse to the foliation. This scheme is applied to two models. The first is
the \emph{two-wave model}, which is perhaps the simplest nonintegrable,
nonlinear, $1\tfrac12$ degree-of-freedom Hamiltonian system. This model
was also studied by MacKay using the vertical foliation \cite{mackay_criterion_1989};
here we will compare several foliations. The
second model is a family of incompressible, Beltrami flows that were called \emph{Q-flows} by 
Zaslavsky \cite{zaslavsky_dynamical_1991}. These
vector fields can have crystalline or quasi-crystalline spatial symmetries. 
Zaslavsky studied anomalous diffusion in these systems, 
but converse KAM theory has not been previously applied.

Using these two models, we investigate methods for choosing an appropriate 1D foliation
and reveal some of the consequences of such a choice. In particular, if the
1D foliation is generated from the gradient flow of a function $J:M \to \R$ , we
are able to test for the nonexistence of invariant tori
that are not of null-homology on the manifold $M\setminus\Sigma$, where
$\Sigma$ is the set of singular leaves of $J$.

In \S\ref{sec:Background} we introduce the concepts and background for converse
KAM theory. The main result of MacKay, Thm.~\ref{thm:MacKay}, is recalled and a
simplification for Cartan-Arnol'd type systems is obtained.
The two models that we study are described in \S\ref{sec:Examples}. 
In \S\ref{sec:Foliations} we discuss the construction of a number of foliations 
for each model and give a description of techniques used to numerically
implement Thm.~\ref{thm:MacKay}.
The results of numerical experiments are given in \S\ref{sec:results} and discussed in \S\ref{sec:discussion}.

\section{Converse KAM Theorem: Background and Theory}\label{sec:Background}

In this section we will expand on some details of \cite{mackay_finding_2018},
prove several helpful results, and discuss the generation of foliations from the gradient flow of an approximate integral.

We will study the flow $\phi$ of a $C^1$ vector field $v$ on an
orientable manifold $M$: 
\begin{equation}\label{eq:flow}
	\frac{d}{dt} \phi_t(x) = v(\phi_t(x)), \quad \phi_0(x) = x ,
\end{equation}
That is, we assume that $\phi$ is a complete flow on $M$. In our applications, $M$ will be three-dimensional
and $v$ will either be Hamiltonian or divergence-free. However, at first,
we more generally assume that $M$ is $n$-dimensional and $v$ is a general vector field on $M$.

\subsection{Surfaces Transverse to a Foliation}\label{sec:MainTheorem}

Let $\cF$ be a given 1D foliation on $M$ and denote the leaf at a point $x$ by $\cF_x$. 
We begin by recalling the definition of a transverse surface to $\cF$.

\begin{definition}[Transverse Surfaces]
  Let $\cF$ be a one-dimensional foliation of a manifold $M$.
  A codimension-one surface $S$ in $M$ is \textit{transverse} to $\cF$
  if, for every point $x\in S$, the tangent space $T_xS$ does not contain the
  tangent space of the foliation $\cF$ at $x$: $T_x\cF \nsubseteq T_xS$.
\end{definition}
Locally, there always exists an $n-1$ dimensional submanifold $B$ transverse to
the leaves of $\cF$. In other words, $\cF$ is locally
diffeomorphic to the trivial fibre bundle $B\times \R$. By definition, a surface
$S$ is locally transverse to the leaves of $\cF$ if and only if
$S$ is locally a graph over $B$. This follows from the converse of the implicit
function theorem; indeed, if there is some point $x \in M$ for which the tangent space
$T_xS$ contains the tangent space of the foliation at $x$, $T_x\cF$
then $S$ is not locally a $C^1$ graph over $B$. 

When $S$ is invariant under the flow \eqref{eq:flow}, we
can use the vector field $v$ to help detect transversality:

\begin{lemma}
  \label{lem:Conditions}
  Suppose that
  $S$ is a surface that is invariant under \eqref{eq:flow}; 
  that is, for each $x \in S$, $v(x) \in T_xS$.
  Let $\eta$ be any nonzero vector field in $\cF$ and
  $y$ be a point in $S$ for which $v(y) \neq 0$. Then if either
 \begin{enumerate}
   \item $\eta_y = c v(y)$, for some $c \neq 0$ or,
  \item  there exists some $s \in T_yS$ independent of $v(y)$ such that $\eta_y =
    s + c v(y)$,
 \end{enumerate}
  then $S$ is not transverse to $\cF$.
\end{lemma}
\begin{proof}
  Since invariance implies that $v(x)\in T_xS$, 
  if there is a vector $s\in T_yS$
  such that $ \eta_y= s + c v(y)$ for some $c \in \R$ then $T_y\cF_y\subset T_yS$. 
  This is true whenever either (1) or (2) is true.
\end{proof}

Note that if the first condition in Lem.~\ref{lem:Conditions} is true, it follows that $\cF_y$ is
not transverse to the vector field $v(y)$. This condition is easy to check
analytically. Unfortunately, the second condition can only be checked if we already know
$S$. However, knowing $S$ in advance is clearly inconsistent with our goal of deciding
if such a surface actually exists! What is needed is a systematic way of ruling out
potential vectors $s$.

Under the assumption that the invariant surface $S$ and the foliation $\mathcal{F}$ are both
orientable we can make progress toward eliminating
potential tangent vectors.
We start by choosing some volume form $\Omega$ in a tubular neighbourhood of $S$. The existence of $\Omega$ is guaranteed as $S$ is orientable, thus, the tubular neighbourhood is orientable, and all orientable manifolds have a volume form. This form need
not be invariant under the flow $\phi$.
Take a point $x_0 \in S$ and let $x_t:=\phi_t(x_0)$ denote its orbit.
At each point $x_t$ consider a basis $s_i^t, i = 1\ldots n-2$, so that
\begin{equation}\label{eq:tangentSpace}
	 \Span\{s_1^t,\dots,s^t_{n-2},v(x_t)\} = T_{x_t}S, 
\end{equation}
and let $\eta_t\in T_{x_t}\cF$ be a nonzero vector field tangent to the orientable foliation. 
We assume that both the bases $s_i^t$ and the vector $\eta_t$ are at least $C^0$ for $t\in I \subset \R$ some interval.

If the orientable surface $S$ is transverse to $\cF$, then an orientation of $\eta_t$ can be chosen such that
\begin{equation}
  \label{eq:etacondition}
  \Omega(\eta_t,s_1^t,\dots,s_{n-2}^t,v(x_t)) > 0,\quad \forall t\in I.
\end{equation}
Indeed, continuity of the basis and vector field $\eta_t$
implies that the sign of \eqref{eq:etacondition} must remain constant.
If condition \eqref{eq:etacondition} were to
fail, then $\cF$ is not transverse to the tangent
space \eqref{eq:tangentSpace} at the point $x_t$, and so
the surface $S$ cannot be transverse to $\cF$.

\begin{figure}[ht]
  \centering
   \includegraphics[width=0.4\textwidth]{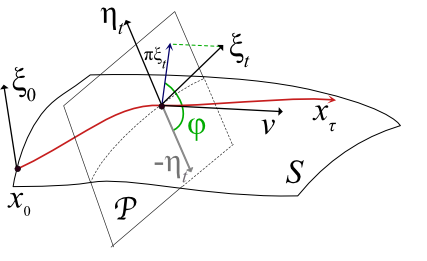}
   \caption{A sketch of a potential invariant surface $S$, the evolution of a vector $\xi_0$,
   to time $t$, and the local tangent $\eta_t$ to the foliation at $x_t$. A surface $\cP$
   transverse to $v(x_t)$ gives a projection $\pi\xi_t$ along $v$,
   and the angle $\varphi$ \eqref{eq:angle}. Here, $\cP$ is chosen
 to contain $\cF$.}
  \label{fig:Sketch}
\end{figure}

Since $S$ is also invariant, a second requirement can be deduced.
Let $\xi_0\in T_{x_0}\cF$ and
\begin{equation}{eq:linearization}
	\xi_t := \phi_{t*}\xi_0  = D_{x_0}\phi_t(x_0) \xi_0 \in T_{x_t}M;
\end{equation}
be the pushforward of $\xi_0$ under the flow $\phi_t$, as sketched in \cref{fig:Sketch}.
Note that, in general, $\xi_t \notin T_{x_t}\cF$. The
invariance of $S$, together with orientability, guarantee that
\begin{equation}
  \label{eq:xicondition}
  \Omega(\xi_t,s_1^t,\dots,s^t_{n-2},v(x_t)) > 0.
\end{equation}
If condition \eqref{eq:xicondition} were to fail then either
$S$ is not orientable or is not transverse to $\cF$.

These two conditions, \eqref{eq:etacondition} and \eqref{eq:xicondition}, can be
combined to imply
\begin{equation}
  \label{eq:etaxicondition}
  \Omega(c\eta_t + \xi_t,s_1^t,\dots,s^t_{n-2},v(x_t)) >0,\qquad \forall c >0.
\end{equation}
From the previous arguments, failure of \eqref{eq:etaxicondition} to hold
implies that $x$ does not lie on an invariant, orientable surface $S$ transverse to $\cF$ and 
with the given tangent space $\Span\{s_1^t,\dots, s_{n-1}^t\}$.

In particular, if there is a time  $t\in I$, such that $\xi_t = -c\eta_t + s$ where
$\eta_t\in T_{x_t}\cF$ and $s\in T_{x_t}S$, then condition
\eqref{eq:etaxicondition} will fail. As we don't know a priori if such an $s$
exists, all that can be checked is the case $s = d v(x_t)$, that is, when
$\xi_t = -c\eta_t + d v(x_t)$ for some $c > 0, d \in \R$.

This last result is still not computationally friendly, as one would need to
check all possible coefficients $c$ and $d$ for each $t\in I$. To help avoid this,
suppose that we can construct some $n-1$ dimensional foliation $\cP$
that is transverse to $v$: we might think of $\cP$ as a
\textit{Poincar\'e foliation} as it could lead to a local Poincar\'e section,
recall \cref{fig:Sketch}. If we define the projection along $v(x)$ by
$\pi:T_xM\to T_x\cP_x$, then the condition $\xi_t = -c \eta_t + d v(x)$
becomes
\begin{equation}
  \label{eq:condition}
  \pi \xi_t = -c \pi\eta_t.
\end{equation}

If $M$ has a Riemannian metric, the corresponding inner product $\langle \cdot , \cdot \rangle$ on each tangent
space can be used to instead consider the angle
\begin{equation}
  \label{eq:angle}
  \varphi_t := \cos^{-1}\left( \frac{\langle\pi\xi_t, -\pi\eta_t\rangle}{\|\pi\xi_t\|\|\pi\eta_t\|} \right)
   \in [0,\pi],
\end{equation}
as sketched in \cref{fig:Sketch}.
Now \eqref{eq:condition} is equivalent to $\varphi_t = 0$. Essentially we have
the following extension of \cite[Thm.~1]{mackay_finding_2018}.

\begin{thm}\label{thm:MacKayExtended}
  Assume that $M$ is a  manifold with dimension $n \ge 3$, $\cF$ is an orientable 1D foliation,
  and $\cP$ is a codimension-one Poincar\'e foliation containing $\cF$. Given $x\in M
  $, a positively oriented tangent $\xi$ to $\cF_x$ and the
  projection $\pi$, if there is some $t\in I\subset \R$ for which $\varphi_t = 0$, then
  $x$ does not lie on any oriented, invariant $n-1$ dimensional sub-manifold of
  $M$ that is transverse to $\cF$.
\end{thm}

However, checking the hypotheses of Thm.~\ref{thm:MacKayExtended}
is only reasonable for $n=3$. The issue is one of
codimension: since each leaf of $\cP$ is $n-1$ dimensional, checking
that $\pi\xi_t = -c\pi\eta_t$ is equivalent to checking if a point crosses a
codimension $n-2$ surface. Numerically, this can only be ascertained if the
codimension of the crossed surface is $1$, that is, when $n = 3$.

When $n = 3$, we can refine Thm.~\ref{thm:MacKayExtended}. Let,
\begin{equation}\label{eq:theta}
	\theta_t  = \angle(\pi\xi_t,-\pi\eta_t)\in (-\pi,\pi]
\end{equation}
be the oriented angle between
$\pi\xi_t$ and $-\eta_t$. Note that computing the full-range angle is only possible if the two
vectors are always in some given two-dimensional plane.
If this angle crosses zero at some time, then we must have
$\pi\xi_t = -\pi\eta_t$ at an earlier time. This leads to a mild extension of the original theorem of
\cite{mackay_finding_2018}.

\begin{thm}{\cite{mackay_finding_2018}}
  \label{thm:MacKay}
  Assume that $M$ is a 3D manifold. Given $x\in M$, positively oriented tangent
  $\xi$ to an orientable foliation $\cF$ at $x$ and the projection $\pi$, if there is a
  $t \in I\subset\R$ such that the oriented angle $\theta_t$ \eqref{eq:theta} changes sign,
  then $x$ does not lie on any oriented, invariant 2D
  sub-manifold of $M$ that is transverse to $\cF$.
\end{thm}

There is another method for adequately detecting nonexistence of invariant surfaces through the failure of condition \eqref{eq:etaxicondition}. If a 3D manifold $M$ can, in addition, be assumed to be orientable, then necessarily there is some globally defined volume form $\Omega$. Checking whether at some time $t_c$ that $\xi_{t_c} = -c \eta_{t_c} + d v(x_{t_c})$ for some $c,d\in \R$ is equivalent to finding a time $t_c$ in which the volume,
\begin{equation}
	K(t) =	\Omega(v(x_t), \xi_t,\eta_t ),
\end{equation}
vanishes. In order for \eqref{eq:etaxicondition} to then be satisfied, it is required that $c \geq 0$. This can be guaranteed provided $K(t_c) = 0$ whilst the inner product $\langle \eta_{t_c},\,\xi_{t_c} \rangle < 0$.
In summary, we have shown the following additional nonexistence condition.

\begin{thm}
	\label{thm:CKAMVolume}
	Assume that $M$ is a 3D, orientable manifold, $x_0\in M$, 
	 $\xi_0 = \eta_0$ is a positively oriented tangent to $\cF$ at $x_0$, and $\Omega$
	is a volume form. If for some time $t \in I\subset\R$,
	$K(t):= \Omega(v_t,\eta_t,\xi_t)$ changes sign while $\langle \eta_t,\,  \xi_t \rangle < 0$
	then $x_0$ does not lie on any oriented, invariant 2D
	sub-manifold of $M$ that is transverse to $\cF$.
\end{thm}

\subsection{Converse KAM for a Cartan-Arnol'd System}\label{sec:specialCases}
One difficult facet of Thm.~\ref{thm:MacKay} or Thm.~\ref{thm:CKAMVolume} is the construction of an appropriate
foliation $\cP$ or of the volume form $\Omega$. 
With some additional structure on $M$ this problem can
be avoided. Specifically, as pointed out in \cite{mackay_finding_2018}, this
can be done if one has a so-called Cartan-Arnol'd one-form for the vector field $v$ on $M$.

\begin{definition}
  A one-form $\alpha$ for a vector field $v$ on a 3-manifold $M$ is
  \emph{Cartan-Arnol'd} if $\iota_vd\alpha = 0 $. If such a one-form exists
  the differential equation associated to $v$ is said to be a \emph{Cartan-Arnol'd system}.
\end{definition} 

Such forms naturally arise in many physical problems. For example, the
Poincar\'{e}-Liouville one form,
\begin{equation}
  \label{eq:PoincareLiouville}
  \alpha = p dq - H(q,p,t) dt ,
\end{equation}
of a time-periodic, one degree-of-freedom Hamiltonian system $H(q,p,t)$
(i.e. a system with $1\tfrac12$ degrees of freedom)
is Cartan-Arnol'd.
Here $M$ is the product of the 2D phase space and a circle for time $t$. 
A similar situation pertains to a two degree-of-freedom Hamiltonian system 
restricted to any regular, compact energy level, $H = E$.

Volume-preserving flows on $\R^3$ give another example:
for any divergence free vector field $B$, let $A$ be the vector
potential such that $B= \nabla\times A$. In this case, 
$\alpha = A\cdot dx$ is a
Cartan-Arnol'd form; indeed $i_B d\alpha = 0$ because $B \times B = 0$.

More generally, any Cartan-Arnol'd vector field has a natural, conserved volume form
\begin{equation}
  \label{eq:VolumeForm}
  \Omega = \tfrac{1}{|v|^2} v^\flat \wedge d\alpha.
\end{equation}
Here, if one takes the Riemannian metric $g $ to be the Euclidean metric, then $v^\flat = v \cdot dx$ is the one-form induced by $v$, so that $i_v v^\flat = |v|^2$.
For this choice, since $i_v d\alpha = 0$ by definition, then $i_v \Omega = d\alpha$.
As a consequence $\cL_v \Omega = 0$, so $v$ is incompressible with respect to $\Omega$.

Given a foliation with tangent vector $\eta$, and  taking $\xi_t,\eta_t$ as in \S\ref{sec:MainTheorem}, 
the sign change of $\varphi_t$ \eqref{eq:theta} for $t \in I$, some interval, 
can be determined by a change in sign of
\begin{equation}
  \label{eq:KCond}
  K(t) = \Omega(v_{\phi(t)}, \xi_t, \eta_t) = d\alpha(\xi_t,\eta_t).
 \end{equation}
Geometrically $K(t)$ represents an area. Indeed, note that, since
$\iota_vd\alpha = 0$, then for any point $x_t$
\begin{equation}
  d\alpha(\xi_t,\eta_t) = d\alpha(\pi\xi_t,\pi\eta_t),
\end{equation}
where $\pi$ is the projection along $v$ onto any complementary subspace $W$ to
$\ker(d\alpha)$. In particular, $W$ could be the tangent space of some Poincar\'{e} foliation $\cP$. 
For any $W$, $d\alpha$ is an area
form (thus a volume form) on $W$. Hence, $K(t) = 0$ implies that the area
spanned by $\pi\xi_t$ and $\pi\eta_t$ is zero, that is, $\pi\xi_t = \pm \pi
\eta_t$. What is remarkable is that $K(t)$ gives an area independent of the choice of $W$.
Thus, there is no need to construct a Poincar\'{e} foliation. 

As pointed out in \cite{mackay_finding_2018}, the only issue with using $K$ to
detect sign changes is that this also happens when $\pi\xi_t = + \pi\eta_t$. 
In line with Thm.~\ref{thm:CKAMVolume}, one way to ensure that $\varphi_t$ crosses $0$ over
some interval $I$ is by requiring that $\langle \eta_t, \xi_t\rangle < 0$. 
We will make use of this fact in the computations throughout the paper.

\subsection{Choosing Foliations}
A crucial prerequisite for the results in the previous subsections is the choice of a one-dimensional
foliation $\mathcal{F}$. It is clear, from the definition of tranversality, that
the converse KAM condition is dependent on this choice (see also \S\ref{sec:results}).
We will say that a point lying on an invariant surface that nevertheless satisfies the converse KAM condition
with respect to a given foliation, does so (only) \emph{dependently}, 
and that point is a \emph{dependent} point.

One goal is to select a foliation such that there are no dependent points; in this
case if a point satisfies converse KAM then it is guaranteed to not lie on an invariant surface.
For example, if the Poincar\'e map of a $1\tfrac12$
degree-of-freedom Hamiltonian system satisfies the twist condition, then its
rotational tori are graphs, and we can choose the vertical foliation in order
to rule these out. This is what was done in the original studies of twist maps
\cite{Mather84,MacKay85} and for the flows studied in
\cite{mackay_criterion_1989}. In this case, every point that satisfies the
converse KAM condition cannot lie on a rotational invariant torus. Such a
foliation will not, however, correctly detect librational tori, such as those in
the island chains of twist maps. Thus there can be dependent points that satisfy converse KAM condition,
even though they lie on librational tori. The cited papers did not address how to select the best
foliation that could detect such structures.

One way to build a foliation is to choose a smooth function $J:M \to \R$, and define
a one-dimensional foliation by the integral curves of the 
gradient vector field $\eta = \nabla J$.
A nice feature is that the hypotheses of Thm.~\ref{thm:MacKay}
only require the vector field $\eta$, and not its integral curves, so one need not actually
construct $\mathcal{F}$: it is sufficient to have $\nabla J$.
Such a foliation will have singularities on critical sets of $J$,
where $\eta = 0$. Any level set of $J$ containing such a critical point is a 
singular leaf of $J$. Removing these singular leaves from the phase space $M$, leads
to additional structure, that (as we discuss below) can distinguish between
tori of different topological classes. For such a foliation, we have to be a
bit more careful in our discussion of the implications of the converse KAM condition.

One way to ensure that $\mathcal{F}$ be predominantly transverse to invariant surfaces
of a vector field $v$ is to choose $J$ to be an approximate invariant of $v$, i.e.,
$v \cdot \nabla J \approx 0$. Indeed, if $J$ were a true invariant of the system, then its gradient
would necessarily be transverse to any invariant tori.
Such a function can be easily found when there is a small parameter so that $v$ is nearly integrable or has an
adiabatic invariant. We will use this idea in \S\ref{sec:Foliations}. 

Having more information about the foliation $\mathcal{F}$ allows one to
determine what types of surfaces on which a dependently satisfying point may lie. Let
us focus on determining the existence of invariant tori, a question pertinent to
Hamiltonian and volume-preserving systems. Let $\Sigma$ be the set of all
singular leaves of $J$ and consider $\cM = M\setminus \Sigma$. 
When $\Sigma$ is nonempty, the homology of $\cM$ can be
distinct from $M$.

There are three distinct issues that we must consider when we find
a point $x$ that satisfies the converse KAM condition for $J$.
Firstly, it is possible that $x$ lies on an invariant torus $S$ of null-homology on
$\mathcal{M}$. Since $S$ is contractable in $\cM$, but the levels sets of $J|_\cM$ are not,
such a torus cannot be everywhere transverse to the
foliation generated from $J$. In this case, $x$ only satisfies
the converse KAM condition dependently on the choice of $J$.

Secondly, it is possible that $x$ lies on a torus $S$ with nonzero homology on $\cM$.
In this case, it could still satisfy the converse KAM condition dependently when $S$
is not everywhere transverse to the foliation, i.e. $S$ is not locally a graph over 
level sets of $J$.

Finally, it is also possible that an invariant torus of $v$ intersects
the set of singular leaves $\Sigma$.  This could be the case when $J$ is an approximate
invariant: tori could still exist near $\Sigma$. In such a circumstance, it is more difficult to say whether a
point $x$ on this invariant surface will be dependent or whether it will not
pass the converse KAM condition.

\begin{figure}[ht]
    \centering
    \includegraphics[width=0.3\textwidth]{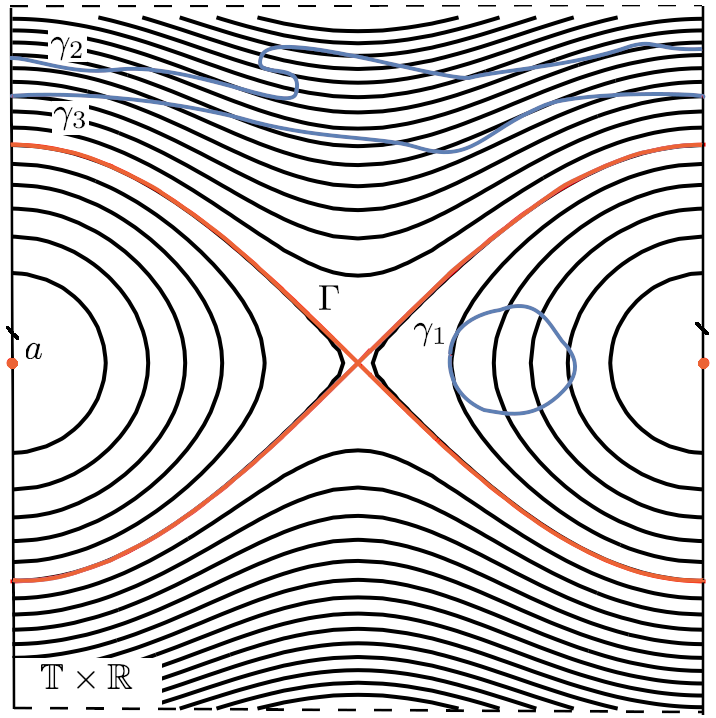}
    \caption{Sketch of the intersection of potentially invariant tori of \eqref{eq:perturbPend}
    on a Poincar\'e section of $M$, a cylinder formed by identifying the right and left edges, at fixed time. 
    Here the singular leaves of $J$ on the section are $\Sigma = \{a\}\cup\Gamma$, 
    where $\Gamma$ is the separatrix}
    \label{fig:homology}
  \end{figure}

Consider, for example, the case of a $1\tfrac12$ degree-of-freedom Hamiltonian 
\begin{equation}\label{eq:perturbPend}
	H_\epsilon (q,p,t) = \tfrac12 p^2 - \mu \cos(2\pi q) + \epsilon H_1(q,p,t).
\end{equation}
Here we assume the time-dependence is periodic, so that $M = \T \times \R \times\T$.
In this case each constant $t$ surface will be a global
Poincar\'{e} section. When $\epsilon \ll 1$, the dynamics is a
weakly-perturbed pendulum, thus a reasonable choice for  $J:M \to
\R$ is one that when restricted to $t=0$ gives the pendulum $J|_{t=0} = H_0$. 
The intersection of the level sets of $J$ with the Poincar\'{e} section are sketched in
\cref{fig:homology}. 
The singular leaves of $J$ intersect the Poincar\'e section at 
the elliptic point $a$ and the separatrix $\Gamma$, the red sets in the figure.
Removing $\Sigma = \{a\} \cup \Gamma$ from $M$ yields a manifold $\cM$ of homology distinct from $M$:
for example, tori that enclose $a$ are no longer contractable.
Three curves have been sketched in \cref{fig:homology} to demonstrate cross sections
of three potentially invariant tori of $v$. The torus $\g_1$ has
null-homology in $\cM$ and hence any point $x\in \g_1$ that satisfies the converse KAM condition
does so dependently. Such a situation would pertain if, for example, the perturbed
system $H_\epsilon$ still has an elliptic periodic orbit, but it moves away from $a$. 

The rotational tori $\g_2$ and $\g_3$ are not of null-homology in $\cM$. However, $\g_2$ does not
remain transverse to the foliation generated by $J$, and thus will contain
dependent points. This situation can happen if the dynamics of \eqref{eq:perturbPend}
no longer satisfied a twist (or tilt) condition.
Finally $\g_3$ is everywhere transverse to the foliation,
and thus no orbit on this set will satisfy the converse KAM theorem.

\section{Examples}\label{sec:Examples}

In this section, we describe two volume-preserving systems and derive several
foliations generated from adiabatic invariants or approximate integrals.
These systems will be tested against the converse KAM condition with
respect to the generated foliations in \S\ref{sec:results}.

\subsection{Two-Wave Model}\label{sec:DoubleWave}
As a first example, we will follow the original work of MacKay \cite{mackay_criterion_1989} and 
study the motion of a charged particle with momentum
$p$ in the potential of two longitudinal electrostatic waves. Without loss of
generality, we can choose one of the waves to have phase velocity zero and the
other to have phase velocity one and normalize the mass of the particle to 
one \cite{escande_renormalization_1981}.  This model, of the form
\eqref{eq:perturbPend}, can be described by a nonautonomous Hamiltonian system,
\begin{equation}
  \label{eq:twowaveHam}
  H(q,p,t) = \tfrac12 p^2 - \mu \cos(2\pi q) - \mu \nu \cos(2\pi k (q-t)),\qquad \mu,\nu\in\R,\, k\in\N.
\end{equation}
Thus the first wave has amplitude $\mu$, and the second has amplitude $\mu\nu$
and wavenumber $k$. By considering $q \mod 1$ and $k(q-t) \mod 1$, the phase space can be thought
of as $\T \times \R \times \T$. 

In fact, using the Poincare-Liouville one-form \eqref{eq:PoincareLiouville}
we see that the system is Cartan-Arnol'd, that is, the Hamiltonian vector field
$v$ associated to \eqref{eq:twowaveHam} also satisfies $\iota_v\alpha =
0$. Thus we can use the sign of $K$, \eqref{eq:KCond}, in our computations.

The Hamiltonian \eqref{eq:twowaveHam} is perhaps the simplest nonintegrable,
nonlinear Hamiltonian with two resonant terms.
There are two obvious integrable limits:
\begin{enumerate}
\item $\mu = 0$ yields the Hamiltonian for the motion of a free particle,
  \begin{equation}
    \label{eq:freePartHam}
    H(q,p,t) = \tfrac12 p^2.
  \end{equation}
\item $\nu = 0$ yields the Hamiltonian for a standard pendulum,
  \begin{equation}
    \label{eq:pendulumHam}
    H(q,p,t) = \tfrac12 p^2 - \mu\cos(2\pi q).
  \end{equation}
\end{enumerate}
If either $\mu\ll 1, \nu\sim O(1)$ or $\nu\ll 1, \mu\sim O(1)$  the dynamics is 
``nearly integrable'', and thus, by KAM theory, there will be a large set of invariant two-tori.

\begin{figure}[ht]
  \centering
  \begin{subfigure}[b]{0.48\textwidth}
     \includegraphics[width = \linewidth]{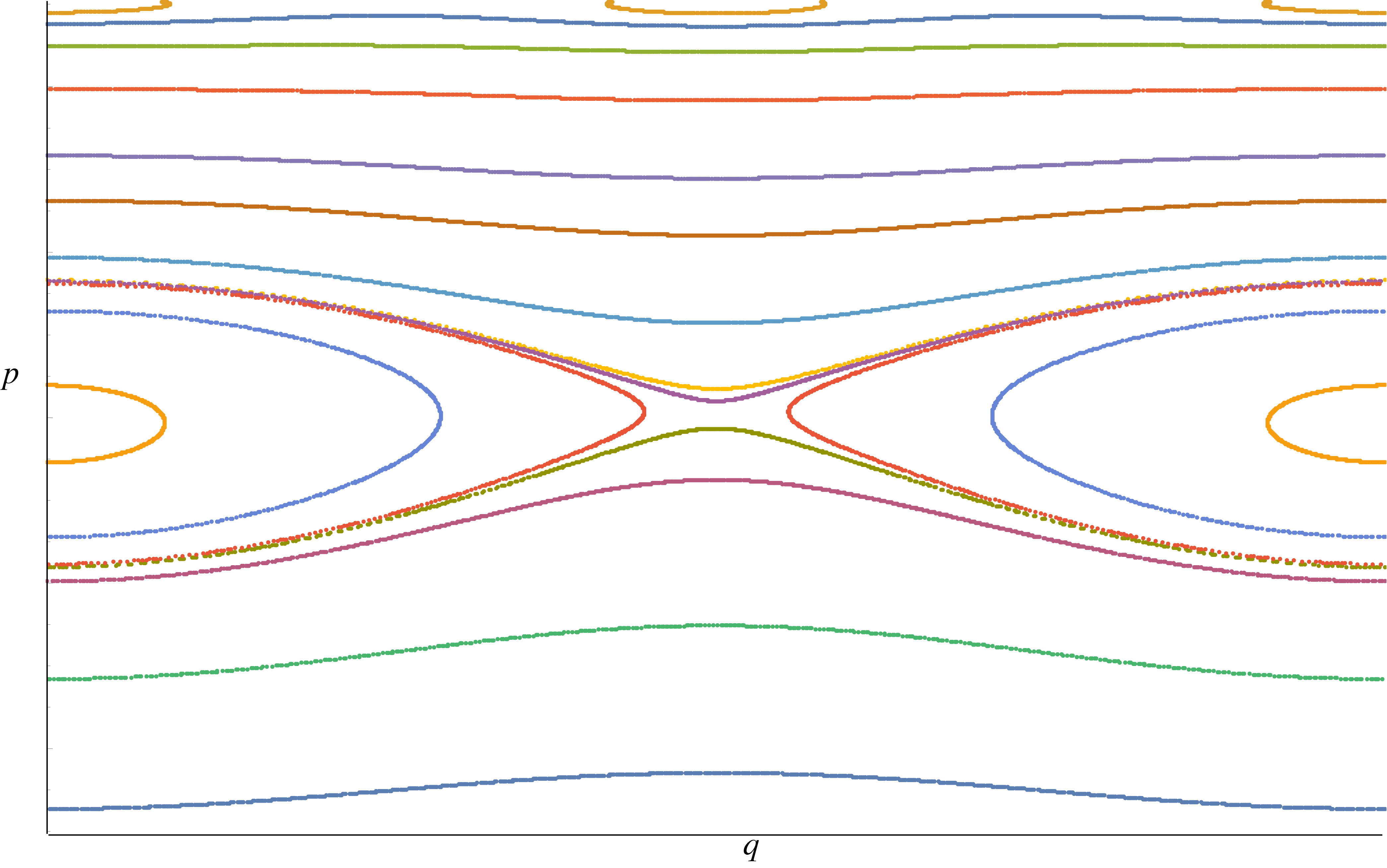}
     \caption{$\mu = 0.0075$}
    \label{fig:plot0075}
  \end{subfigure}
  \hfill
  \begin{subfigure}[b]{0.48\textwidth}
    \includegraphics[width = \linewidth]{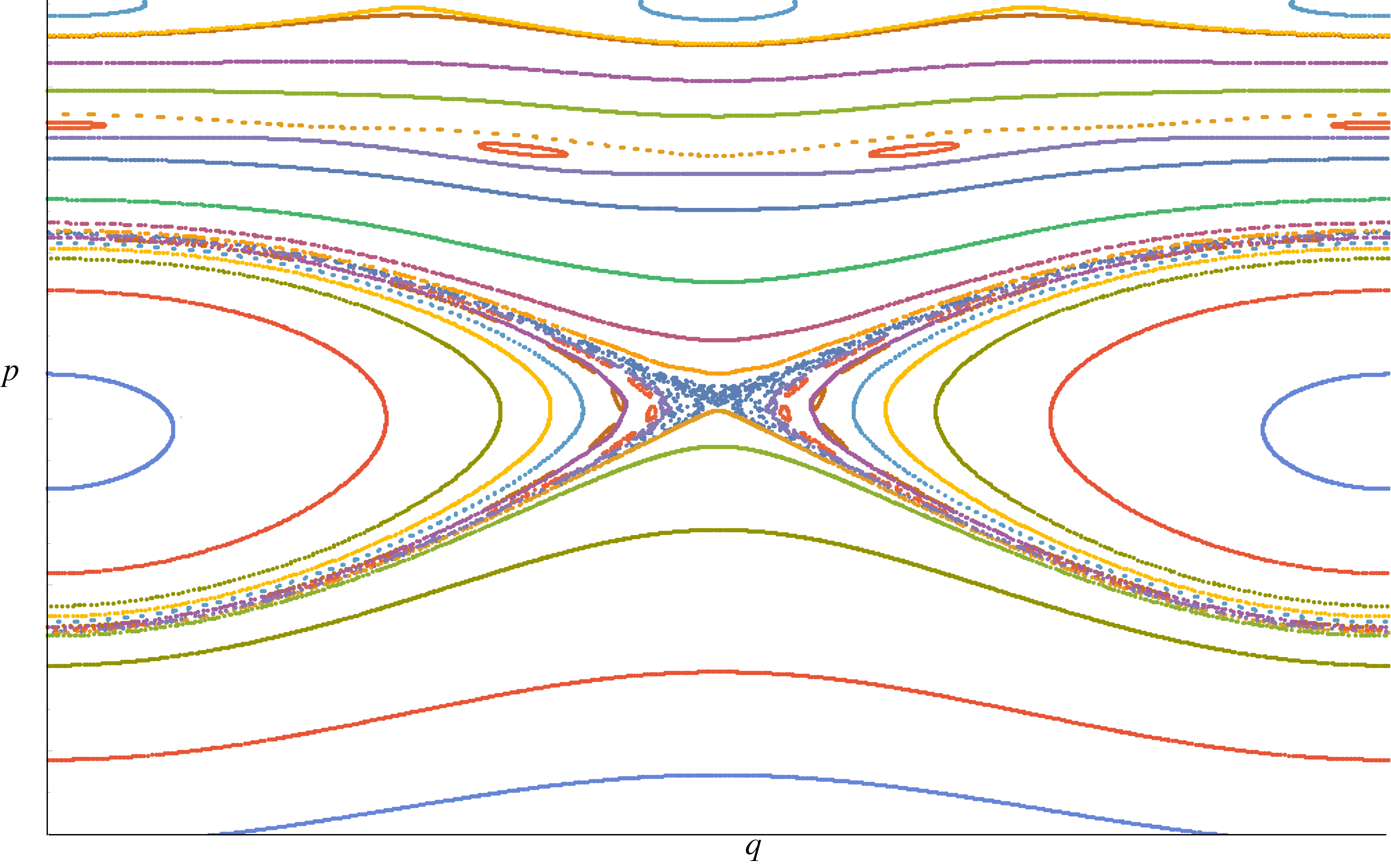}
     \caption{$\mu = 0.015$}
    \label{fig:plot015}
  \end{subfigure}
  \begin{subfigure}[b]{0.48\textwidth}
    \includegraphics[width = \linewidth]{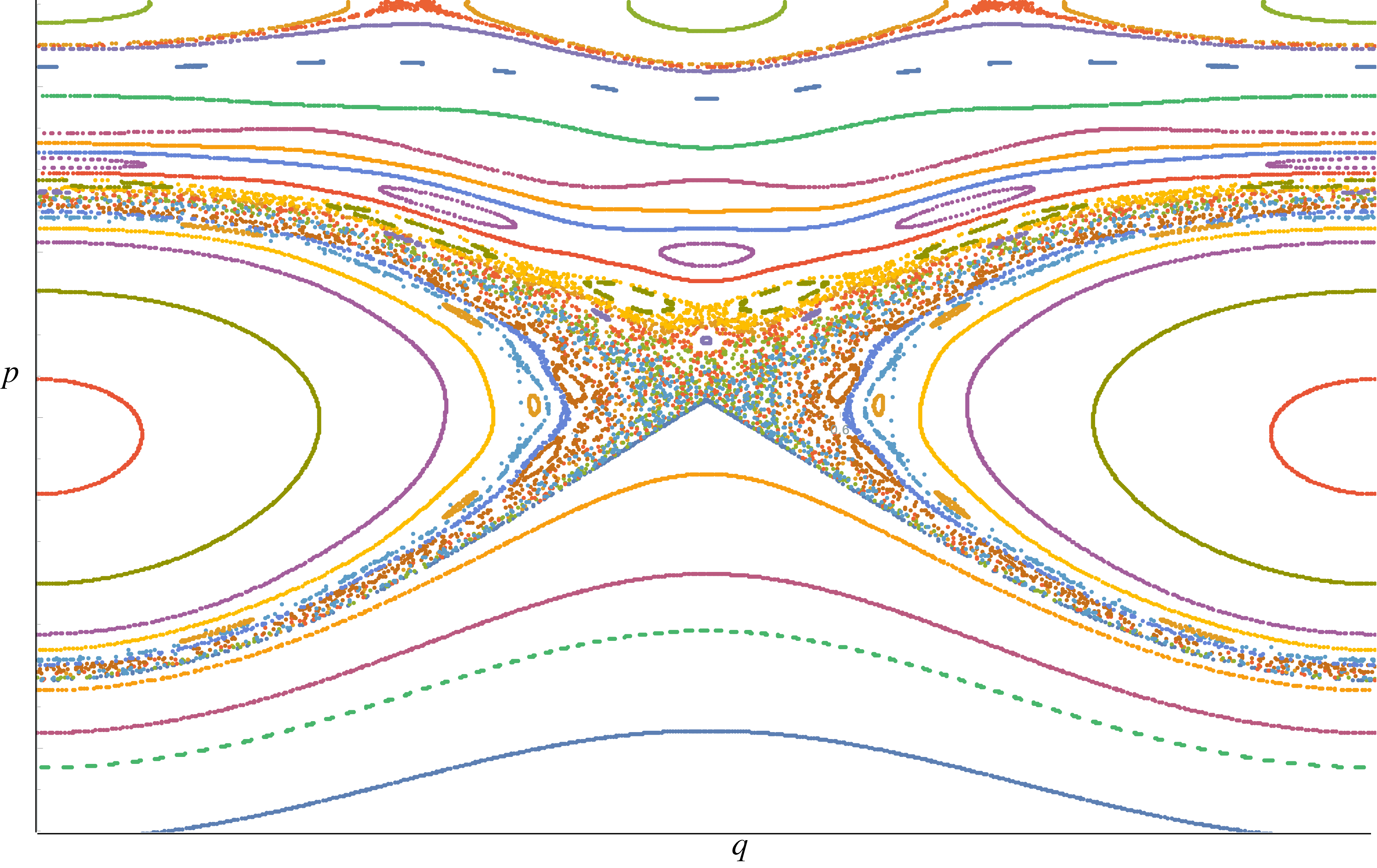}
    \caption{$\mu = 0.0225$}
    \label{fig:plot0225}
  \end{subfigure}
  \hfill
  \begin{subfigure}[b]{0.48\textwidth}
    \includegraphics[width = \linewidth]{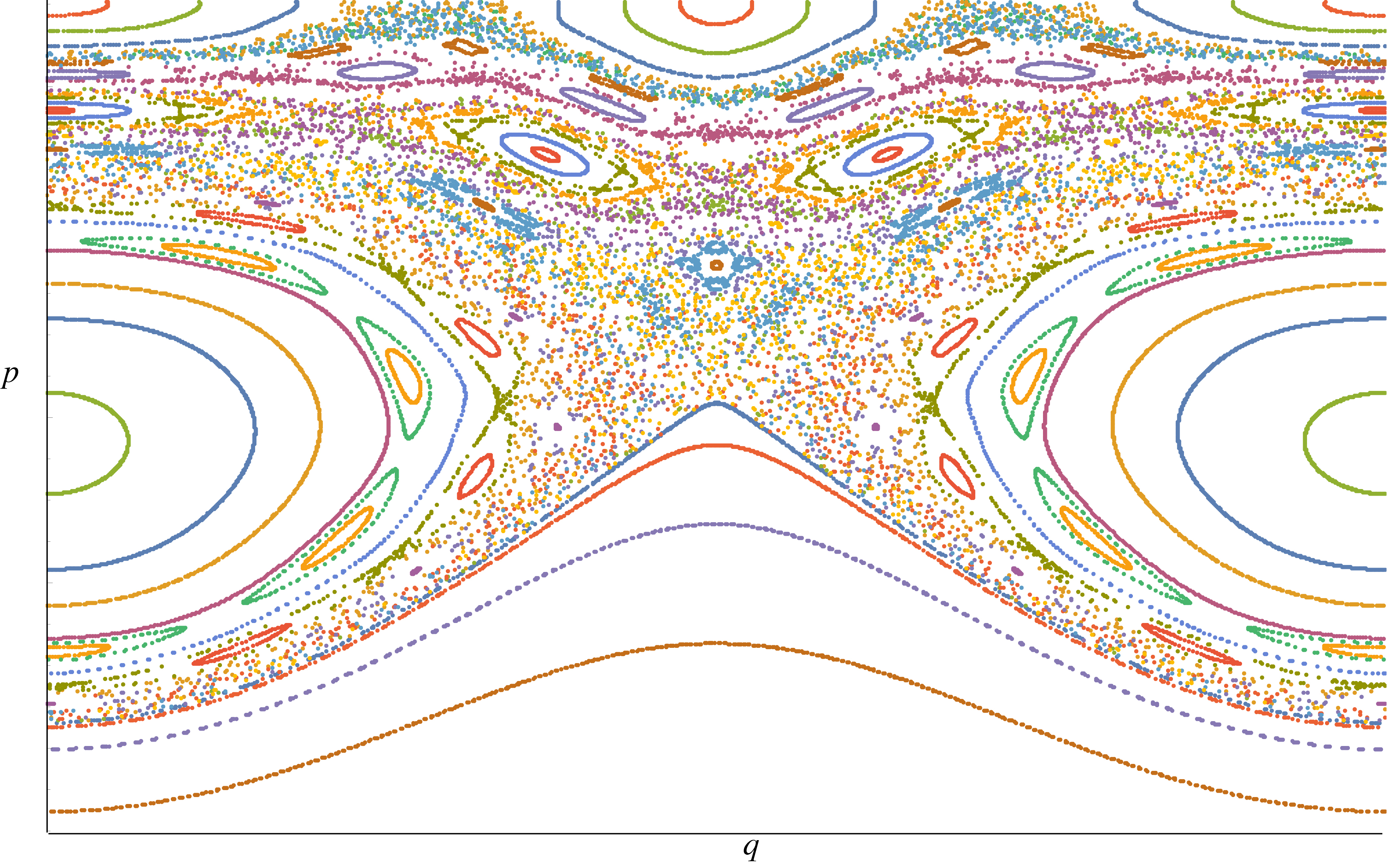}
    \caption{$\mu = 0.03$}
     \label{fig:plot03}
  \end{subfigure}
  \caption{Poincar\'{e} sections for \eqref{eq:twowaveHam} at $t = 0$ for four values of $\mu$ when $\nu=k=1$.
  The domain of the figures is $(0,1) \times (-\frac12, \tfrac12)$.}
  \label{fig:Msmall}
\end{figure}

Following \cite{mackay_criterion_1989}, we will study \eqref{eq:twowaveHam} with $\nu = k = 1$. 
Poincar\'{e} sections for four values of $\mu$, are shown in
\cref{fig:Msmall}. Note that since the system is periodic in $t$, the section
can be taken at $t = 0 \mod 1$.
Near $\mu = 0$, most orbits lie on rotational invariant tori, but
as $\mu$ increases, the section exhibits two primary
elliptic fixed points near $(0,0)$ and $(0,1)$ and two hyperbolic fixed
points near $(\tfrac12,0),(\tfrac12,1)$. The figures only show the lower wave, near $p=0$.
Of course, these fixed points correspond to elliptic and hyperbolic
periodic orbits of the nonautonomous flow. The tori encircling
the primary elliptic orbits are \textit{librational}; they
correspond to particles trapped in one of the two electrostatic waves.

By $\mu = 0.015$, chaos is clearly evident near $(\tfrac12,0)$ in \cref{fig:Msmall}(b).  
As $\mu$ increases a secondary island
appears near $p = \tfrac12$ that generates its own chaotic region.
As $\mu$ increases further, more resonant rotational tori break up, creating island
chains and chaos near their hyperbolic orbits. This is clearly evident in
\cref{fig:Msmall}(c). By $\mu = 0.03$, in \cref{fig:Msmall}(d), the chaotic region around 
the secondary island merges with that of the original island: Chirikov's resonance overlap.

Converse KAM theory was applied to this model in
\cite{mackay_criterion_1989}. MacKay looked for orbits that
do not lie on rotational invariant tori and
gave several related ``criteria'' for the existence of these tori. 
The most relevant to numerical application relies on finding conjugate points.
When $H$ has positive definite kinetic energy, this criterion is equivalent
to the foliation condition of \cite{mackay_finding_2018} upon choosing a 
foliation $\cF$ that is the union of vertical lines:
\[
	\cF = \bigcup_{(q_0,t_0) \in \T^2} \{(p,q_0,t_0): p \in \R\}.
\]
Of course this can also be thought of as the foliation generated by the function $J = \tfrac12 p^2$.

Below, we will extend MacKay's results by choosing several different foliations. 
Our aim is to to capture some of the librational tori around the main resonances 
at $p = 0$ and $1$, as well as those in the driven resonance near $p = \tfrac12$.

\subsection{Zaslavsky's Q-Flows: A Beltrami Example}
In this section, to contrast the results of the two-wave model, 
we will study Zaslavsky's so-called \textit{Q-flows} \cite{zaslavsky_dynamical_1991}.
While these dynamical systems are similar to $1\tfrac12$ degree-of-freedom Hamiltonian
systems, they do not have an obvious global Poincar\'e section. We will argue that they are more
properly thought of as vector fields obtained from a Cartan-Arnol'd one form.

\begin{figure}[ht]
  \centering
        \centerline{
         \includegraphics[width=0.45\textwidth]{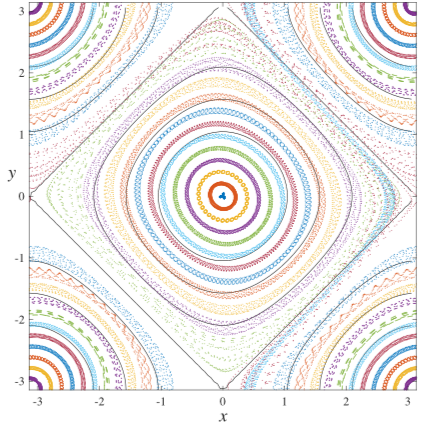}
         \hskip 0.25in
         \includegraphics[width=0.45\textwidth]{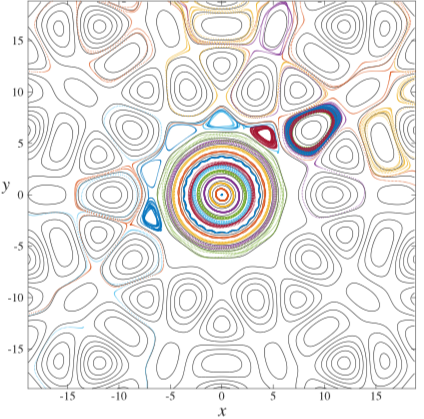}}
 \caption{Projections of the flow of \eqref{eq:qflow} onto the $(x,y)$ plane for $q=4$ (left),
 $q=5$ (right) and $\varepsilon = 0.15$ for 25 initial conditions along the line $y=x$ for $z=0$. 
 The contours of $\psi$ are shown as thin black curves. For $q=4$, periodicity is used in the $(x,y)$ plane.}
  \label{fig:qflows}
\end{figure}

We begin with a vector field 
$v_0 = \hat{z} \times \nabla \psi_q = (-\partial_y \psi_q,\partial_x \psi_q, 0)$, where the
stream function is
\begin{equation}
  \label{eq:qpotential}
  \psi_q (x,y) := \sum_{j=1}^q \cos\left(x \cos(2\pi j/q) + y \sin(2\pi j/q)\right).
\end{equation}
Note that $\psi_q$ is an invariant of the flow of $v_0$ since
\[
	v_0 \cdot \nabla \psi_q  = 0 .
\]
Thus the orbits lie on the surfaces $\psi_q = const$.
These stream functions have discrete, $q$-fold rotational symmetry on $\R^2$. For
$q = 1,2,3,4,6$, the contours of $\psi_q$ form a periodic lattice; otherwise they
have quasi-crystal symmetry, see \cref{fig:qflows}.

The Q-flows correspond to a 3D perturbation of $v_0$, with the vector field
\begin{equation}
  \label{eq:qflow}
  v = (\partial_y \psi_q + \varepsilon \sin(z), -\partial_x \psi_q + \varepsilon \cos(z), \psi_q)
\end{equation}
for each $\varepsilon \in \R$.
This vector field maintains the $q$-fold symmetry about the $z$-axis. 
The case $q=4$ is a special case of the ABC vector field of \cite{Dombre86}.

The vector field \eqref{eq:qflow} is incompressible and satisfies
the Beltrami condition:
\begin{equation} \label{eq:BeltramiCondition}
  \div v = 0,\qquad \curl v = \kappa v,
\end{equation}
with $\kappa = 1$, because $\nabla^2 \psi_q = -\psi_q$,

Zaslavsky \cite{zaslavsky_dynamical_1991} shows that the Q-flows can be formally
thought of as a Hamiltonian system if one uses $z$ as an effective time variable.
Unfortunately, this gives a singular vector field on the surfaces $\psi = 0$.
It is more appropriate to interpret $v$ as a Cartan-Arnol'd vector field
%
using the form
\begin{equation}
  \label{eq:QFlowOneForm}
  \alpha_{q} =  \left( \int^y \psi_q dy \right) dx - H(x,y,z) dz,
\end{equation}
with the effective Hamiltonian
\begin{equation}
  \label{eq:QFlowHam}
  H(x,y,z) = \psi_q(x,y) + \varepsilon ( y \sin(z) - x \cos(z)).
\end{equation} 
For this case the volume form \eqref{eq:VolumeForm} is $dy \wedge dx \wedge dz$ and
the vector field $v$ is the unique vector field satisfying
\[
	\iota_v (dy\wedge dx \wedge dz) = d\alpha_{q} = \psi_q dy \wedge dx - dH\wedge dz,
\]
Thus $\iota_v d\alpha_{q} = 0$, and the form $\alpha_{q}$ generates the vector field $v$, up to a scalar scaling.

We will study in more detail the cases $q=4$ and $q=5$ as representatives of spatially
periodic and quasi-periodic flows, respectively.
When $\varepsilon = 0$ and $q = 4$, the $(x,y)$
phase portrait contains a square lattice of elliptic equilibria surrounded by a lattice of hyperbolic equilibria.
Note that $\psi_4(x,y) = 0$ on the hyperbolic points and their heteroclinic orbits so $z(t)$ is constant. 
For the case $q = 5$, the contours of $\psi_5$ are a quasi-lattice 
(the contours in \cref{fig:qflows}(b)), but
there is still an elliptic orbit at the origin surrounded by a set of invariant tori.
When $\varepsilon$ grows from zero, an increasing degree of chaos in the near the separatrices of 
these lattices is observed. Projections of the dynamics for $\varepsilon = 0.15$ onto the $(x,y)$ plane 
are shown in \cref{fig:qflows}.

When chaotic orbits of the Q-flows are unbounded they typically exhibit anomalous diffusion---that is,
orbits travel along the channels associated with heteroclinic orbits, escaping
from any bounded domain. The mean square drift grows with some non-integer power of $t$. 
This phenomenon is the motivation
for the work in \cite{zaslavsky_dynamical_1991}. Here, we are
more concerned with the breakup of tori, changes in topology, and birth of chaos
as $\varepsilon$ increases. As we discuss next, to apply the converse KAM
method we must start by constructing adequate foliations to
detect these phenomenon.

\section{Foliations}\label{sec:Foliations}
The primary aim of this section is to describe the various foliations
to numerically test Thm.~\ref{thm:CKAMVolume} and determine the nonexistence of
invariant tori transverse to these foliations.

\subsection{Two-Wave Model}\label{sec:WaveFoliations}
For the two-wave model \eqref{eq:twowaveHam}, five foliations are
considered, each of which is generated from the gradient flow of an integrable
Hamiltonian.

\begin{figure}[ht]
  \centering
  \begin{subfigure}{0.48\linewidth}
    \includegraphics[width = \linewidth]{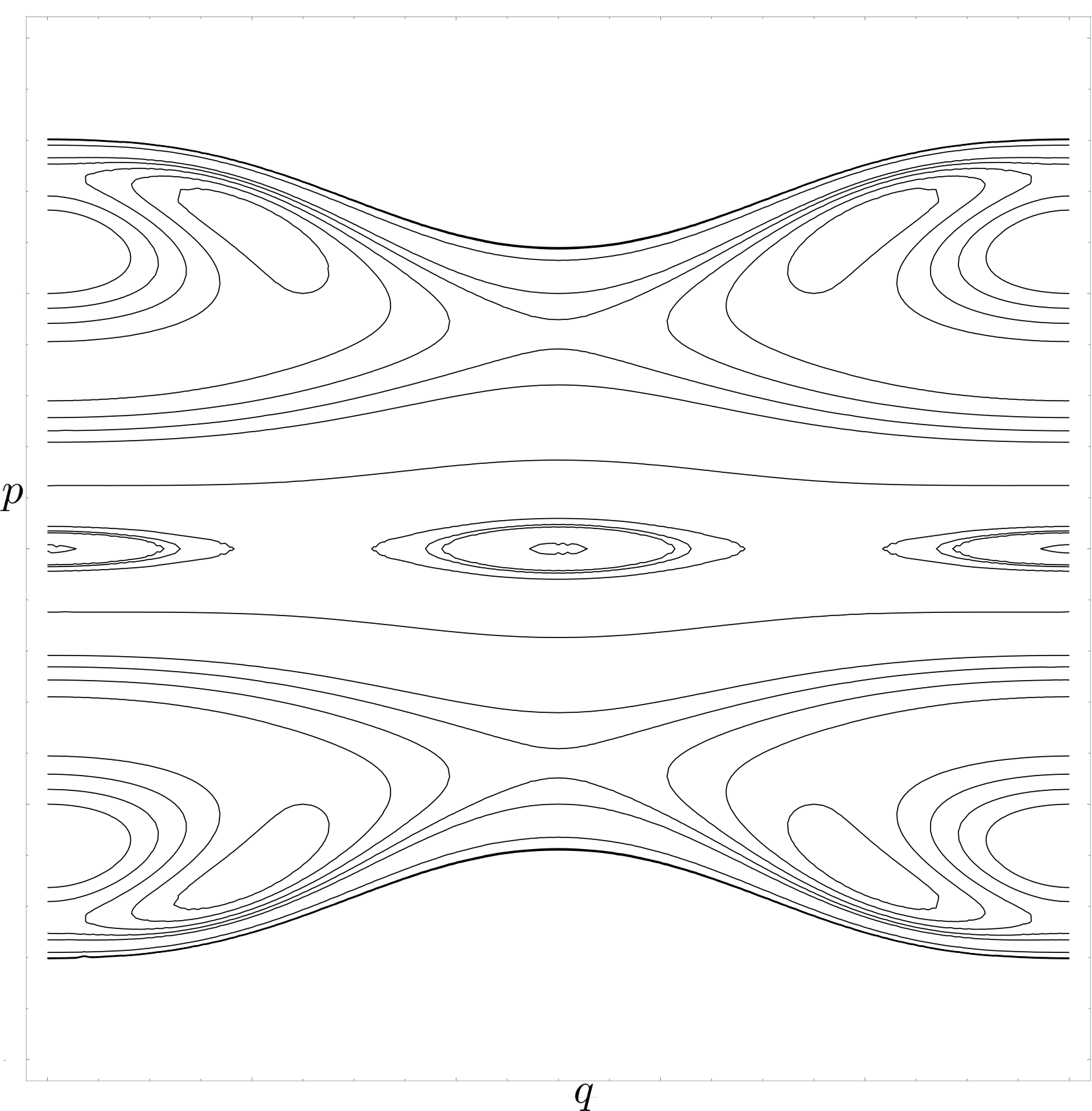}
  \end{subfigure}
  \begin{subfigure}{0.48\linewidth}
    \includegraphics[width = \linewidth]{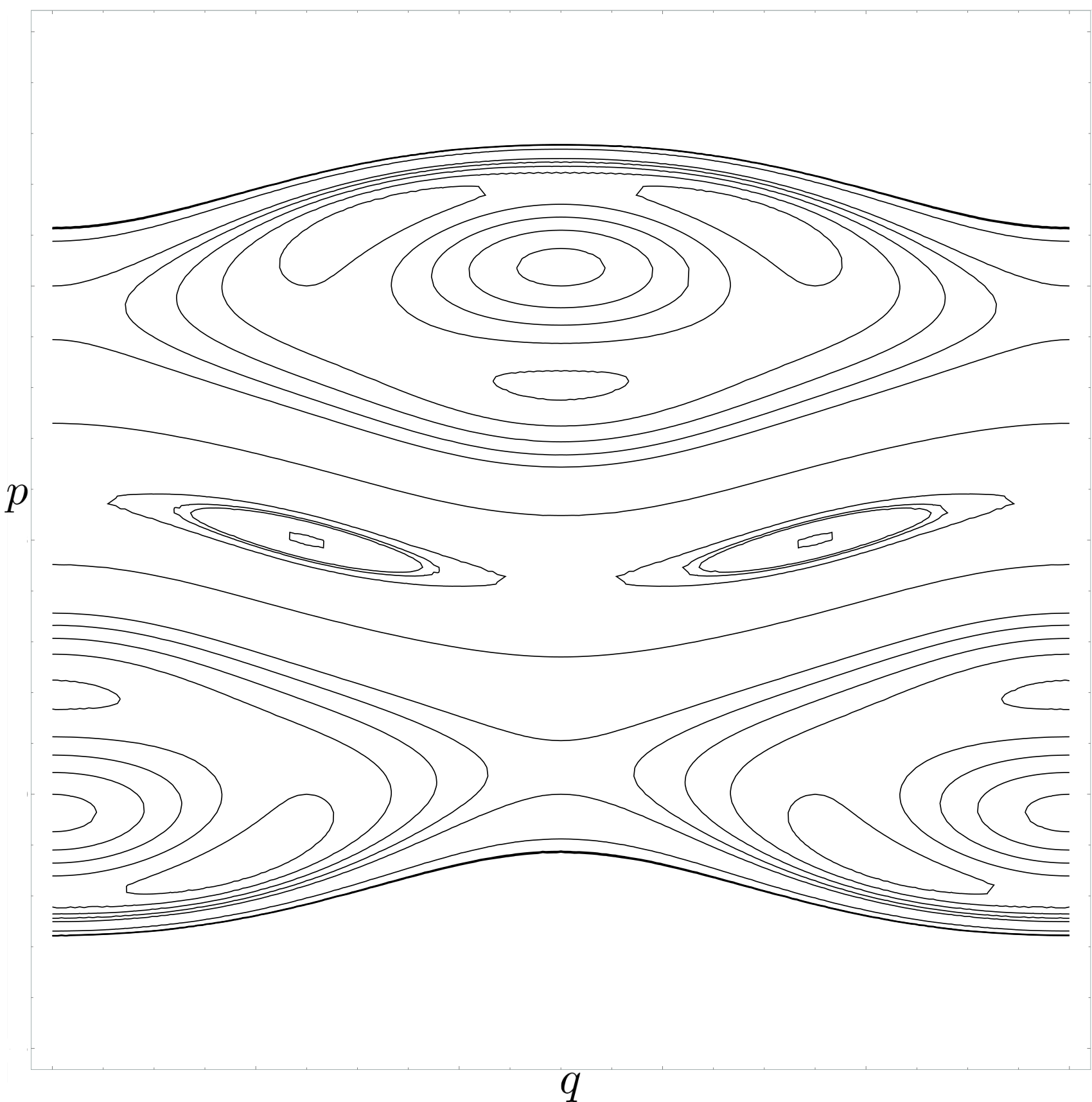}
  \end{subfigure}
  \caption{Level sets for the foliation $\cF^{s2}_{0.02}$ at time $t=0$ (left)
    and $t=\tfrac12$ (right). The domain of the figures is $(0,1)\times(-0.5,1.5)$.}
  \label{fig:s2foliationLevelSets}
\end{figure}

\begin{enumerate}
\item \emph{r-foliation}: Generated from the gradient flow of the free particle
  Hamiltonian \eqref{eq:freePartHam}. It is designed to capture so called
  rotational tori. Specifically, it is the foliation given by
\begin{equation}\label{eq:rFoliation}
  \cF^r = \bigcup_{q_0,t_0} \{(p, q_0, t_0) | p \in \R\}.
\end{equation}

\item \emph{l-foliation}: Generated from the gradient flow of the harmonic
  oscillator Hamiltonian $H = \tfrac12(p^2+q^2)$. It is designed to capture the
  so called librational tori around the elliptic fixed point $(0,0)$ of the
  Poincar\'{e} section.
  This foliation is given by rays of constant angle from the origin in the $(q,p)$ plane:
\begin{equation}\label{eq:lFoliation}
 	\cF^l = \bigcup_{\theta,t_0} \{(r\cos\theta, r\sin\theta, t_0) | r \in \R^+\}
\end{equation}

Note that this foliation is singular along $(q,p) = (0,0)$. Consequently, the
curve $(q,p) = (0,0)$ needs to be removed from the phase space. Any orbit that
intersects this removed curve is removed from consideration of the converse KAM
theorem.
  
\item \emph{p-foliation}: Generated 
  from the gradient flow of the pendulum Hamiltonian \eqref{eq:pendulumHam}. 
  It is designed to capture both rotational
  and librational tori. This foliation $\cF^p_\mu$ changes with the value of
  $\mu$, and is thus technically a family of foliations. 
  This foliation reduces to $\cF^{r}$ when $\mu = 0$.

  As in the l-foliation, the p-foliation has singular points. These occur along
  $(q,p) = (0,0)$ and $(q,p) = (\tfrac12,0)$, and are elliptic and hyperbolic,
  respectively. Any orbit that passes though either of these singularities must
  be removed from consideration under the converse KAM theorem.
  
\item \emph{s1-foliation}: Generated by the gradient flow of a first order
  global invariant,
  \begin{equation}
    \label{eq:approxInv}
    J = \tilde{H}_\mu(q,p,t) = -\tfrac12 p^2 + \tfrac13 p^3 
  	- \mu \left[(p-1) \cos (2 \pi q) + \nu p \cos (2 \pi k (q-t))\right].
  \end{equation}
  Taking the orbits of the gradient flow of \eqref{eq:approxInv} for each value
  of $\mu$ gives the family of foliations $\mathcal{F}^{s}_\mu$. The invariant
  is computed using a method called `global removal of resonances'
  from \cite{lichtenberg_regular_1992}. The details of this method and its
  application to the two-wave model are given in Appendix \ref{sec:Invariants}.

  The foliations $\cF^{s1}_\mu$ are designed to capture both rotational and
  librational tori around $(0,0)$, as well as librational tori about the
  elliptic fixed point $(0,1)$, and track the movement of the two elliptic fixed
  points starting at $(0,0)$ and $(0,1)$ as $\mu$ varies.

\item \emph{s2-foliation}: Generated by the gradient flow of a second order
  global invariant, \eqref{eq:invariant2}. It is designed to capture all
  features of the s1-foliation as well as librational tori arising from the
  $1:2$ resonance. The family of foliations is denoted by $\cF^{s2}_\mu$. The
  invariant \eqref{eq:invariant2} is computed using a higher order method of the
  global removal of resonances from \cite{lichtenberg_regular_1992}. This higher
  order method was first realized in \cite{mcnamara_superconvergent_1978} and
  its details and application to the two-wave model are given in
  Appendix \ref{sec:Invariants}.

  The level sets of \eqref{eq:invariant2} at the slice $t = 0$ are given in
  Figure \ref{fig:s2foliationLevelSets}. 
\end{enumerate}

\subsection{Q-Flow}
For the Q-flows defined in \eqref{eq:qflow} we will consider only two
foliations:
\begin{enumerate}
\item \emph{l-foliation}: Generated from the gradient flow of the Hamiltonian
  for the harmonic oscillator $H(x,y,z) = \tfrac12(x^2+y^2)$. The foliation is
  expected to capture only the primary librational tori that continue from those
  that encircle the elliptic orbit at $(0,0)$ for $\varepsilon = 0$.
\item \emph{$\psi$-foliation}: Generated from gradient flow of the approximate
  invariant $\psi_q(x,y)$ \eqref{eq:qpotential}. This foliation
  should capture the librational tori around any elliptic point in the periodic
  or quasi-periodic lattice of $\psi_q$.
\end{enumerate}

\section{Results}\label{sec:results}

In this section we will compute the converse KAM condition for the two-wave and Q-flow models.
In each case we will use several foliations, as discussed in \S\ref{sec:Foliations}.
The main results are shown below in \cref{fig:rfoliation,fig:lfoliation,fig:pfoliation,fig:s1foliation,fig:s2foliation,fig:q4plots,fig:q5plots}.

\subsection{Numerical Implementation}
For each model and each foliation the condition of Thm.~\ref{thm:CKAMVolume}
will be tested numerically. Since both models can be formulated as Cartan-Arnol'd systems,
we will use the techniques in \S\ref{sec:specialCases}. Thus the
converse KAM condition is satisfied for the two-wave or Q-flow models only if,
\begin{equation}\label{eq:KqfAnddw}
	K_{tw}(t) := d\alpha(\xi_t,\eta_t),\quad \mbox{or} \quad 
	K_{qf}(t) := d\alpha_{q}(\xi_t,\eta_t),
\end{equation}
changes sign, respectively.

To compute the sign change we choose a proposed point $u_0 \in M$
to check for nonexistence of an invariant surface, and 
numerically find the flow $\phi_t(u_0)$.
Since each foliation is determined by a function $J$, we choose
\begin{equation}
  \eta_t = \grad J(\phi_t(u_0)),\qquad \xi_0 = \eta_0 = \grad J(u_0).
\end{equation} 
so that  $\eta_t$ is in the tangent space to the foliation at each point
and $\xi_0$ begins parallel to $\eta_0$.
The vector $\xi_t$ is then obtained by integration of the linearized equations
 \eqref{eq:linearlization} about the orbit $\phi_t(u_0)$.  

As discussed in \S\ref{sec:specialCases}, $K(t)$ changing sign is necessary
but not sufficient to check the hypothesis of Thm.~\ref{thm:MacKay} as it
equally detects $\varphi_t$ crossing $\pi$. 
To remedy this, we will only check
the change of sign of $K(t)$ when $\eta_t \cdot \xi_t < 0$, using the
Euclidean dot product, as this guarantees that $\theta_t < \pi/2$.

Both models and their tangent orbits, are integrated
numerically using the $5/4$-Runge-Kutta scheme of Tsitouras
\cite{tsitouras_rungekutta_2011}; specifically we use the algorithm \texttt{Tsit5()} of
the \texttt{DifferentialEquations.jl} library in Julia. A callback function is
used to determine sign changes in $K(t)$ when $\eta_t\cdot \xi_t < 0$, and
linear interpolation to find a more accurate value of $t_c$, the time for which
$K(t_c) = 0$.

\subsection{Two-Wave Model}\label{sec:DWResults}

We computed the converse KAM condition for the two-wave model 
for the five different foliations $\cF^r$, $\cF^l$, $\cF^p$, $\cF^{s1}$ and $\cF^{s2}$
defined in \S\ref{sec:WaveFoliations}. 
We typically chose a grid of $500$ values for $\mu \in [0,0.03]$, 
and of $500$ initial points on the line $(0,p_0,0)$ for $p_0 \in [0,1]$.
For each initial point we integrated the orbit until either the converse KAM condition was satisfied at
some time $t_c$ or $t$ reached time $150$.

\begin{figure}[ht]
  \centering
  \includegraphics[width = 0.5\linewidth]{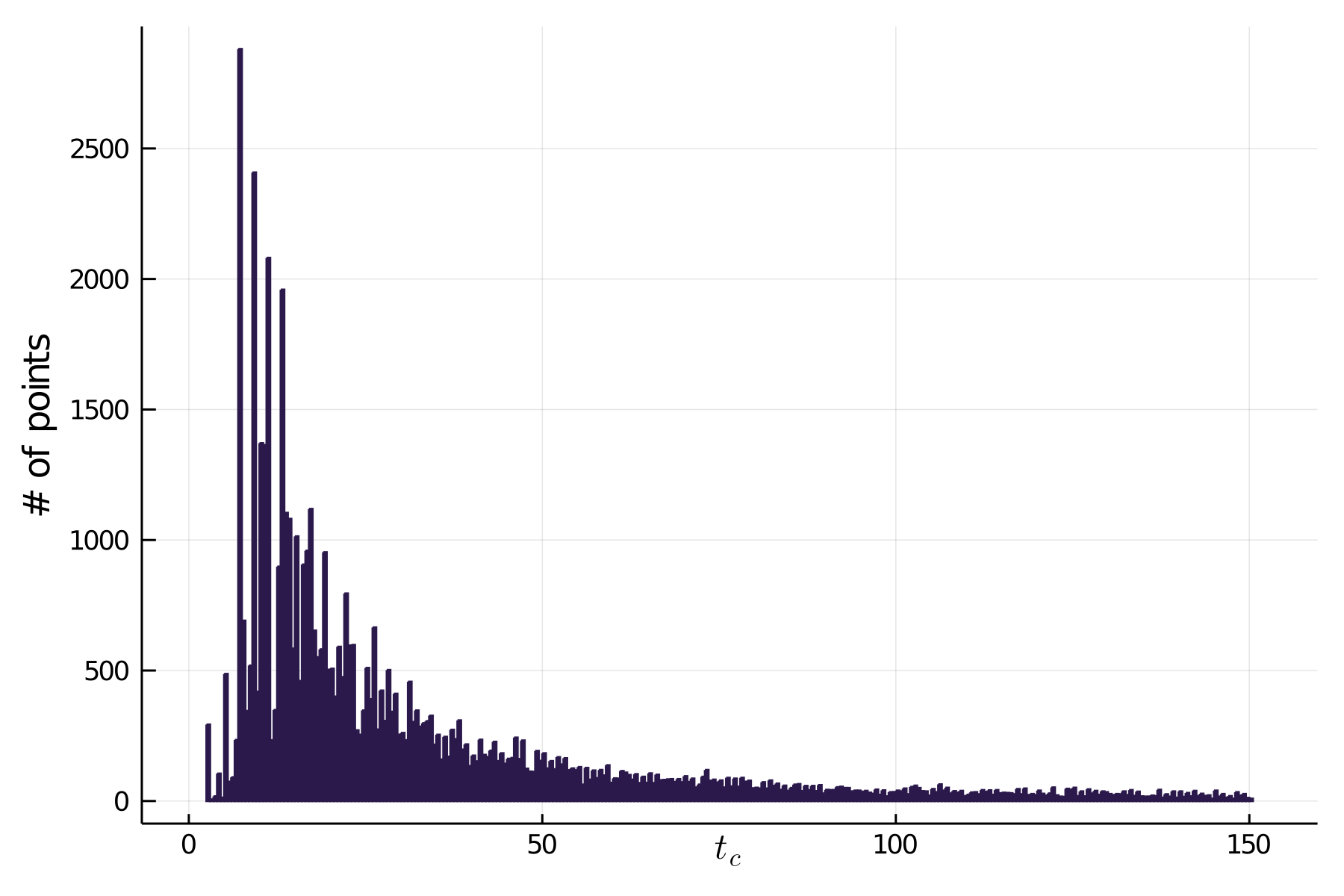}
  \caption{Histogram of the number orbits transverse to the s1-foliation
    at time $t_c$ for $250,000$ orbits of the two-wave model with $\mu \in [0,0.03]$.}
  \label{fig:convergencePts}
\end{figure}

However, first we address the choice of maximum integration time.
The choice $t_c \le 150$, seemed to be sufficient to capture most of the points
that satisfy the converse KAM condition. This is indicated by \cref{fig:convergencePts},
which shows a histogram of the number of points captured at
time $t_c$ for the s1-foliation over the entire $500 \times 500$ grid.
It is evident that by $t = 150$ the vast majority of points that would eventually
satisfy the converse KAM criterion are captured.
A similar convergence is also seen for the other foliations.

\begin{figure}[ht]  \centering
  \begin{subfigure}{0.49\linewidth}
    \includegraphics[width=\linewidth]{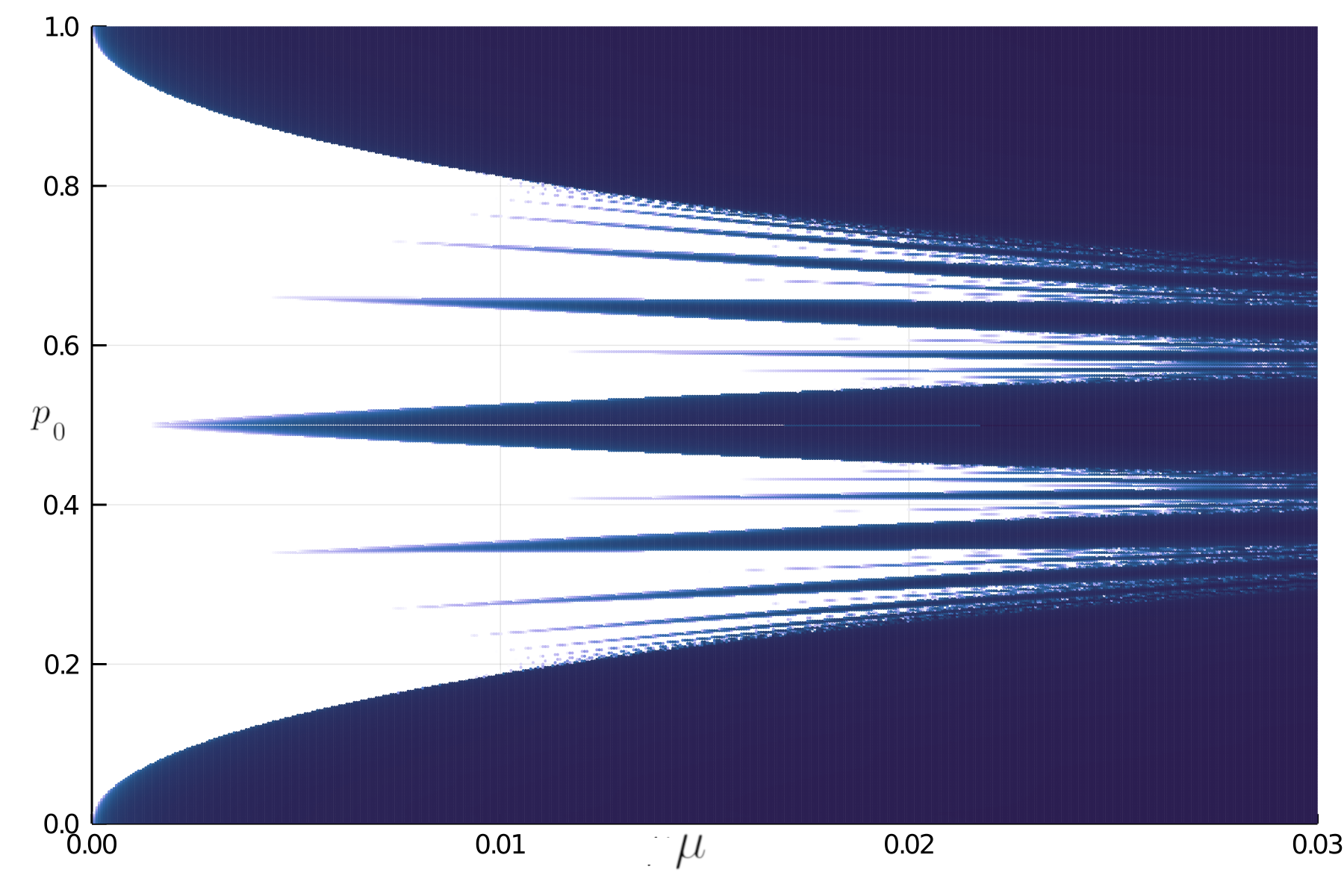}
  \end{subfigure}
  \begin{subfigure}{0.49\linewidth}
    \includegraphics[width=\linewidth]{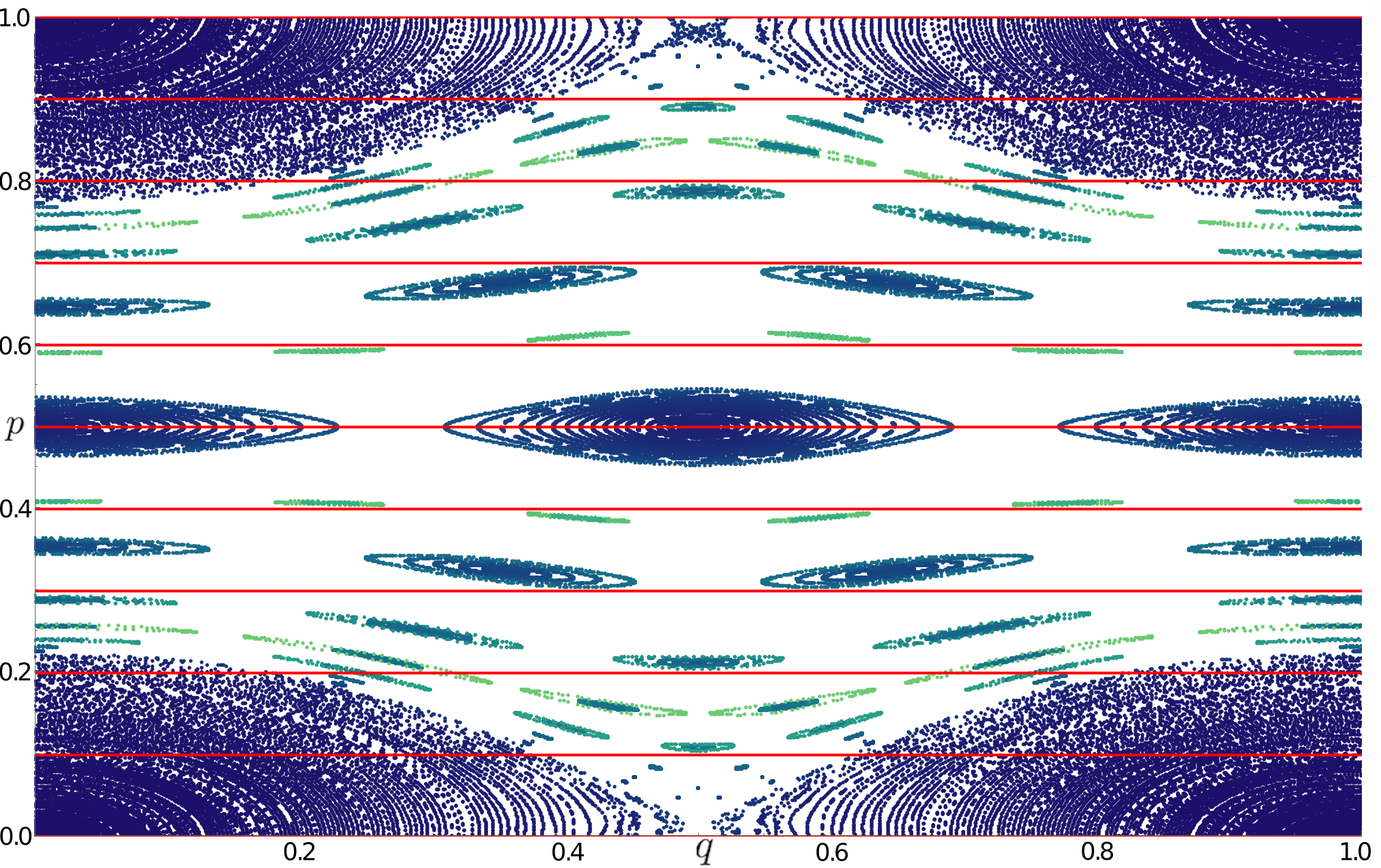}
  \end{subfigure}
 \caption{(r-foliation) (a) The left panel depicts the values of $(\mu,p_0)\in[0,0.03]\times[0,1]$ that
     pass the converse KAM condition within the time $t_c \leq 150$ for the r-foliation. 
     The points are
     colored according to the time taken to transversality violation, with a
     darker blue indicating a shorter time. 
     (b) The right panel shows a $t= 0$ Poincar\'e section
     of the orbits from (a) for $\mu = 0.015$. 
     The red curves are a cross-section of the level sets
     of the generating invariant for the foliation, namely $p^2$.}
 \label{fig:rfoliation}
 \end{figure}

\begin{figure}[ht]
  \centering
  \begin{subfigure}{0.49\linewidth}
    \includegraphics[width=\linewidth]{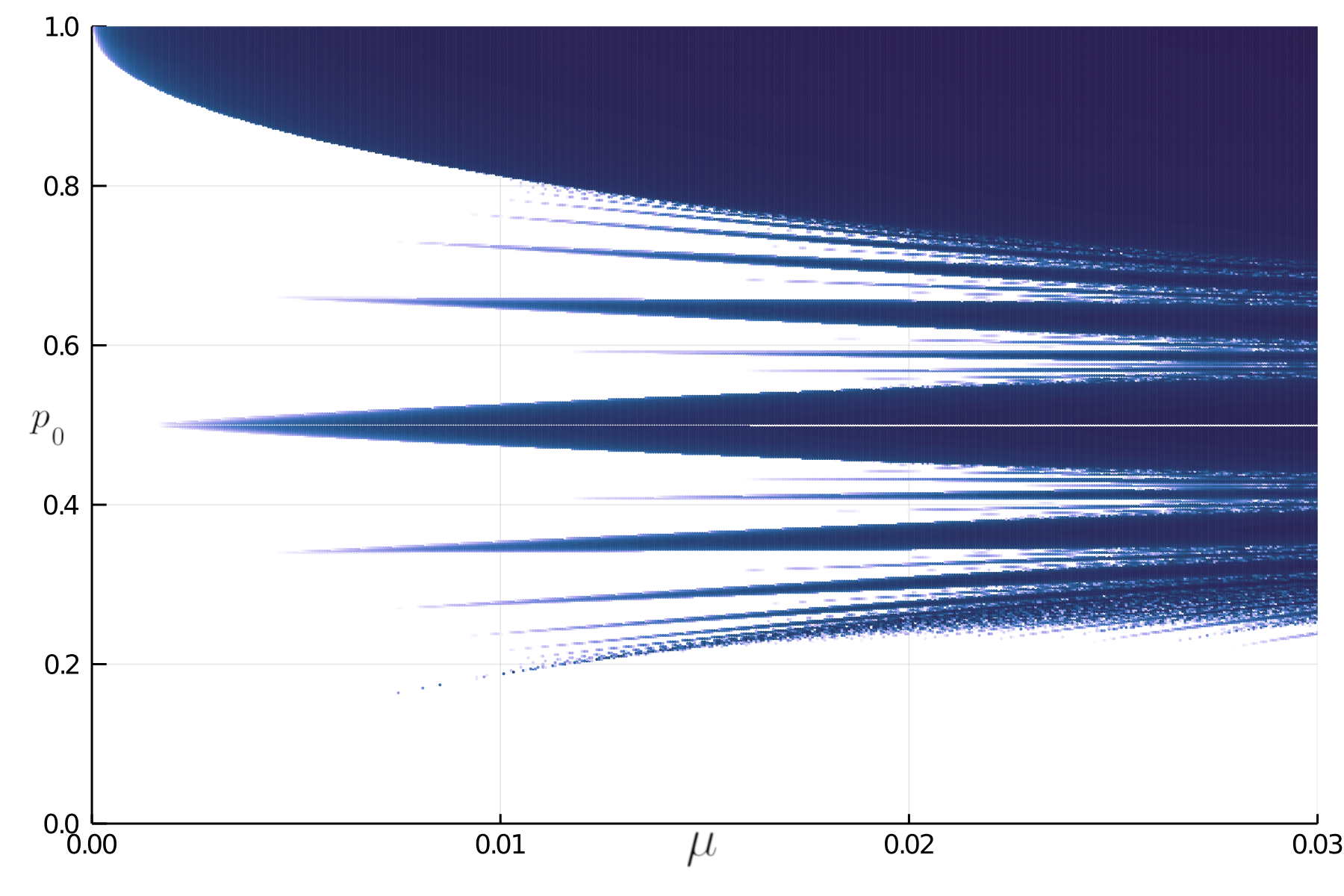}
  \end{subfigure}
  \begin{subfigure}{0.49\linewidth}
    \includegraphics[width=\linewidth]{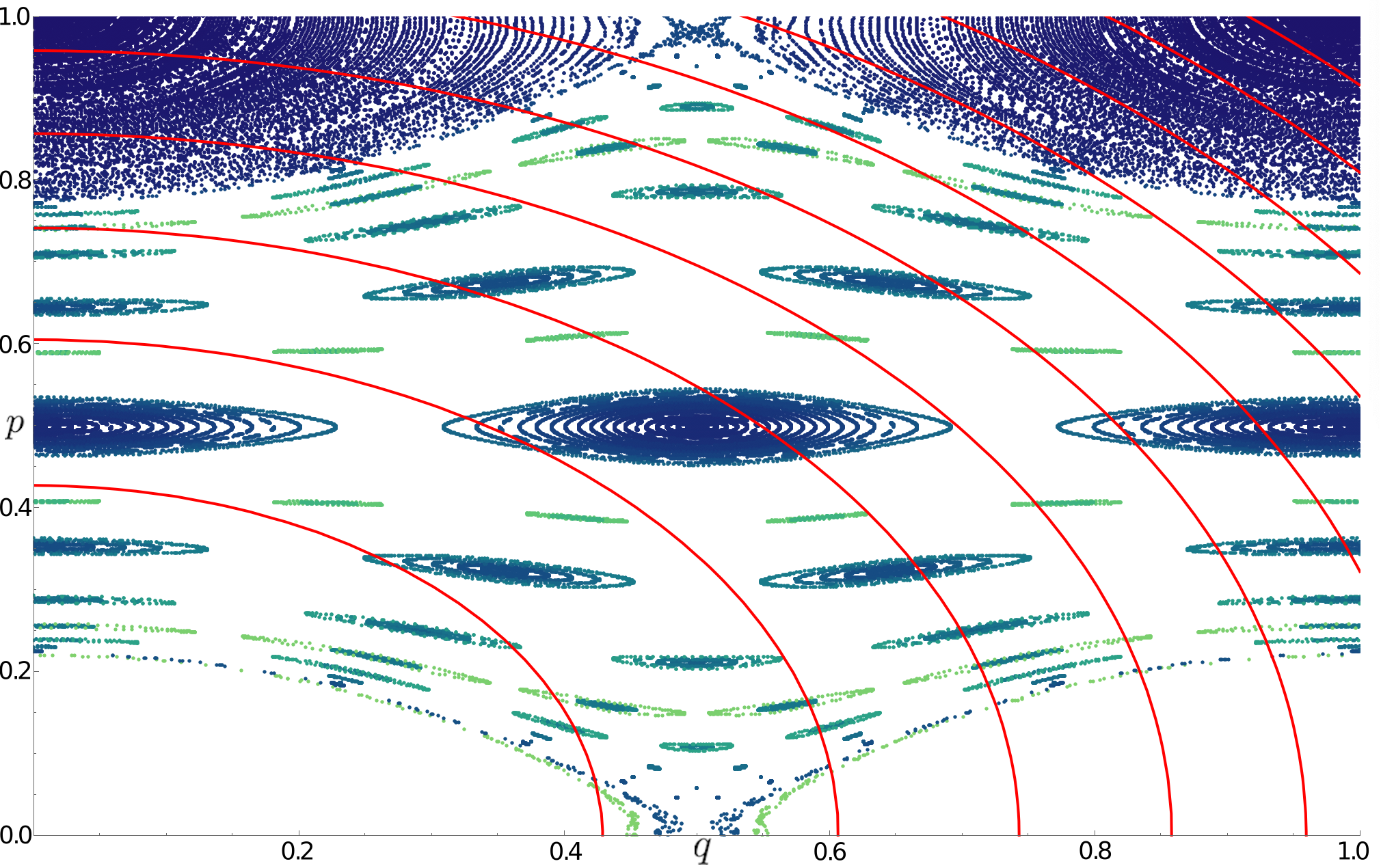}
  \end{subfigure}
 \caption{(l-foliation) Orbits satisfying the converse KAM condition for
 the l-foliation (details as in the caption to \cref{fig:rfoliation}).
 The red curves (right panel) correspond to the Hamiltonian $\tfrac12(p^2+q^2)$.}
 \label{fig:lfoliation}
 \end{figure}

\begin{figure}[ht]
  \centering
  \begin{subfigure}{0.49\linewidth}
    \includegraphics[width=\linewidth]{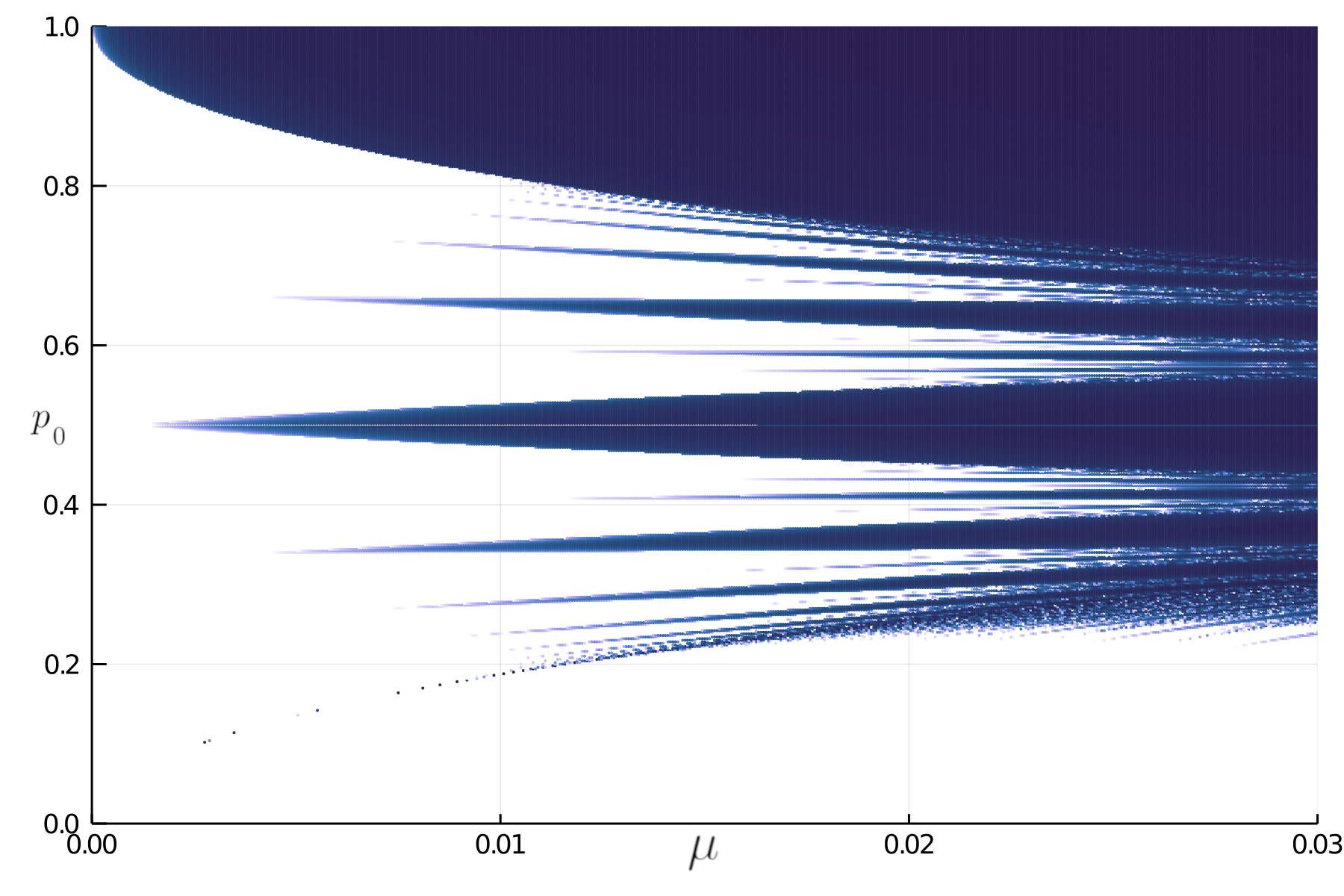}
   \end{subfigure}
  \begin{subfigure}{0.49\linewidth}
    \includegraphics[width=\linewidth]{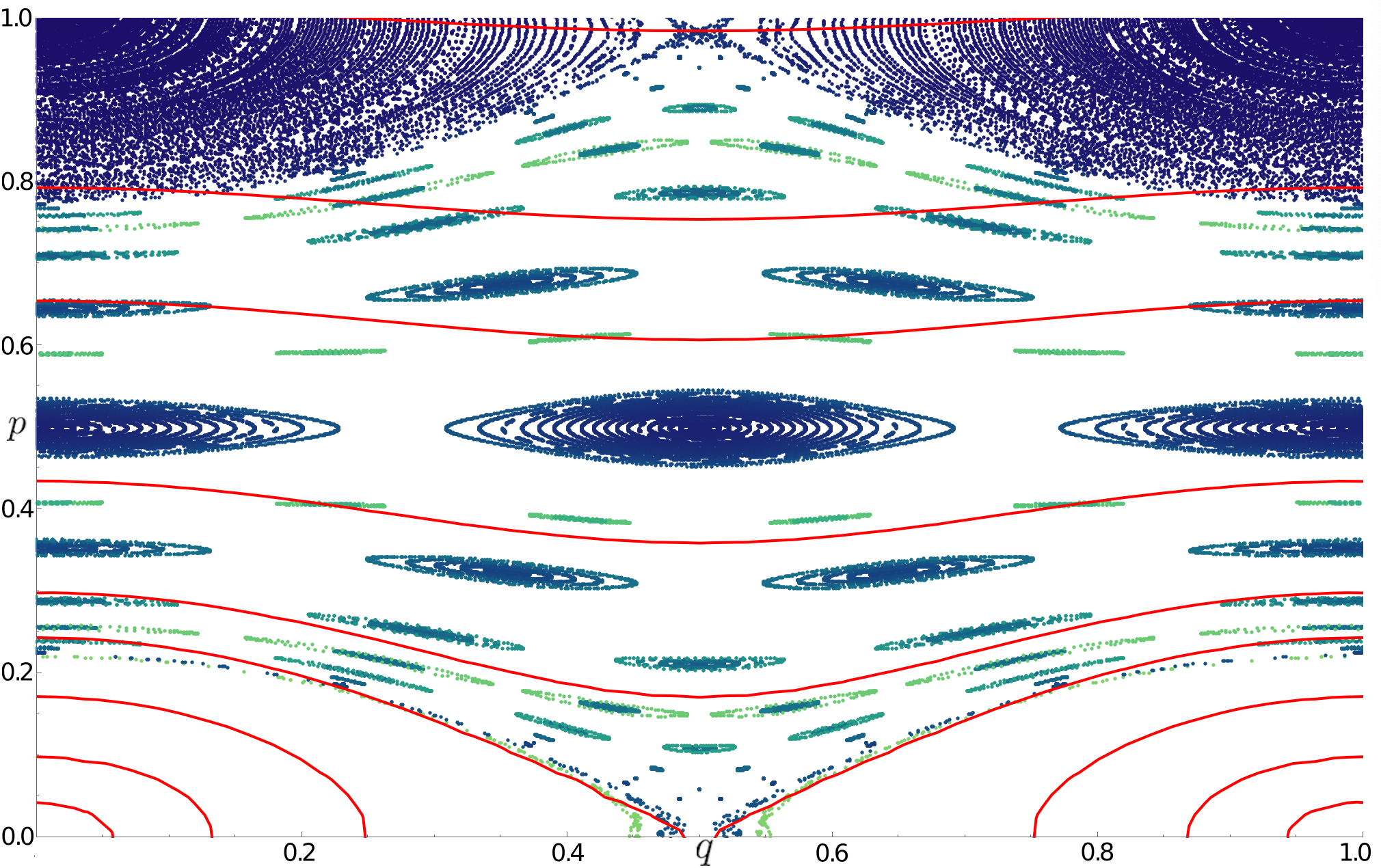}
  \end{subfigure}
 \caption{(p-foliation) Orbits satisfying the converse KAM condition for
 the p-foliation (details as in the caption to \cref{fig:rfoliation}).
 The red curves (right panel) correspond to the Hamiltonian \eqref{eq:pendulumHam}.}
 \label{fig:pfoliation}
\end{figure}

\begin{figure}[ht]
	\centering
	\begin{subfigure}{0.49\linewidth}
		\includegraphics[width=\linewidth]{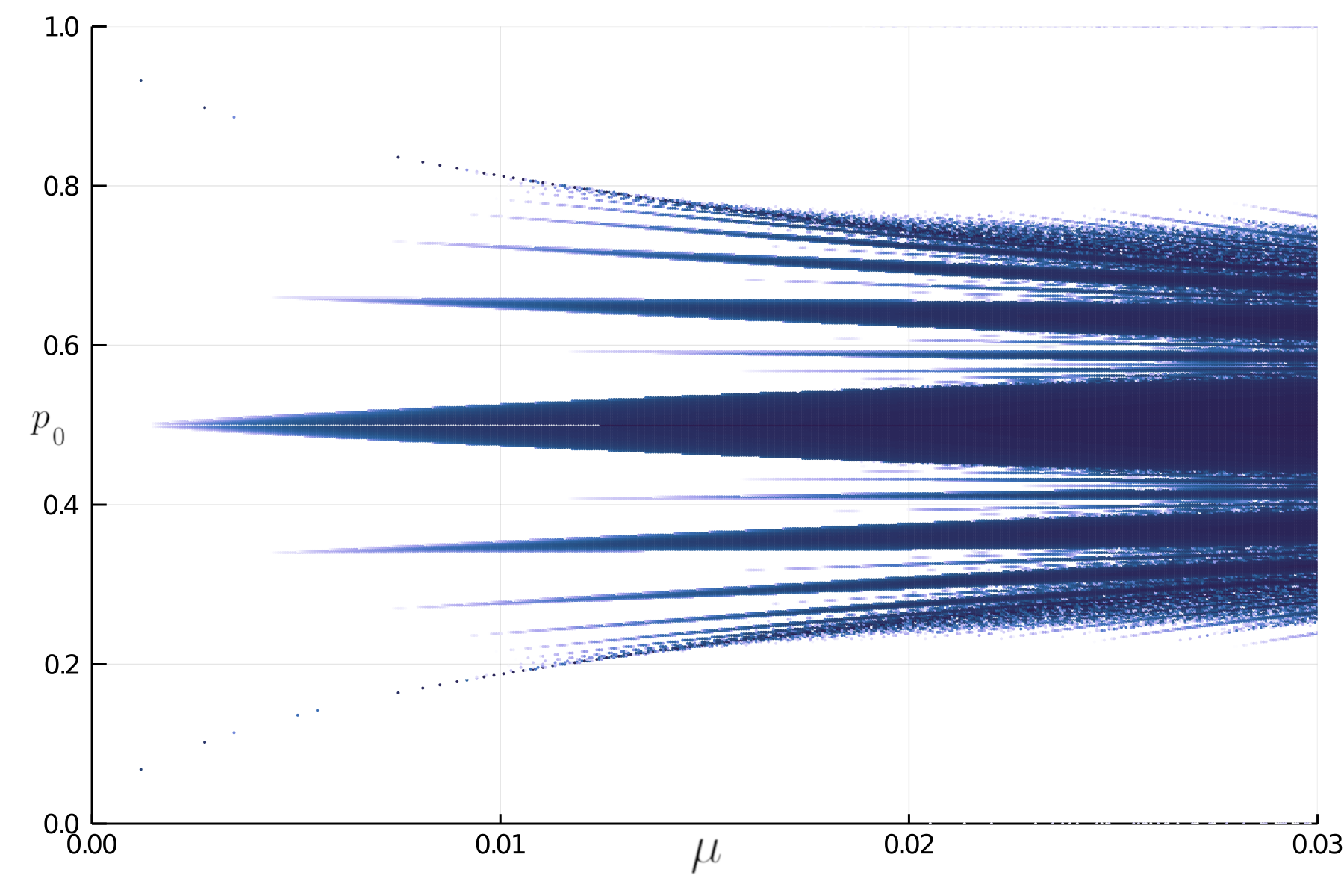}
	\end{subfigure}
	\begin{subfigure}{0.49\linewidth}
		\includegraphics[width=\linewidth]{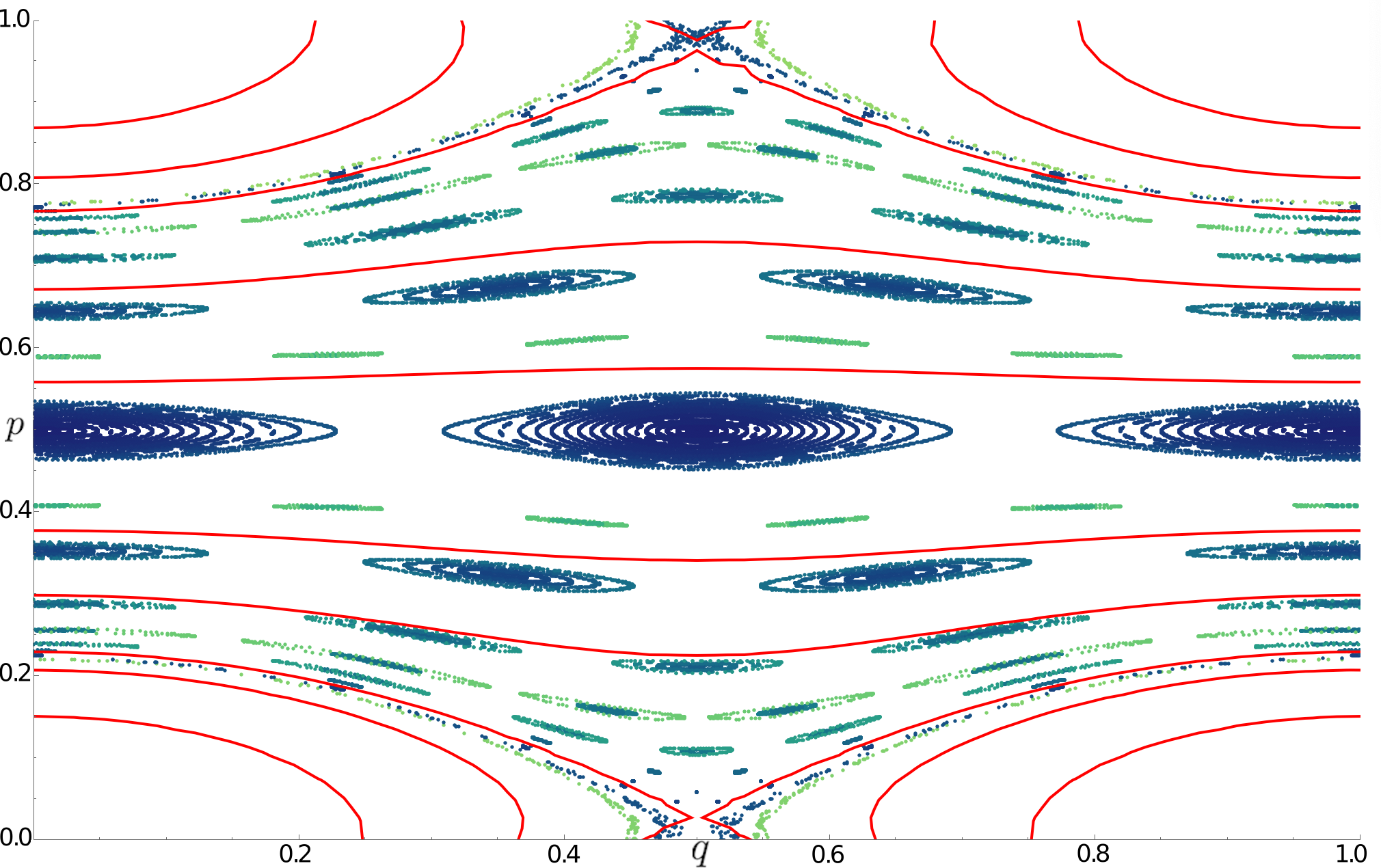}
	\end{subfigure}
	\caption{(s1-foliation) Orbits satisfying the converse KAM condition for
		the s1-foliation (details as in the caption to \cref{fig:rfoliation}).
		The red curves (right panel) correspond to the Hamiltonian \eqref{eq:approxInv}.}
	\label{fig:s1foliation}
\end{figure}

The results of the computations are shown in \cref{fig:rfoliation,fig:lfoliation,fig:pfoliation,fig:s1foliation,fig:s2foliation}
for the five foliations. As we will argue below, the results for the various foliations differ, especially
with respect to librational tori in various resonances of the dynamics.

The results for the r-foliation $\cF^r$, are shown in \cref{fig:rfoliation}.
In the left panel, which is equivalent to Fig.~4 of \cite{mackay_criterion_1989}, 
each initial condition that satisfies the converse KAM condition, i.e., 
that does not lie on an invariant surface transverse to the foliation, is colored blue, and the darker hues 
indicate a smaller value of $t_c$.

In the right panel of \cref{fig:rfoliation}, we show the $t=0$ Poincar\'e section
for each of the orbits that pass the converse KAM condition for $\mu = 0.015$.
The contours of the foliation are also overlaid on the panel.
Note that any initial conditions that lie within the 0:1 and 1:1 resonance islands do not lie on invariant surfaces
transverse to the r-foliation and thus they pass the converse KAM condition ``dependently''.
These correspond to librational tori around the two primary elliptic orbits that start
near $(q,p,t)=(0,0,0)$ and $(0,1,0)$.
This is expected; clearly librational tori, whose
Poincar\'e sections are circles about their respective elliptic points, fail to be
transverse to planes of constant $p$. This explains the two dark blue regions
in the left panel that grow in width as $2\sqrt{\mu}$, the half-width of these islands. 
Also seen in this panel are several 
tongues, the largest of which has a width that grows linearly in $\mu$ starting near $(\mu,p_0)\approx
(0.002,0.5)$. This corresponds to the driven resonance near $p = \tfrac12$ that was seen in
\cref{fig:Msmall} and is clearly visible in the right panel of \cref{fig:rfoliation} as well.

The l and p-foliations are shown in \cref{fig:lfoliation,fig:pfoliation}, respectively. For the range $0.2<p_0<1.0$, these
figures share many features of the r-foliation. However, the librational tori near $(q,p,t) =
(0,0,0)$ no longer satisfy the converse KAM condition, that is, the computations correctly
indicate that these orbits lie on invariant tori. 
In addition, small spikes close to the curve $p_0 = 2\sqrt{\mu}$ are now evident;
these represent the incursions of resonant islands into the edge of the main resonance.
The p-foliation differs subtly from the l-foliation in that there are slight
differences in some of the resonance tongues.

Though the s1-foliation, in \cref{fig:s1foliation}, is nearly identical to the p-foliation
for $0 < p < 0.8$, it also now identifies points on librational tori in the 1:1 resonance.
Thus the s1-foliation correctly identifies the regular regions about both the $0:1$ and $1:1$ resonances.

\begin{figure}
	\centering
	\begin{subfigure}{0.49\linewidth}
		\includegraphics[width=\linewidth]{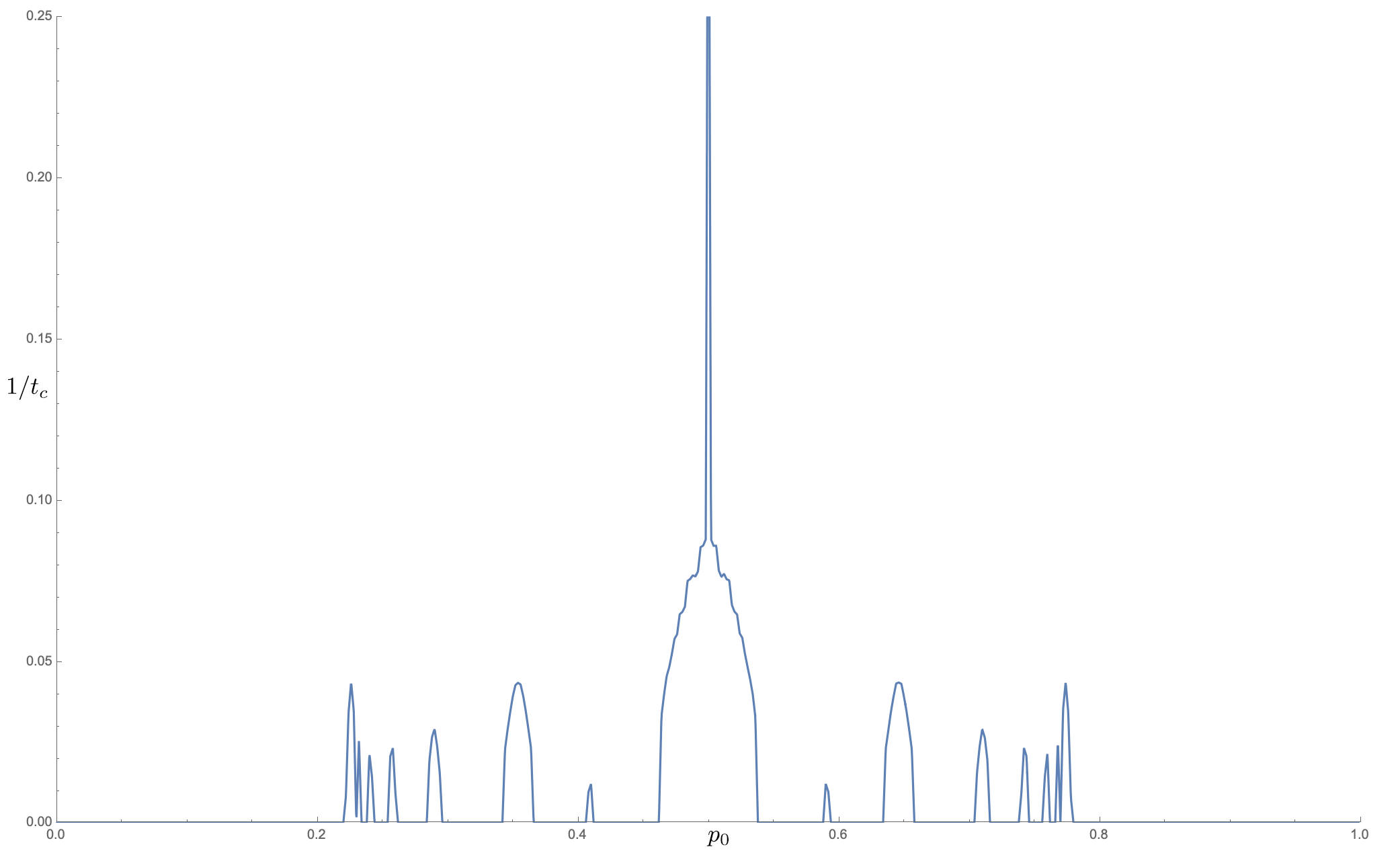}
	\end{subfigure}
	\begin{subfigure}{0.49\linewidth}
		\includegraphics[width=\linewidth]{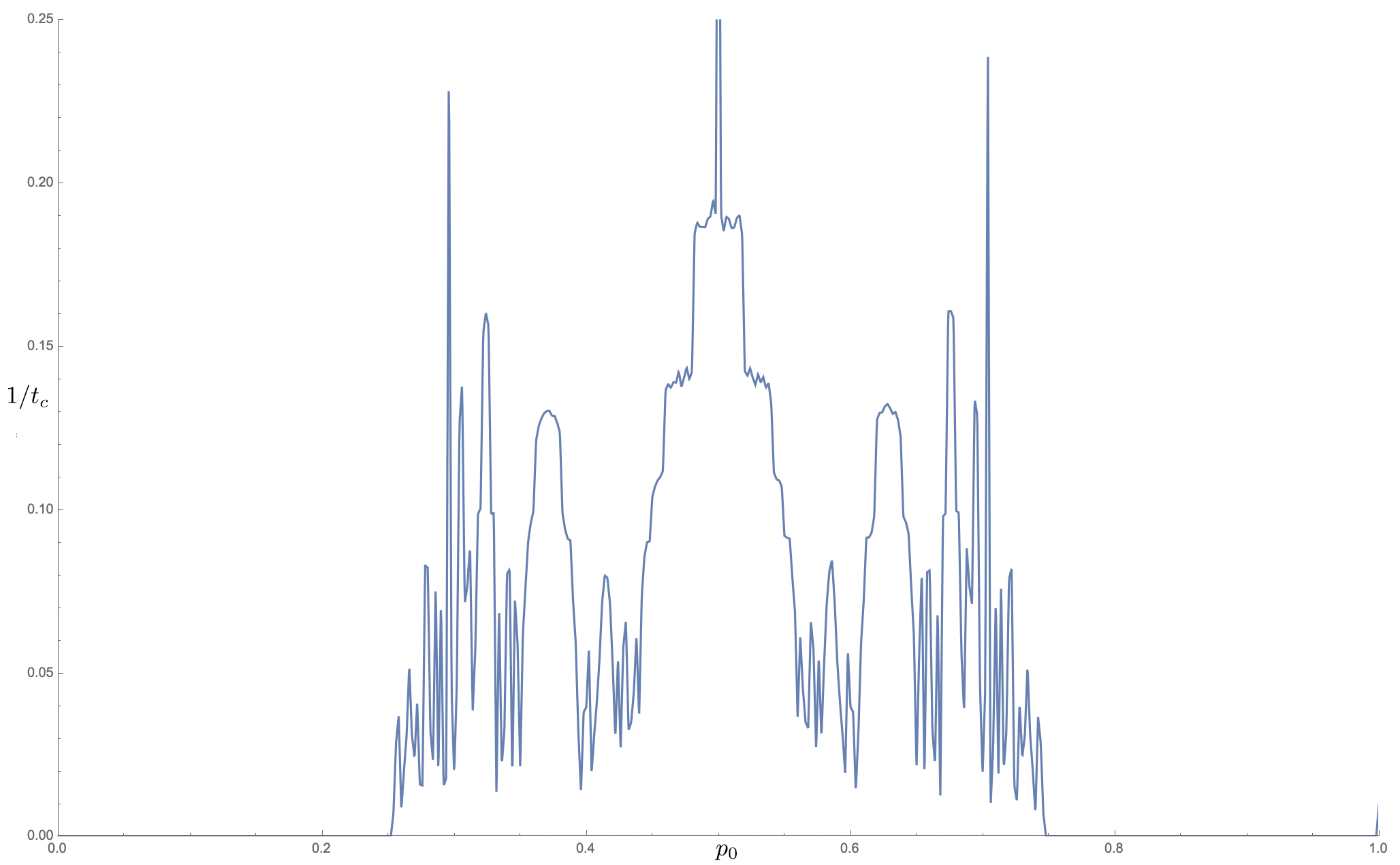}
	\end{subfigure}
	\caption{Plot of the inverse time $1/t_c$ taken to pass the converse KAM condition with respect to the s1-foliation for initial conditions $(0,p_0)$ for (left) $\mu = 0.015$ and (right) $\mu=0.03$. }
	\label{fig:freqPlots}
\end{figure}

\begin{figure}[ht]
  \centering
  \begin{subfigure}{0.49\linewidth}
    \includegraphics[width=\linewidth]{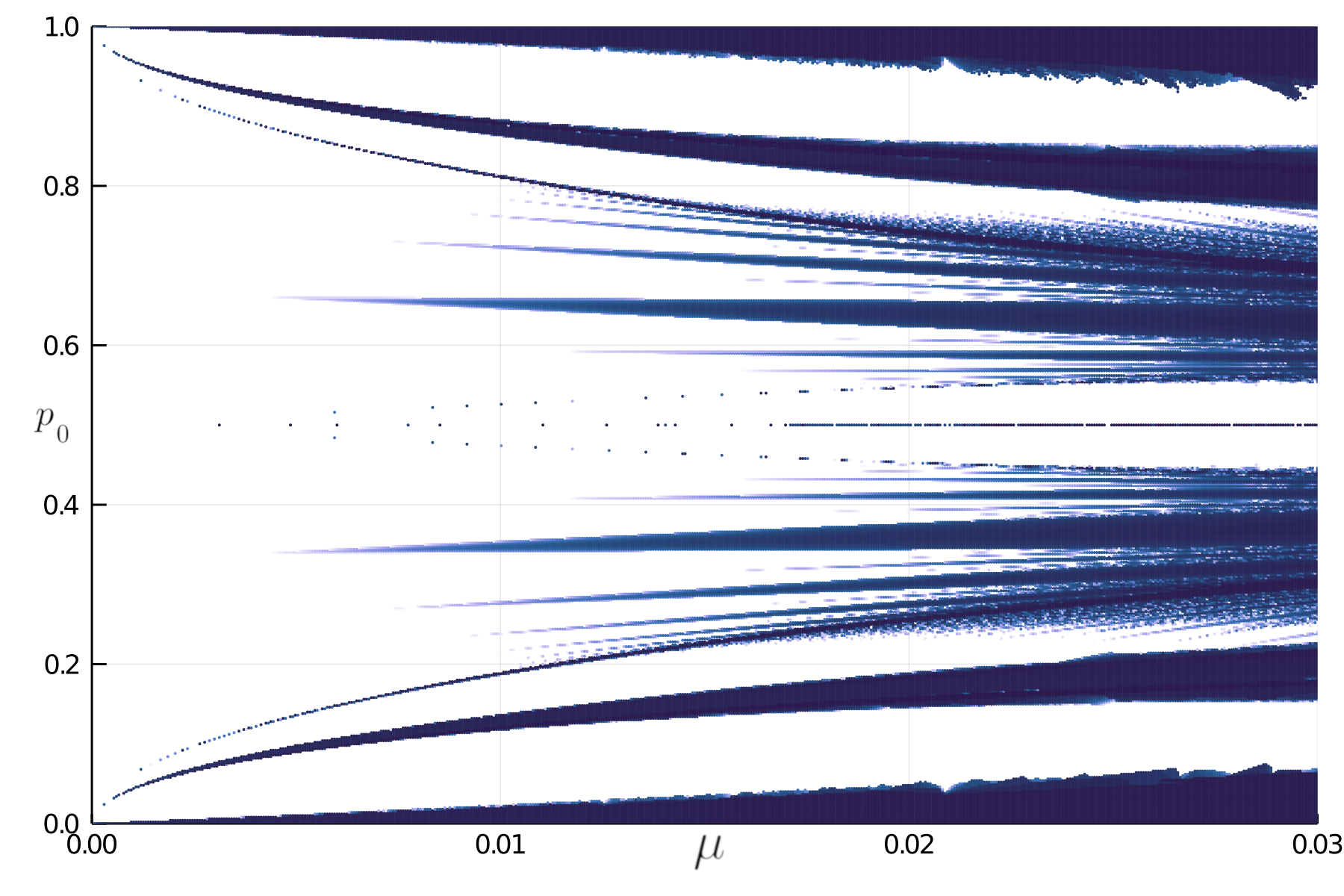}
  \end{subfigure}
  \begin{subfigure}{0.49\linewidth}
    \includegraphics[width=\linewidth]{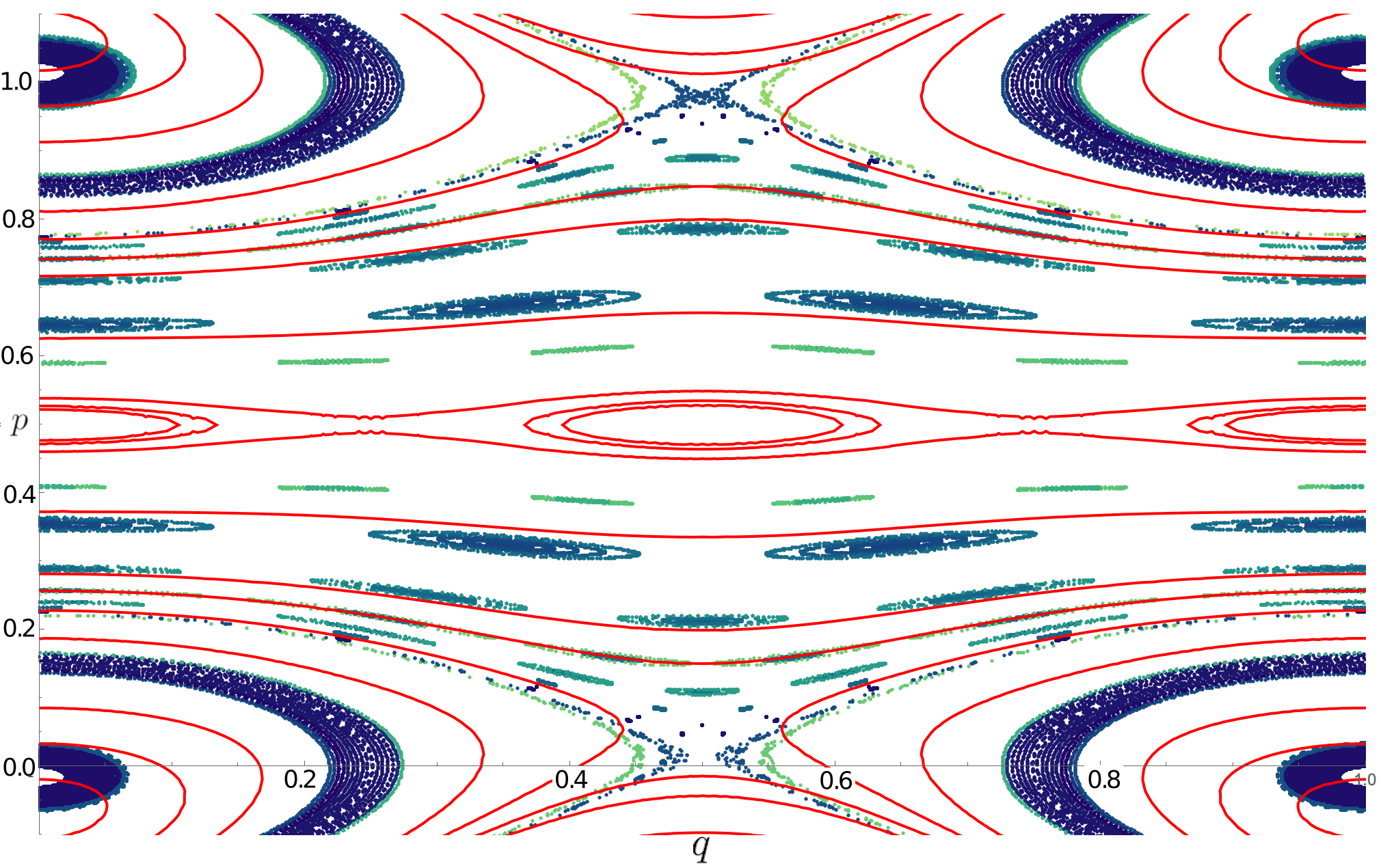}
  \end{subfigure}
  \caption{(s2-foliation)  Orbits satisfying the converse KAM condition for
 the s1-foliation (details as in the caption to \cref{fig:rfoliation}).
 The red curves (right panel) correspond to the Hamiltonian \eqref{eq:invariant2}.}
  \label{fig:s2foliation}
\end{figure}

Two constant-$\mu$ slices though the results for the s1-foliation are shown in \cref{fig:freqPlots}. These plots 
show the inverse of the time taken to pass the converse KAM condition $1/t_c$ as a function of $p_0$ $\mu = 0.015$ (left panel) and $\mu = 0.03$ (right panel). The peaks in each plot lie at the centre of the resonant tori not captured by the s1-foliation, or at the chaotic regions near the primary hyperbolic fixed point. These plots
are analogous to those of \cite{whiteDeterminationBrokenKAM2015}, who computed a rotation frequency for his phase vectors.

Finally, \cref{fig:s2foliation} shows computations for the s2-foliation.
This case is clearly different from the
previous results for initial conditions in the 1:2 resonant island (near $p=\tfrac12$).
Indeed, as we saw in \cref{fig:s2foliationLevelSets}, this foliation has contours
that are topologically similar to this island chain, and therefore many of these librational tori no longer satisfy the converse KAM criterion.
Unfortunately, there are now regions that satisfy converse KAM near the 0:1 and 1:1 elliptic orbits,
the blue regions near $p=0$ and $p=1$. 
These are apparently due to the fact that the contours of \eqref{eq:invariant2},
shown in \cref{fig:s2foliationLevelSets}, have a large period three
island trapped in the main resonances that is not present in the dynamics.

\begin{figure}[ht]
  \centering
    \includegraphics[width=0.49\linewidth]{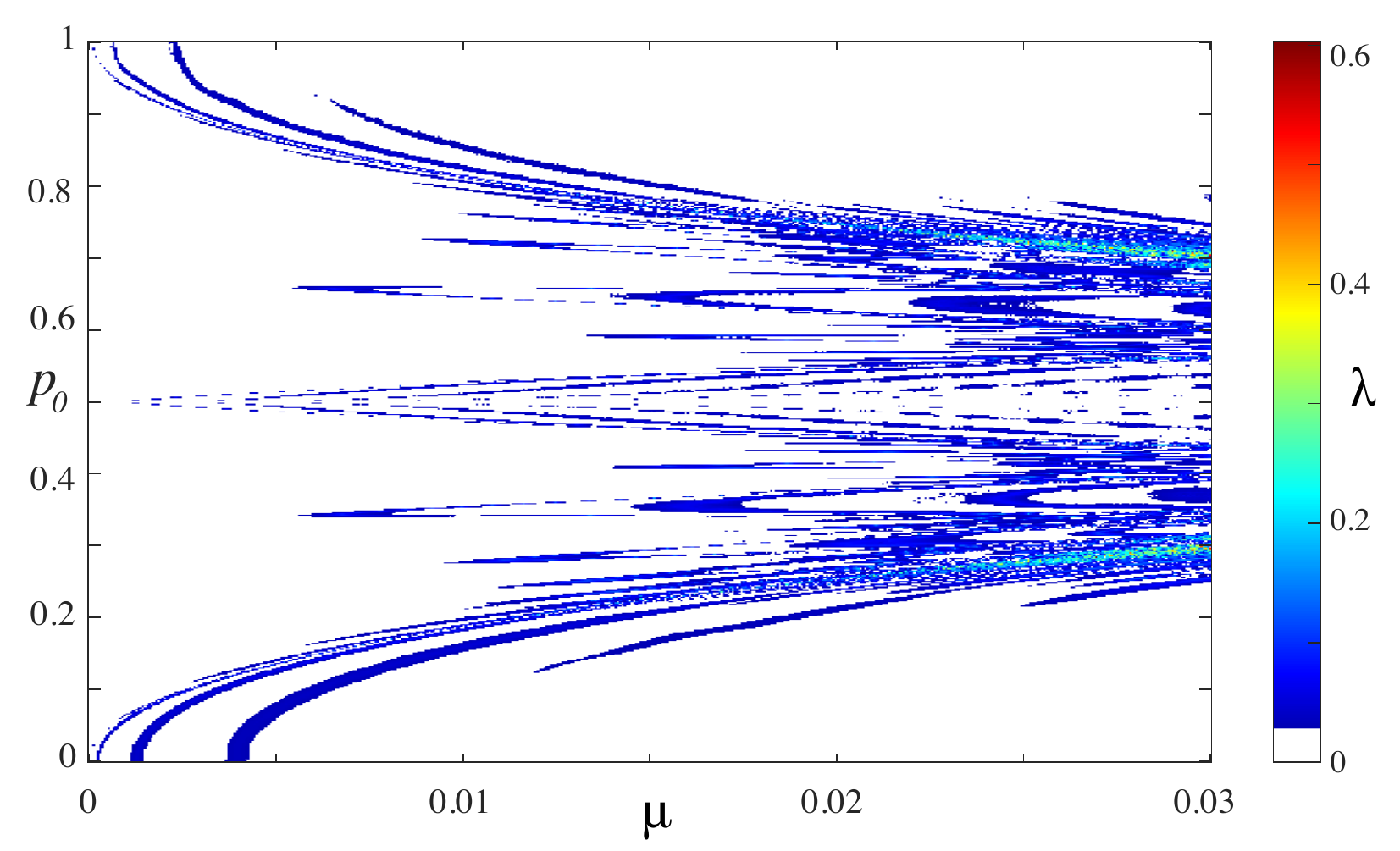}
\caption{(Lyapunov Exponent) The finite time Lyapunov exponent for the same initial conditions
as \cref{fig:rfoliation,fig:lfoliation,fig:pfoliation,fig:s1foliation,fig:s2foliation}.}
\label{fig:s1Lyapunov}
\end{figure}

As a comparision, a computation of the finite time Lyapunov exponent
is shown in \cref{fig:s1Lyapunov}. This shows the relative value of the computed Lyapunov exponent
for each initial condition $(\mu,p_0)$ integrated to $t = 150$. 
Note that almost all of these exponents are
quite small---those with $\lambda < 0.05$ are colored white in the figure. 
However, it would be difficult to declare an orbit chaotic when $\lambda < 0.2$
as it is for most of the parameters in the figure. Nevertheless, the domains where
$\lambda > 0.05$, correlate well with those that satisfy the converse KAM condition
for the l and s1-foliations, apart from the regions in the interior of the 1:2 and higher-order
resonances.

\begin{figure}[ht]
  \centering
  \begin{subfigure}{0.49\linewidth}
    \includegraphics[width=\linewidth]{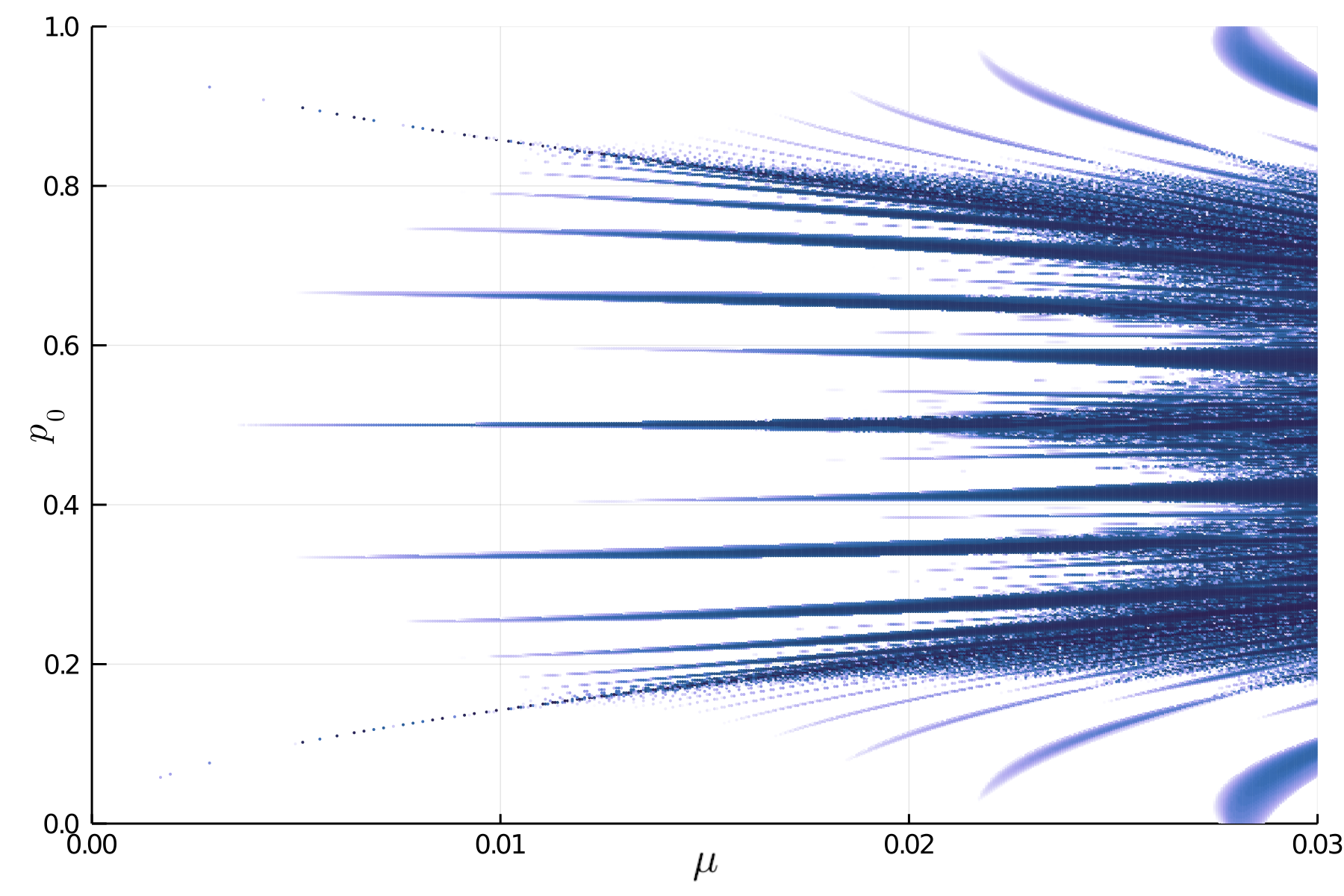}
  \end{subfigure}
  \begin{subfigure}{0.49\linewidth}
    \includegraphics[width=\linewidth]{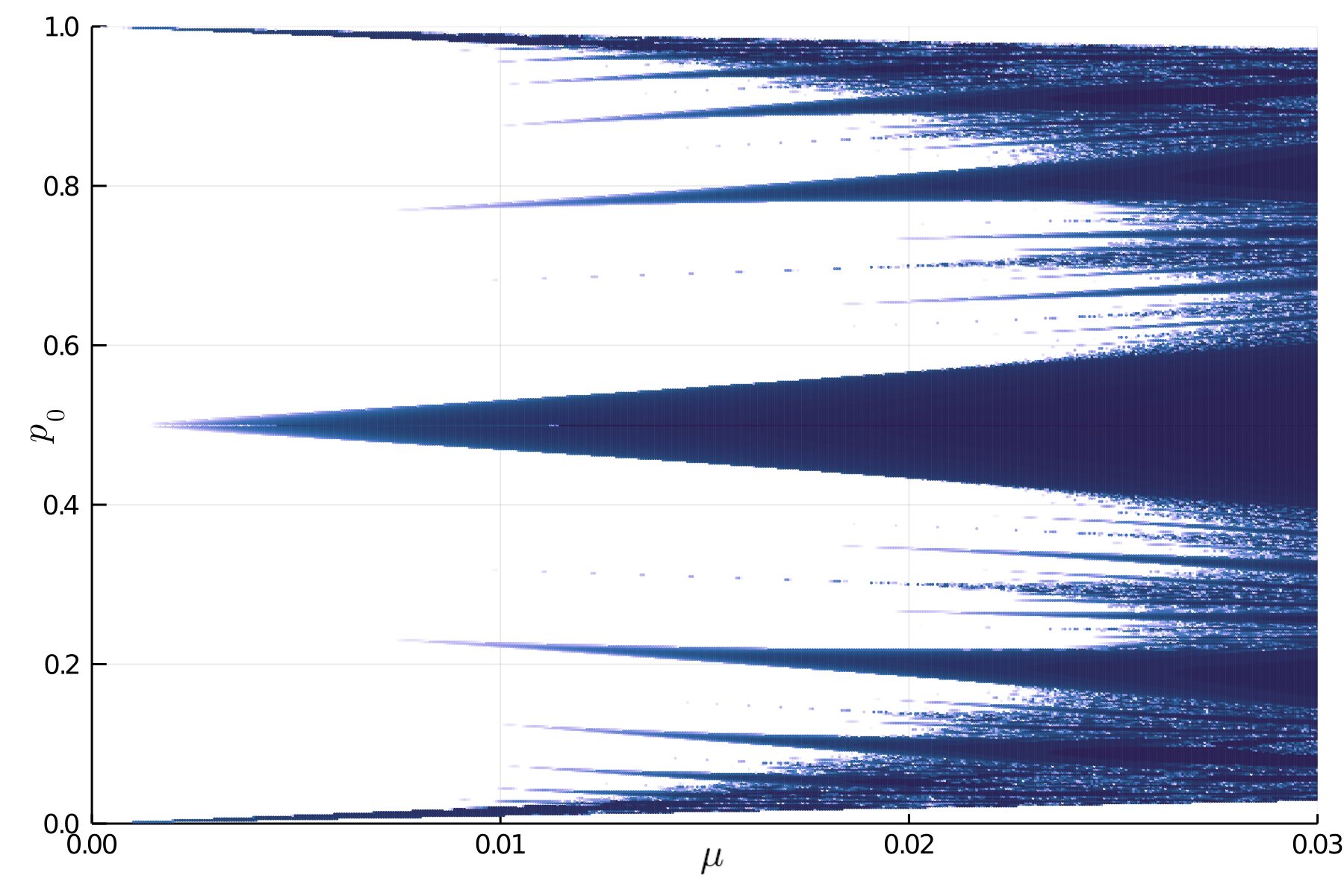}
  \end{subfigure}
  \caption{(s1-foliation)  Points that
    pass the converse KAM condition for the s1-foliation for Poincar\'e sections
    at $t = 0.25$ (left) and $t = 0.5$ (right). Compare these to \cref{fig:s1foliation}, which
    is on the section $t=0$.}
  \label{fig:diffqval}
\end{figure}

To reveal the effect of the choice of section on the results, the
converse KAM criterion with respect to the s1-foliation was tested on initial
conditions $(q_0,p_0,0)$ for $q_0 = 0.25$ and $0.5$. The results are depicted in
\cref{fig:diffqval}.

\subsection{Q-Flows}\label{sec:QResults}

\begin{figure}[ht]
  \centering
  \begin{subfigure}{0.49\linewidth}
    \includegraphics[width=\linewidth]{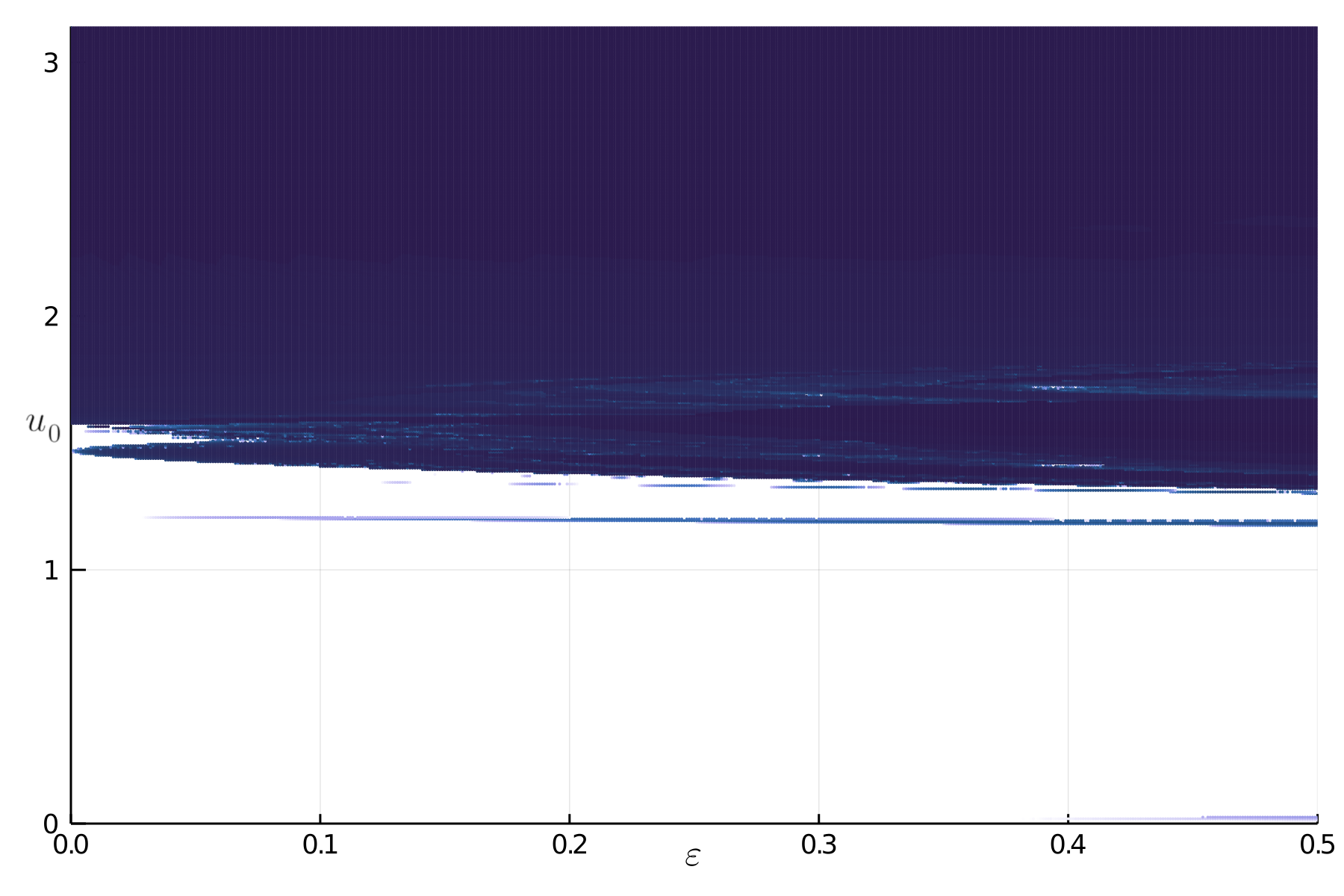}
    \caption{l-foliation, $ (\varepsilon,u_0)\in [0,0.5]\times [0,\pi]$}
    \label{fig:q4-lfoliation-full}
  \end{subfigure}
  \begin{subfigure}{0.49\linewidth}
    \includegraphics[width=\linewidth]{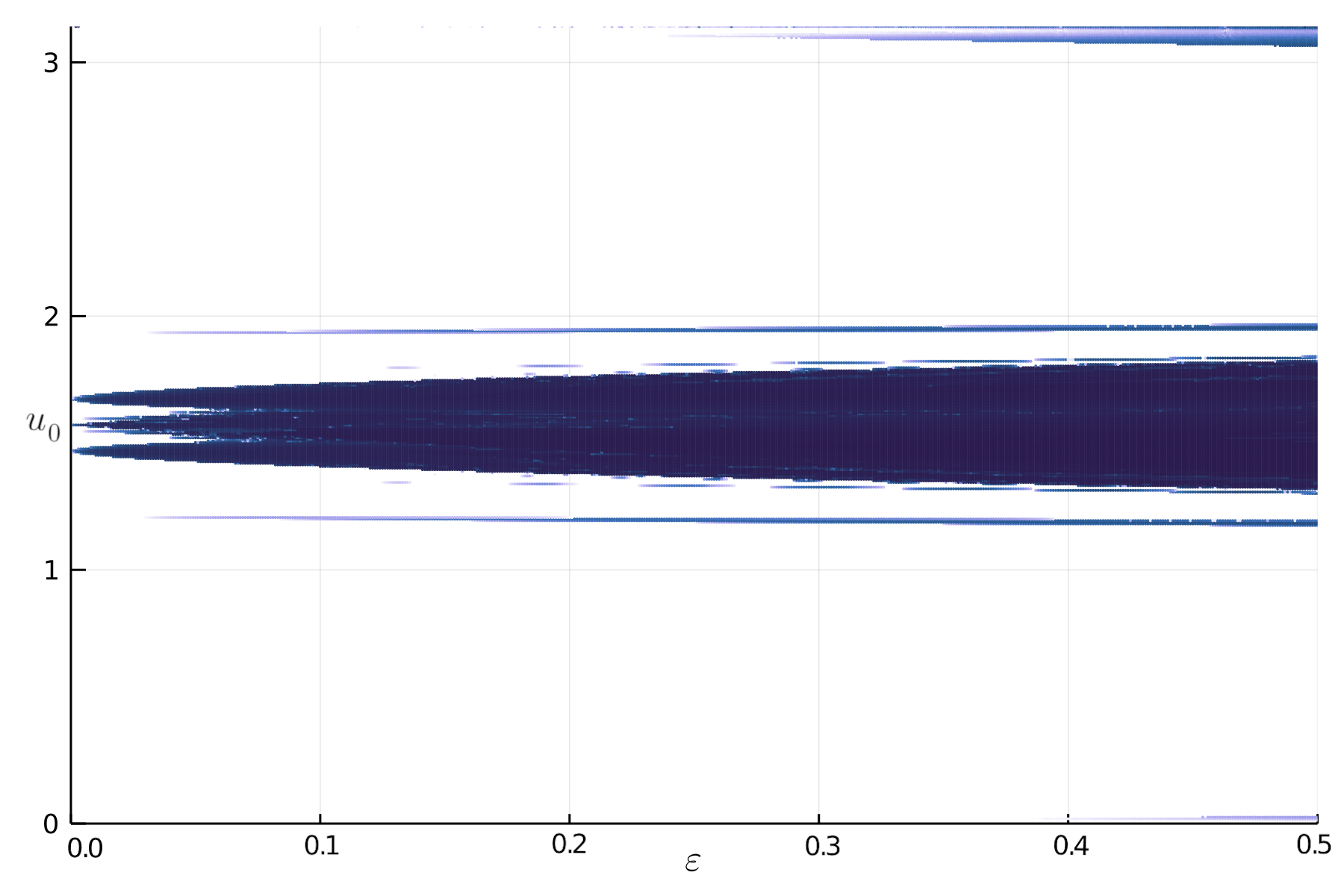}
    \caption{$\psi$-foliation, $ (\varepsilon,u_0) \in [0,0.5]\times [0,\pi]$}
    \label{fig:q4-pfoliation-full}
  \end{subfigure}
  \begin{subfigure}{0.49\linewidth}
    \includegraphics[width=\linewidth]{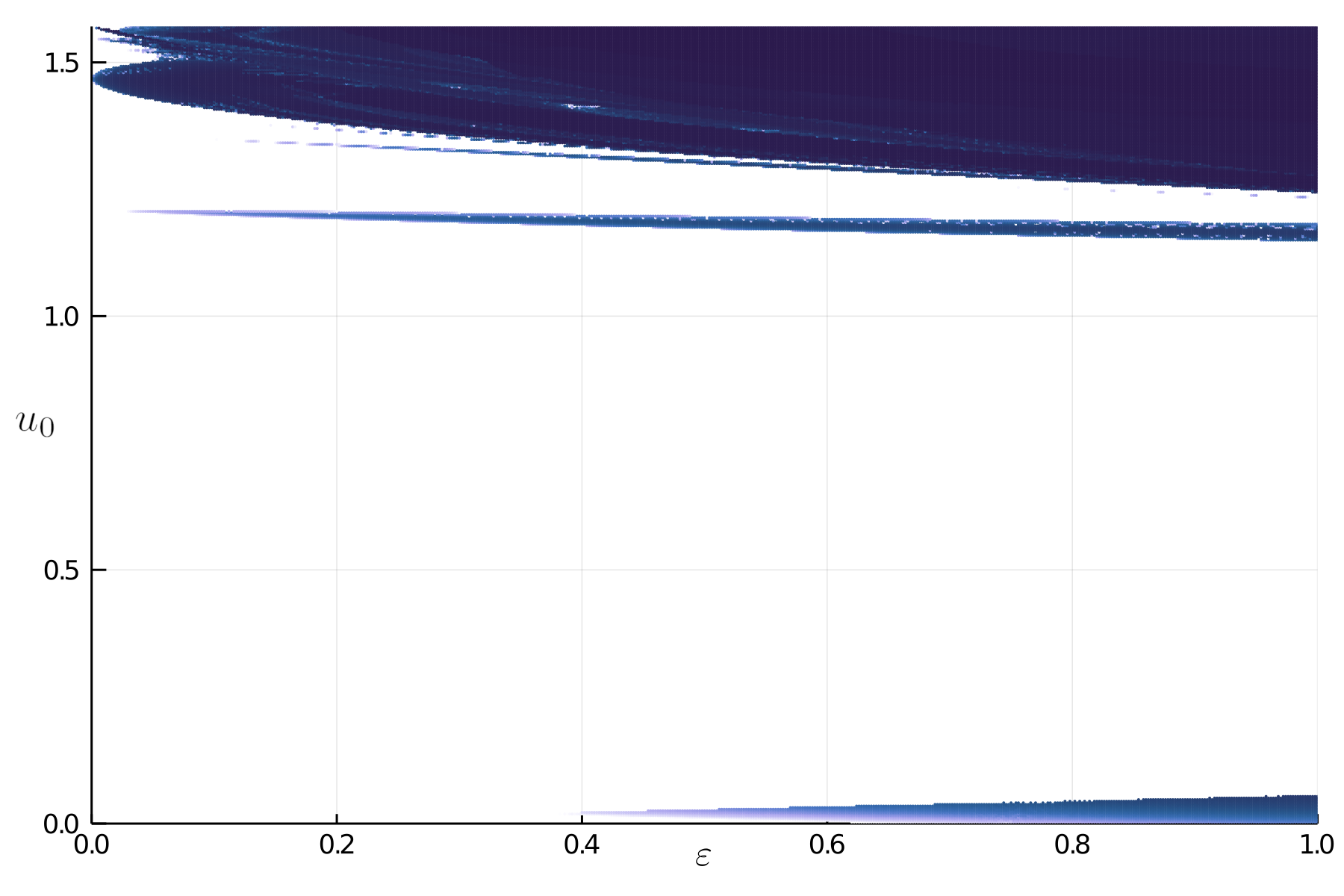}
    \caption{l-foliation, $(\varepsilon,u_0)\in [0,1]\times[0,\pi/2]$}
    \label{fig:q4-lfoliation-small}
  \end{subfigure}
  \begin{subfigure}{0.49\linewidth}
    \includegraphics[width=\linewidth]{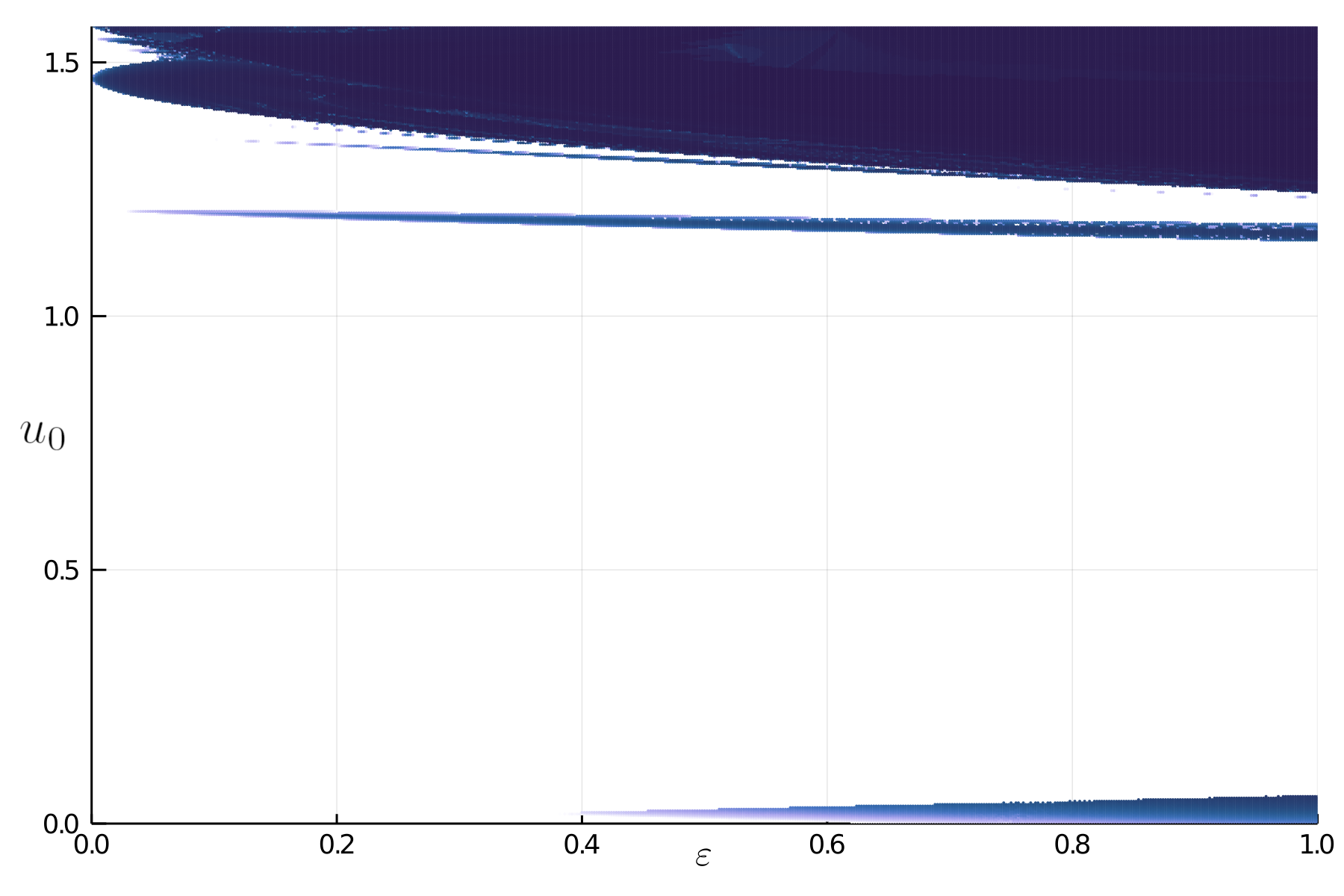}
    \caption{$\psi$-foliation,  $(\varepsilon,u_0)\in [0,1]\times[0,\pi/2]$}
    \label{fig:q4-pfoliation-small}
  \end{subfigure}
  \caption{Orbits of the Q-flow with $q = 4$, 
  that pass the converse KAM condition for 
  time $t_c \le150$ with initial conditions on the line $(u_0,u_0,0)$. 
  Panels (a) and (c) use the l-foliation and (b) and (d) use the $\psi$-foliation.
  The points are colored according to the time taken to
  transversality violation, with a darker blue indicating a shorter time.}
    \label{fig:q4plots}
\end{figure}

\begin{figure}[ht]
  \centering
  \begin{subfigure}{0.49\linewidth}
    \includegraphics[width=\linewidth]{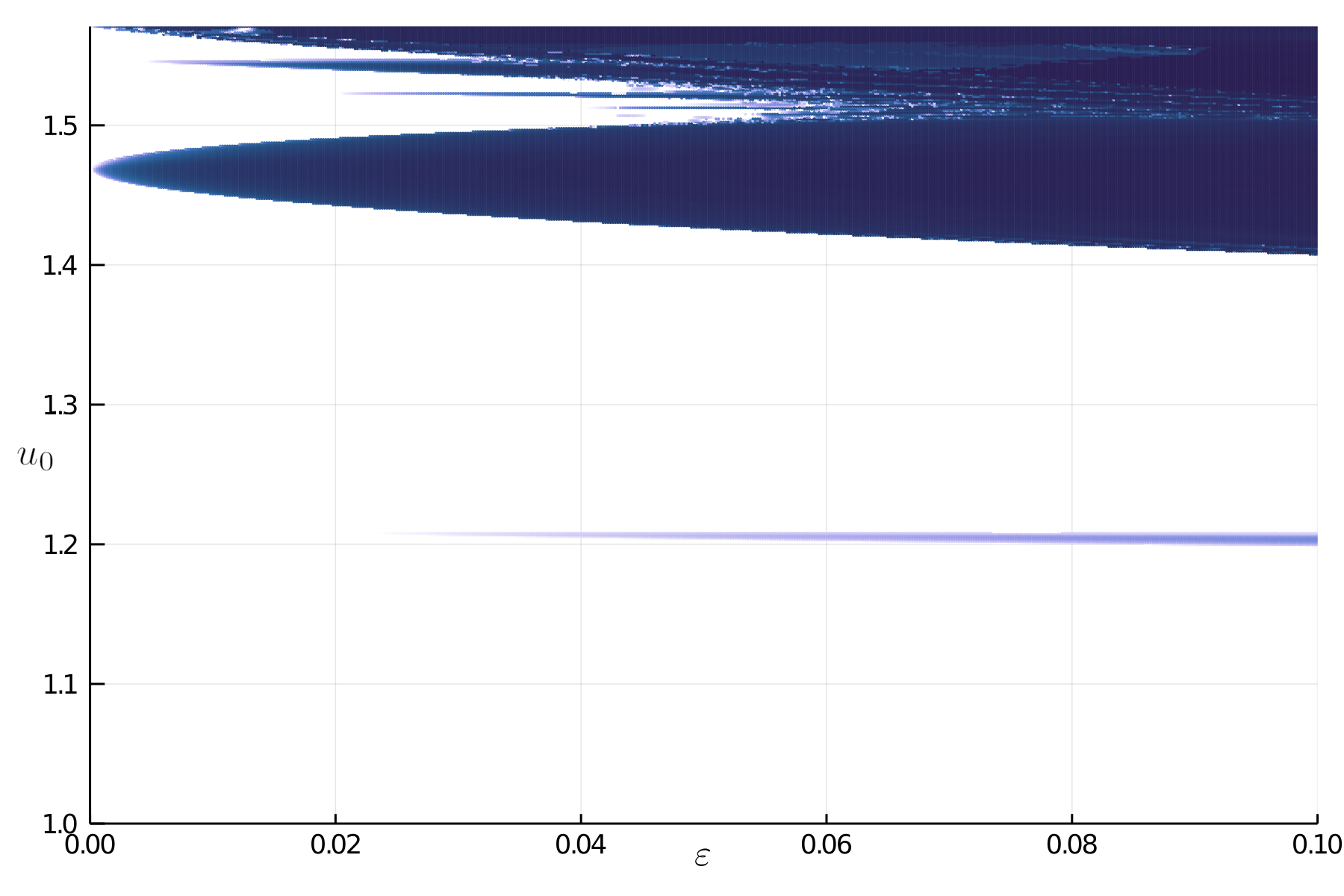}
    \caption{$\psi$-foliation}
    \label{fig:q4-pfoliation-zoom}
  \end{subfigure}
  \begin{subfigure}{0.49\linewidth}
    \includegraphics[width=\linewidth]{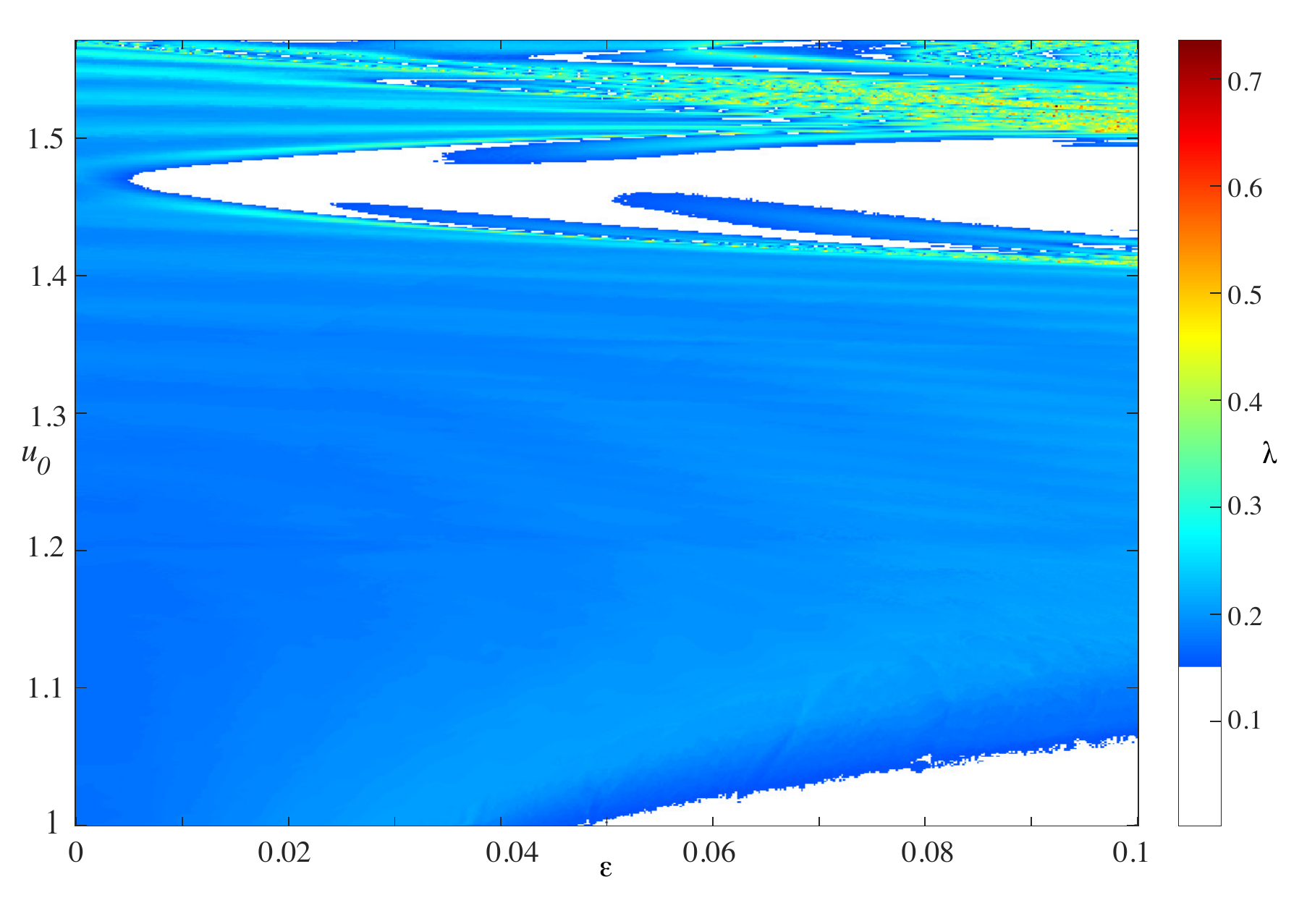}
    \caption{Maximal Lyapunov exponent}
    \label{fig:q4-lyapunovs}
  \end{subfigure}
  \caption{Q-flow with $q = 4$ for $(\varepsilon,u_0)\in [0,0.1]\times[1,\pi/2]$.
  Panel (a) is an enlargement
  of \cref{fig:q4plots}(d). Panel (b) is the maximal Lyapunov exponent for the orbits of (a)
  using the initial tangent vector $v_0 = (0,1,0)$.}
  \label{fig:q4enlarged}
\end{figure}

\begin{figure}[ht]
  \centering
  \begin{subfigure}{0.49\linewidth}
    \includegraphics[width=\linewidth]{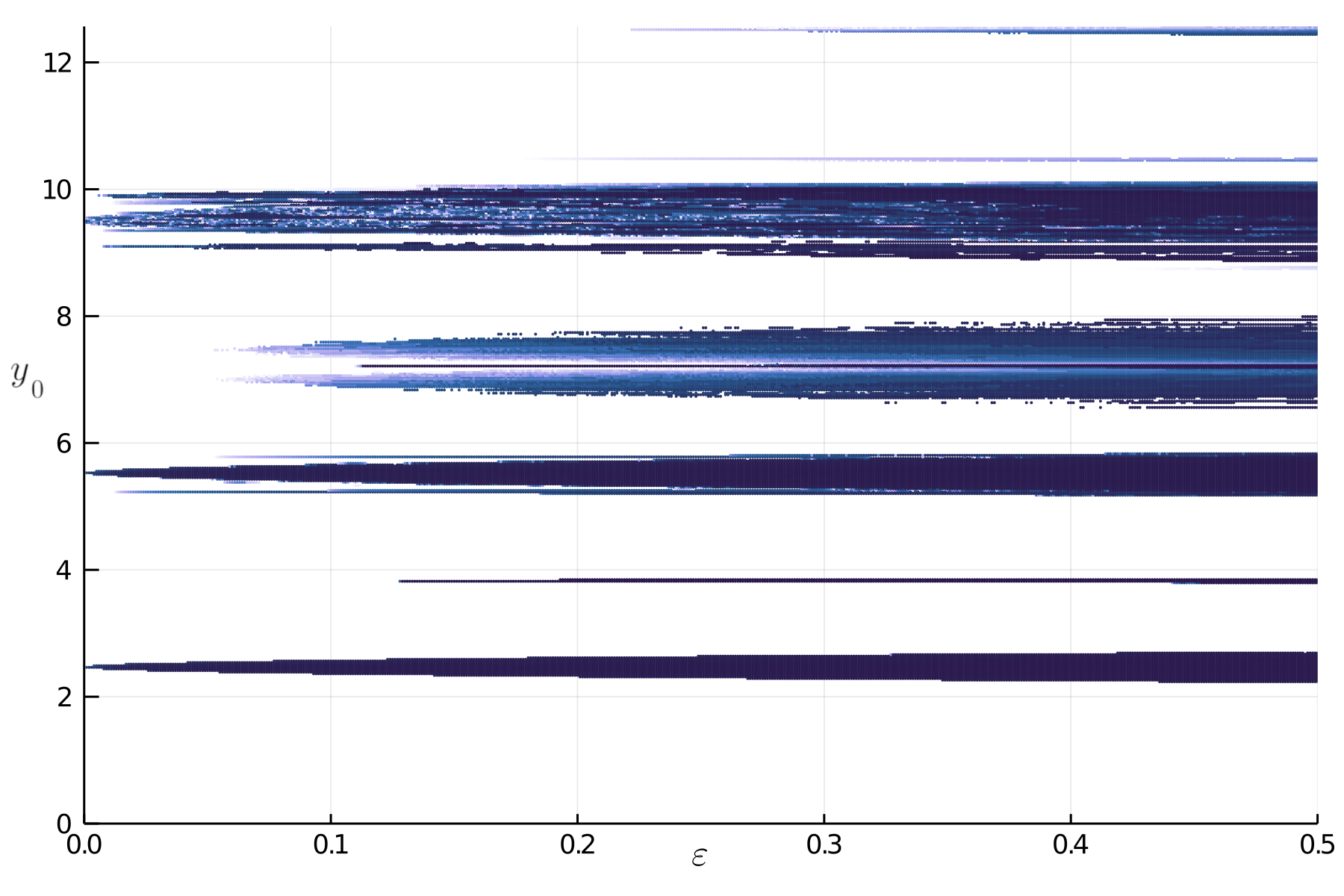}
    \caption{$\psi$-foliation}
    \label{fig:q5-pfoliation-full}
  \end{subfigure}
  \begin{subfigure}{0.49\linewidth}
    \includegraphics[width=\linewidth]{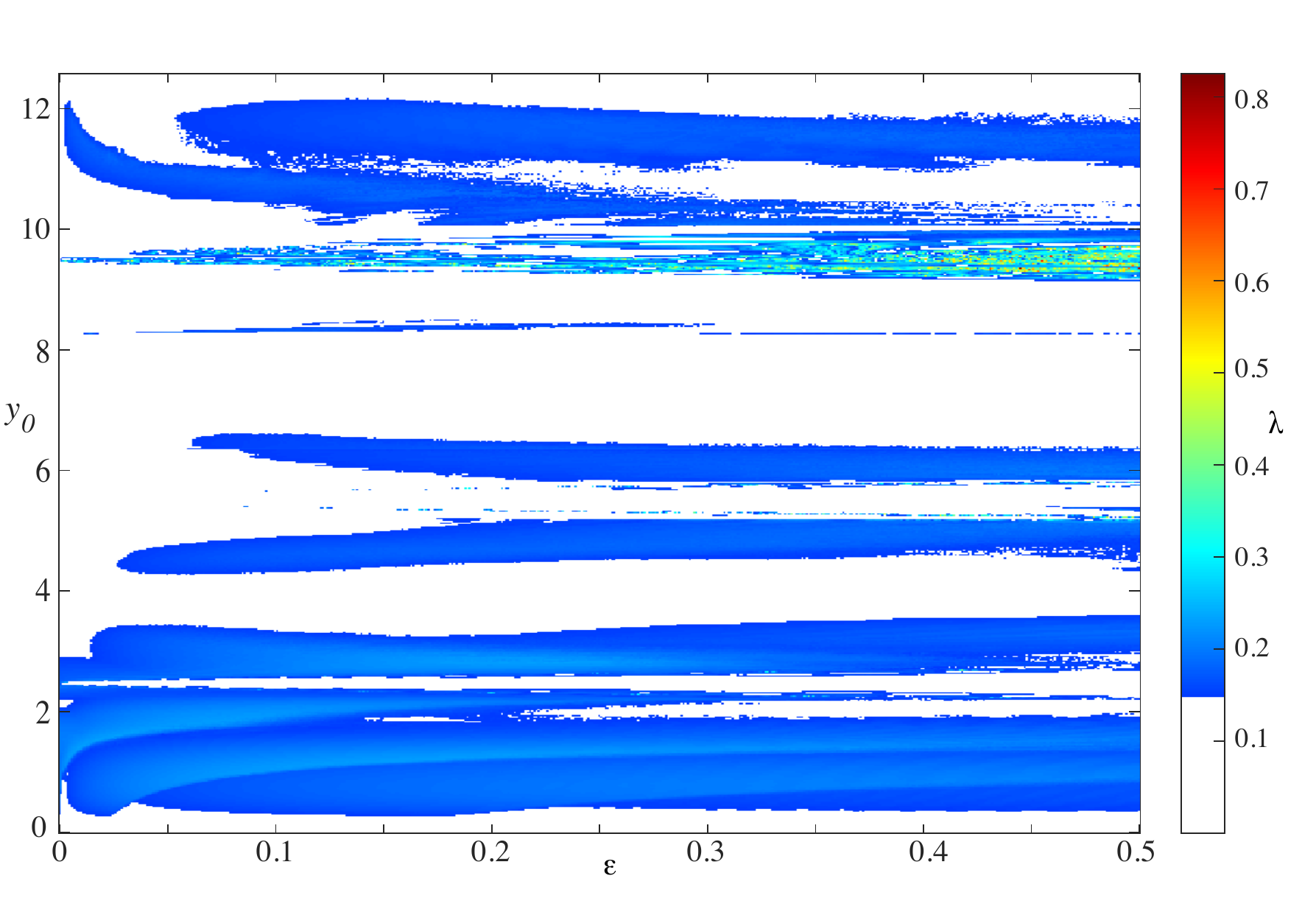}
    \caption{Maximal Lyapunov exponent 
    }
    \label{fig:q5-lyapunovs}
  \end{subfigure}
   \caption{Q-flow with $q=5$  for initial conditions on the line $(0,y_0,0)$ and
    $(\varepsilon,y_0)\in [0,0.5]\times[0,4\pi]$. Panel (a) shows orbits that pass 
    the converse KAM condition for the $\psi$-foliation for $t\leq 150$.
    Panel (b) depicts the maximal Lyapunov exponent for the orbits of Panel (a).}
  \label{fig:q5plots}
\end{figure}

\begin{figure}[ht]
  \centering
  \begin{subfigure}{0.49\linewidth}
    \includegraphics[width=\linewidth]{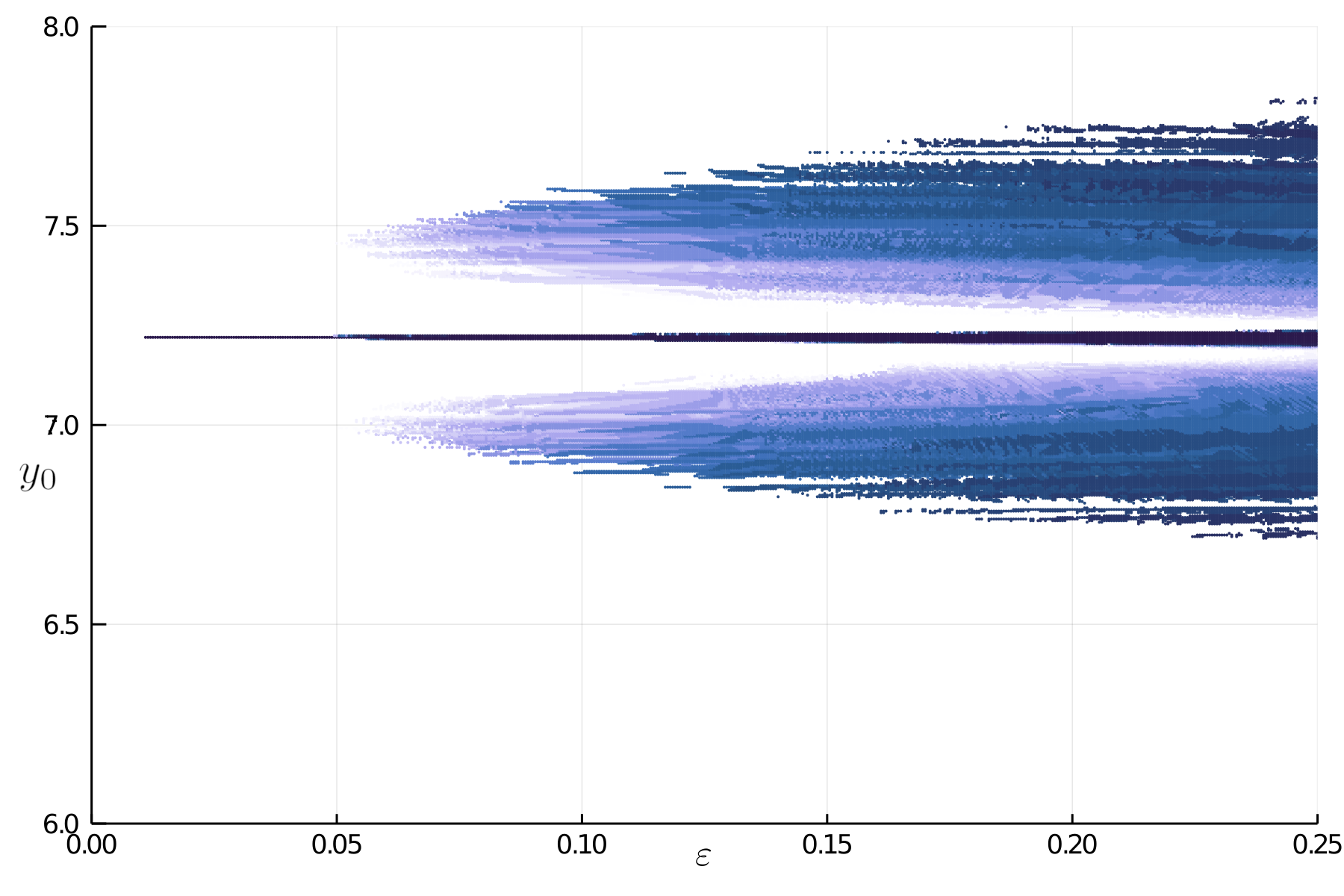}
    \caption{Enlargement of (a)}
    \label{fig:q5-pfoliation-small}
  \end{subfigure}
  \begin{subfigure}{0.49\linewidth}
    \includegraphics[width=\linewidth]{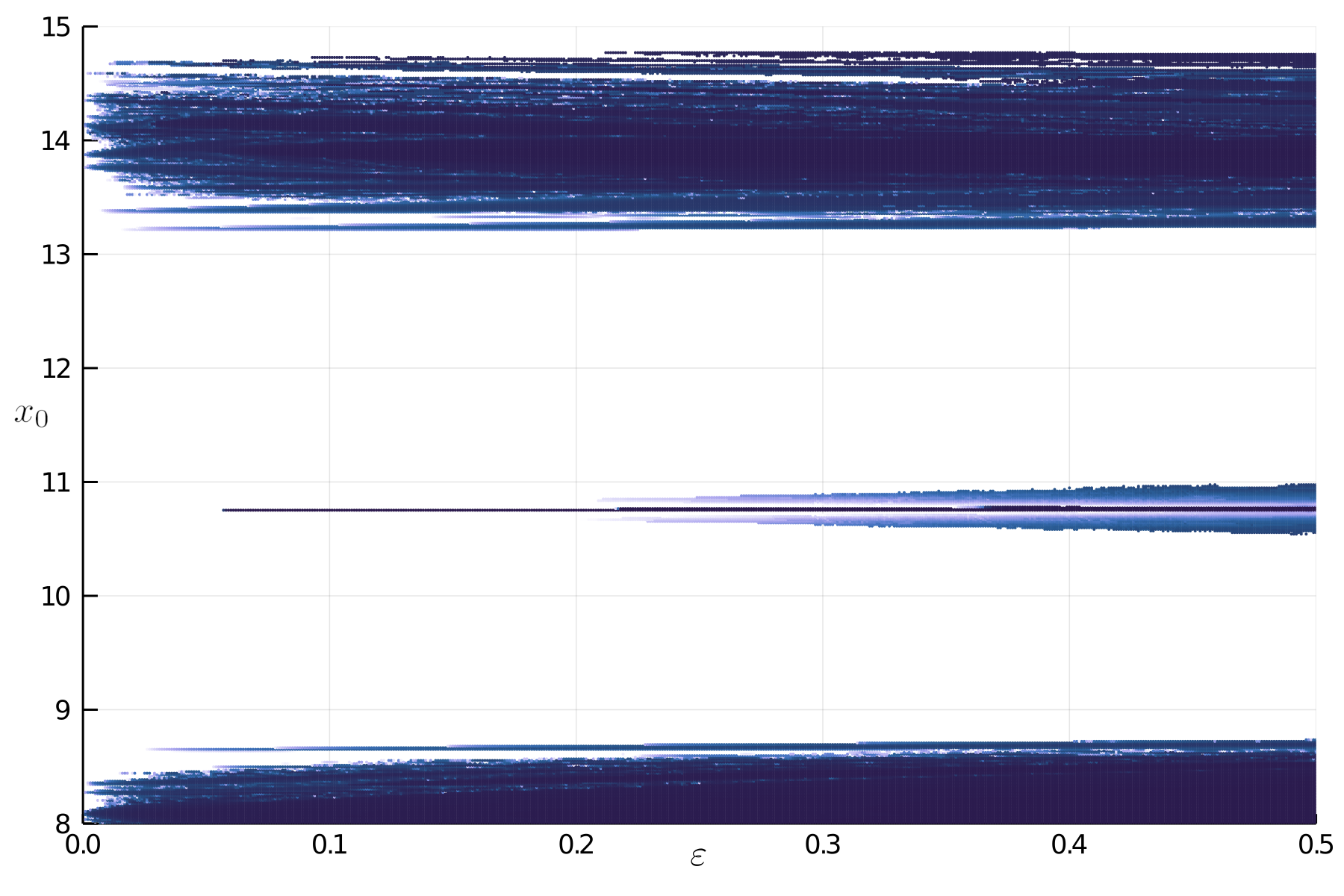}
    \caption{Initial conditions on the line $(x_0,0,0)$.}
    \label{fig:q5-pfoliation-x}
  \end{subfigure}
  \caption{Q-flow with $q=5$. Panel (a) is an enlargement of \cref{fig:q5plots}(a). 
  Panel (c) shows orbits staring on the line $(x_0,0,0)$ that pass the converse KAM condition.}
   \label{fig:q5enlarged}
\end{figure}

Computations of the converse KAM condition for the Q-flows, \eqref{eq:qflow}, are shown in
\cref{fig:q4plots,fig:q4enlarged} for $q=4$ and \cref{fig:q5plots,fig:q5enlarged}
for $q=5$. In each case, we
choose $500$ initial conditions along a line in the phase space at $z = 0$, and check the sign of $K_{qf}$ \eqref{eq:KqfAnddw} for $t_c \le 150$. As for the two-wave results, each point that satisfies the
converse KAM condition is indicated in blue, and a darker blue indicates
a smaller value of $t_c$.

For \cref{fig:q4plots}, the ordinate, $u_0$,
represents initial conditions along the line $(x,y,z) = (u_0,u_0,0)$ and
the abscissa is the perturbation strength, $\varepsilon$. Panels (a) and (b)
compare the l-foliation to the $\psi$-foliation, respectively, for $(\varepsilon,u_0)\in[0,0.5]\times[0,\pi]$.
As was seen in \cref{fig:qflows}(a), when $q=4$ the initial conditions with $u_0 \in
[0,\pi/2)$ correspond the lattice cell encircling the origin. For these initial
conditions, the results of the two foliations are nearly identical. This can
also be seen in panels (c) and (d), which display only orbits with $u_0 \in [0,\pi/2]$ but
extend the range of $\varepsilon$ to $[0,1]$. Both
foliations correctly assess the expected growing chaotic zone that forms around
the separatrix from the saddles at $(x,y) = (0, \pm \pi)$ and $(\pm \pi, 0)$.
They also indicate a narrow tongue near $u_0 = 1.2$ that starts near
$\varepsilon = 0.031$ in which there are no tori appropriate to the radial
foliations.

When $\varepsilon = 0$, the orbits of initial conditions with $u_0 \in (\pi/2,\pi]$ 
lie in a second cell of $\psi$ contours that is centered at $(\pi,\pi)$. 
The implication is that orbits may satisfy
the converse KAM condition with respect to the l-foliation, as
seen in \cref{fig:q4plots}(a), even though they may still lie on tori---now
enclosing the elliptic orbit near $(\pi,\pi)$. Indeed when $\varepsilon \ll 1$,
most of the orbits in this region do lie on such tori. By contrast, the $\psi$-foliation shown in
\cref{fig:q4plots}(b) correctly indicates that there are still tori
in this region, since the foliation in this case is essentially radial with
respect to the origin at $(\pi,\pi)$.

A further enlargement of \cref{fig:q4plots}(d) is shown in \cref{fig:q4enlarged}(a).
Note that there are many small tongues of destroyed tori just above $u_0 = \tfrac{\pi}{2}$. These
results can be compared to the maximal Lyapunov exponent shown in the right panel.
Here initial conditions with $\lambda<0.15$ are---somewhat 
arbitrarily---deemed to be regular and colored white. A large tongue that
emanates from $(u_0,\varepsilon) = (1.47,0)$ is a regular island chain. 
As for the two-wave case, it is difficult, using the finite-time Lyapunov exponent
for $t=150$, to obtain a sharp distinction between regular and chaotic orbits .

Orbits that pass the converse KAM condition for the $q = 5$ case,
which, as we saw in \cref{fig:qflows}(a), has
quasi-crystal symmetry, are displayed in
\cref{fig:q5plots}(a) using the $\psi$-foliation.
Here the initial conditions lie on the line
$(0,y_0,0)$. The computations show that there are many tongues that have no tori relative to this foliation.
This can be compared to the computations of the Lyapunov exponent in panel (b). Note
that, again, the converse KAM condition is much more sensitive to the destruction of
tori as it would require a longer computation to get convergence of the exponent.

Figure \ref{fig:q5enlarged}(a) shows an enlargement (near $y_0 = 7$) of two of the tongues of \cref{fig:q5plots}(a).
A different set of initial conditions,$(x_0,0,0)$ for $x_0$
near $11$ and $\varepsilon \in [0,0.5]$, are shown in panel (b).
This panel gives a closer view of the elliptic orbit
that passes near $(10.75,0,0)$.

\section{Discussion}\label{sec:discussion}

We have shown computations of the converse KAM condition for two models
using a variety of foliations. We summarize some our results here.

\begin{itemize}
\item \textit{For the rotational tori, our results agree with those of MacKay.}

  The results for the r-foliation in the two-wave model
  (\cref{fig:rfoliation}) are identical to those of
  \cite[Fig.~4]{mackay_criterion_1989}. This foliation detects when
  a given point does not lie on a rotational invariant torus that is a graph over $(q,t)$.
    
\item \textit{The selection of an appropriate foliation is crucial.}

  As was discussed in \S\ref{sec:Foliations},
  the results obtained are sensitive to the choice of foliation. 
  In \S\ref{sec:DWResults}, we presented results for four different foliations
  of the two-wave flow, and, in \S\ref{sec:QResults}, results for two foliations of
  the Q-flow. 
 
  For the two-wave case, one interesting aspect of the l-foliation, \eqref{eq:lFoliation}, shown in 
  \cref{fig:lfoliation}, is that points not on librational tori could be expected to dependently 
  satisfy the converse KAM condition. However, it is apparent in the figure that those points on 
  rotational tori are not dependent.
  This agrees with the discussion in \S\ref{sec:Foliations}: the invariant tori of null-homology 
  with respect to $\cM = M\setminus \Sigma$ for the singuarlity set $\Sigma$ 
  will satisfy the converse KAM condition dependently. 
  Note that for the l-foliation, $M = \T\times\R\times\T$ and $\Sigma= \set{(0,0,t)}{t\in\T}$; thus 
  the rotational tori are not of null-homology, and these tori are transverse to the l-foliation.

  Hidden in \cref{fig:lfoliation} and \cref{fig:pfoliation}, for the
  l- and p-foliations, are small tongues around $p_0 = 0$
  that satisfy the converse KAM condition. These result from the elliptic periodic orbit
  that begins at $(q,p) = (0,0)$ for $\mu = 0$, having a nontrivial oscillation when $\mu >0$.
  By contrast, the singular set of the s1-foliation more accurately tracks the drift of this orbit,
  so that $\cM = M\setminus \Sigma$ agrees better
  with the homology of the true orbits of the vector field $v$.

  The s1-foliation (\cref{fig:s1foliation})  also has an elliptic orbit near $(q,p) = (0,1)$.
  Thus, rotational tori around the elliptic orbit of the vector field near $(0,1)$ are 
  not of null-homology with respect to this foliation. This gives a more accurate picture of
  both the rotational tori and the librational tori in the two forced resonances.
  
\item \textit{It is difficult to conclude a point lies on \emph{no} invariant surface.}

  In a typical Hamiltonian or volume-preserving flow resonances are dense. Many of these
  will form island chains, and if one were to find points that do not lie on \textit{any} invariant
  surface, one would have to know the locations of \textit{all} of these resonant island chains.
  
  For example, the s1-foliation is constructed to ensure that points on rotational, 
  and librational tori around $(0,0)$ and $(0,1)$, are not necessarily dependent.
  The s2-foliation was designed to
  also ensure that the 1:2 resonant tori near $(q,p)=(0,\tfrac12)$ are not dependent. However, 
  it fails to appropriately test for librational tori in the primary elliptic regions; indeed,
  there are points passing the converse KAM criterion near 
  $p = 0$ and $p = 1$ in \cref{fig:s2foliation} even though a phase portrait shows there are tori.

  It might be possible to use renormalization procedures to better design a foliation adapted to 
  some resonances, but since there are typically an infinite number of
  homologies of tori, it seems impossible to detect them all.
  The methods discussed here could be used to design a foliation appropriate to a particular homology in
  a localized region of phase space.

\item \textit{The results depend on choice of initial conditions.}

  Figure~\ref{fig:diffqval} shows how the results depend on the
  choice of section for the two-wave model since we only have chosen initial conditions along
  a line $(0,p_0,t_0)$. Though the overall structure
  is similar to the $t=0$ section in \cref{fig:s1foliation}, there
  are some notable differences.
  
  The $t=0.25$ section reveals more details about trapped islands in the librational
  regions about $p_0 =0$ and $1$ because the initial conditions are 
  closer to the hyperbolic points of these resonances.
  
  In the $t=0.5$ section there is almost a complete lack of orbits trapped in the librational
  regions because the line of initial conditions goes through the
  hyperbolic orbits of these main island chains. If one wants to
  understand all possible islands, one should take care to choose appropriate
  initial conditions. 

  The importance of the choice of initial conditions is also evident in the Q-flow for $q=5$
  by noting the various islands not sampled when taking $x_0 = 0$ in
  (\cref{fig:q5-pfoliation-full}) compared to $y_0 = 0$
  (\cref{fig:q5-pfoliation-x}). 
  
\item \textit{Finite-time Lyapunov exponents can be less effective.}

  Finite-time Lyapunov exponents, shown in \cref{fig:s1Lyapunov,fig:q4-lyapunovs,fig:q5-lyapunovs},
  are often used to detect chaotic regions. We have seen that the converse KAM criterion accurately detects 
  such regions as well (though it also detects orbits on non-transverse tori).
  In our computations, we used the same time $t = 150$ in both computations, and the distribution
  of Lyapunov exponents does not give a clear distinction between regular and chaotic orbits.
  Celletti et al \cite{Celletti07} note that ``in contrast to the Lyapunov exponents,
  the nonexistence criterion allows us to distinguish between chaotic motions [and]
  librational tori, and rotational tori, [thus] extending the idea of the fast Lyapunov
  indicators \cite{Guzzo02}.'' This statement is clearly evident as the converse
  KAM method detects the nonexistent of tori of certain homology.
  
\item \textit{The inverse of the converse KAM time gives a measure of librational frequency.}

	As White notes \cite{whiteDeterminationBrokenKAM2015}, 
	the rate of rotation of the phase vector, $\xi_t$, for orbits in islands gives a
	measure of the librational frequency inside the island. This should correlate with a times scale for 
	flattening of a density profile of a passive scalar inside these island chains. 
	As such this could be a useful measure for plasma confinement devices.

  \end{itemize}

In conclusion, the numerical application of Thm.~\ref{thm:MacKay} or Thm.~\ref{thm:CKAMVolume} is
extremely useful if one is looking in particular at the existence of invariant
surfaces of a given homological type.

\subsection*{Future Directions}

Whilst converse KAM theory dates back several decades, we hope that with an appropriate choice of foliation, it can be directly applied to the optimization of fusion confinement devices such as stellarators. For a stellarator to be effective at trapping charged particles, there should be many invariant tori that can act as effective barriers for particle transport. Indeed, integrability or near-integrability is an often used criterion for optimal design of these devices. In this sense, we hope that converse KAM theory can be used to create a rapidly computable measure that would help to design configurations that maximize integrability.

Another potential application is to the construction of divertors; configurations that enhance particle
transport at controlled points on the outer surface of the plasma. These allow better control of the plasma boundary and help to control impurities \cite{Boozer15}. One type of divertor uses a magnetic field with a hyperbolic orbit on the plasma edge \cite{Boozer18}. The positioning of this hyperbolic orbit must be precise. Hence, a foliation could be constructed to test potential magnetic fields for their ability to adhere to the desired divertor configuration.

While converse KAM theory has been applied to detect nonexistence of invariant Lagrangian tori in \cite{mackay_converse_1989}, it does not seem to have been applied more generally to higher-dimensional systems.
The ideas we have discussed can search for the nonexistence of invariant surfaces with a fixed Lipschitz constant. 
The core issue with checking Thm.~\ref{thm:MacKayExtended} for $n > 3$ is the codimensionality of the 
line $\Span(\eta_t)$, which is one for $n=3$ in $\cP$, but in general is $n-2$. Thus, $\xi_t$ will almost never cross the line $\eta_t$. Given a constraint on the Lipschitz constant of possible invariant surfaces, one then asks if $\xi_t$ crosses a cone about $\eta_t$. This is again a codimension-one problem.

A converse KAM condition is said to be  exhaustive if all orbits that do not lie on surfaces transverse to the chosen foliation $\cF$ eventually satisfy the converse KAM condition. This notion is addressed in \cite{mackay_finding_2018} where it is noted that orbits on some invariant sets (for example cantori) may not satisfy the converse KAM condition even though these are not invariant surfaces. The remedy is to amend the procedure with the so-called ``killends" condition \cite{MacKay85}. As argued in \cite{mackay_finding_2018}, this condition allows any orbit in a volume-preserving system not lying on an invariant surface transverse to $\cF$ to be detected. However, in general, exhaustiveness is only conjectured. In the future we hope to develop the killends condition to help address this question.

\appendix
\section{Global Removal of Resonances}\label{sec:Invariants}
The approximate invariant \eqref{eq:approxInv} is the result of a technique used
in \S 2.4d of \cite{lichtenberg_regular_1992}, where it is called ``global
removal of resonances''. 
We assume that $M = (\T \times \R)^2$ with coordinates $(\theta_1,I_1,\theta_2,I_2)$,
where the angles have period one.
The technique works for any two degree-of-freedom Hamiltonian that can be
written as a perturbation of the integrable Hamiltonian $H_0: M \to \R$,
\begin{equation}
  \label{eq:Ham0}
  H_0(I_1,\theta_1,I_2,\theta_2) = S_1(I_1) + I_2.
\end{equation}
Note that this requires that the integrable system have global action-angle variables,
$I_i$ and be separable. We also assume, for simplicity, that the second degree of freedom
is harmonic, and without loss of generality set the harmonic frequency to one. 
Note that such a system also arises from a periodically time-dependent, $1\tfrac12$ degree-of-freedom system,
such as the two-wave model of \S\ref{sec:DoubleWave}, in extended in extended phase space.

Assume we have a Hamiltonian of the form
\begin{equation}
  \label{eq:HamSeries}
  H(\theta_1,I_1,\theta_2,I_2) = H_0(I_1,I_2) + \varepsilon H_1(I_1,\theta_1,I_2,\theta_2) + \dots.
\end{equation}
The key observation is that
when $\varepsilon \ll 1$, $I_1$ will vary slowly; that is,
$I_1$ is an adiabatic invariant. The question then becomes; can we build a
`better' (read `higher order'), adiabatic invariant $J$ from $I$?

A function $J:M \to \R$ is invariant whenever
\begin{equation}
  \label{eq:isInvariant}
 \{H, J\} = 0 ,
\end{equation}
where  $\{\cdot,\cdot\}$ is the canonical Poisson bracket,
\[
	\{F,G\} = \sum_{i=1}^n \frac{\partial F}{\partial I_i}  \frac{\partial G}{\partial \theta_i}
	        -  \frac{\partial F}{\partial \theta_i}  \frac{\partial G}{\partial I_i},
\]
and we use the convention that the Hamiltonian vector field $X_H = \{H, \cdot \}$.

Expanding $J$ in $\varepsilon$,
\[ 
	J = J_0 + \varepsilon J_1 + \dots,
\]
the most obvious direction to proceed is to expand $\{H,J\}$ in $\varepsilon$
and solve the equation at order-by-order in $\varepsilon$. This approach is
sufficient for low order calculations, however, if one wants to compute $J$ to
higher order, a Lie transform approach is more convenient. 

Let us first use the
direct method to compute the lowest order invariant. Expanding
\eqref{eq:isInvariant} to first order in $\varepsilon$ gives the two equations,
\begin{align}
  \label{eq:firstOrderInvariant}
  \{H_0, J_0\} &= 0 \\
  \{H_0,J_1\} &=  -\{H_1, J_0\} = \{J_0,H_1\}.
\end{align}
The first condition is satisfied for any $J_0 = J_0(I_1,I_2)$. For the examples
considered here, it is sufficient to consider $J_0$ a function of $I_1$ only.
Expanding the second equation yields,
\begin{equation}
  \label{eq:secondInvEquation}
  \omega_1(I_1)\partial_{\theta_1} J_1 + \partial_{\theta_2} J_1 = J_0^\prime(I_1) \partial_{\theta_1} H_1.
\end{equation}
where $\omega_1 = S_1'$, is the action-dependent frequency of \eqref{eq:Ham0}.
Now, let,
\begin{equation}
  \label{eq:FSeries}
  H_1 = \sum_{n,m} h_{n,m}(I_1,I_2) \me^{2\pi i(n\theta_1 + m\theta_2)},\qquad 
  J_1 = \sum_{n,m} k_{n,m}(I_1,I_2) \me^{2\pi i(n\theta_1+ m\theta_2)} .
\end{equation}
Substituting these Fourier series into \eqref{eq:secondInvEquation} yields,
\[
    (\omega_1 n + m) k_{n,m} = J_0^\prime n h_{n,m}.
\]
The solution for $k_{n,m}$ is singular whenever $\omega_1 n +  m = 0$, the
so called `resonances'. However, McNamara had the central idea that, provided
one knows which resonances appear in $H_1$ a priori, then, through judicious
choice of the free function $J_0^\prime$, it is possible to ensure the right hand side vanishes at
resonance, thus obtaining a nonsingular value for $J_1$.

Let us apply this theory to the two-wave model. Here we
extend the phase space by adding the (negative) energy variable
$e$, so that $H \to H +e$ and the canonical form becomes $\omega = dp\wedge dq +
de\wedge dt$. Setting $(\theta_1,I_1,\theta_2,I_2) = (q,p,t,e)$, 
$\varepsilon = \mu$ and
\begin{align*}
	H_0 &= \tfrac12 p^2 + e, \\
	H_1 &= -\cos(2\pi q) -\nu\cos(2\pi k(q-t)). 
\end{align*}
so that $e$ is the negative of the time-dependent energy of the original system \eqref{eq:twowaveHam}.
Taking,
\[
	J_1 = A(p)\cos(2\pi q) + B(p)\cos(2\pi k (q-t)),
\]
to match the harmonics of $H_1$, then \eqref{eq:secondInvEquation} is solved
provided,
\begin{equation}
  \label{eq:ABsol}
  A = -\frac{J_0^\prime(p)}{p},\qquad B = -\nu \frac{J_0^\prime(p)}{p-1}.
\end{equation}
It is then evident that choosing, for instance, $J_0^\prime(p) = p(p-1) $
will yield a nonsingular $J_1$. In summary, we obtain---to first order---that $J$ is
given by \eqref{eq:approxInv}.

Of course, it is important to ask whether one can use a similar method to obtain
higher orders in $\varepsilon$. This is possible to do directly as above;
however, making use of some normal form theory and Lie transforms allows for
some ease in describing how to do this. This higher order version was first
realized in \cite{mcnamara_superconvergent_1978}.

The primary aim of this higher order invariant is to transform the Hamiltonian
\eqref{eq:HamSeries}
to the Hamiltonian averaged over the $H_0$ trajectories,
\[\bar{H}(\phi_1,\bar{I}_1,\phi_2,\bar{I}_2) = H_0(\bar{I}_1,\bar{I}_2) + \varepsilon \bar{H}_1(\bar{I}_1,\bar{I}_2) + \dots,\]
at least to some order. If this is achieved then the resulting coordinates
$\bar{I}_1,\bar{I}_2$ will be nearly invariant and writing $\bar{I}_1 =
\bar{I}_1(\theta_1,I_1,\theta_1,I_2)$ gives the desired invariant in the old
coordinates.

A convenient method for obtaining the averaged Hamiltonian is through normal
form theory and Lie transforms. Given a function $W:M\to \R$ one can generate
the associated Hamiltonian vector field $X_W$ through 
\[
	X_W:=\{W,\cdot\}.
\]
The time-$\varepsilon$ flow of such a vector field, $\phi_W^\varepsilon =
e^{\varepsilon X_W} $, is a canonical transformation from the old coordinates
$z$ to new coordinates, say $Z = \phi_W^\varepsilon (z)$.
This gives a remarkable connection between functions and near-identity canonical
transformations. Moreover, pulling back a function $F$ by
$\phi_W^{\varepsilon}$, the so called Lie transformation, can be computed from
the exponentiation of the vector field,
\[ \phi_W^{\varepsilon *} F = \me^{\varepsilon X_W}F.\] This is a consequence of
the fact that the space of Hamiltonian vector fields $X$ is the Lie algebra
associated to the Lie group of canonical transformations. The pullback of a
function $F$ by $\phi_W^{\varepsilon *}$ can be written in the form
\[ 
	\phi_W^{\varepsilon *} F = \me^{\varepsilon \cL_W}F. 
\]
The Lie transformation view gives a quick method for computing transformed
a Hamiltonian under a canonical transform generated from some $W$. In
particular, to carry out a normalization at order $\varepsilon^i$, we write,
\[H = H_0 + \varepsilon H_1 + \dots, \qquad W = \varepsilon^{i-1}W_{i} ,\]
then the Lie transform method yields the transformed Hamiltonian as,
\begin{align*}
  \phi_W^{*}H &= \me^{X_W}H \\
              &= \me^{\varepsilon^{i}\{W,\cdot\}}(H_0 + \varepsilon H_1 + \dots) \\
              &= H_0 + \varepsilon H_1 + \dots + \varepsilon^{i}(\{W_i,H_0\} + H_i) + \dots .
\end{align*}
Hence, if $\phi_W^{*}H := \bar{H}$ then the $\bar{H}_j = H_j$ for $j < i$ and at
order $i$ we have,
\begin{equation}\label{eq:Homological}
	\{H_0,W_i\} = H_i - \bar{H}_i,
\end{equation}
the so call homological equation.

This observation gives a method for computing the averaged Hamiltonian $\bar{H}$
by sequentially solving each order $i$. The operator $\cL_{H_0}:= \{H_0,\cdot\}$ is
known as the cohomological operator. For the two-wave example it is given by,
\[
	\cL_{H_0} = p \partial_q + \partial_t.
\]
Equation \eqref{eq:Homological} can be solved by taking Fourier series as in \eqref{eq:FSeries}.
However the solution will be singular at resonance.

To obtain an integral that is not singular, despite the averaged Hamiltonian
being singular, we can apply the inverse transformation generated by the $W_i$
to an arbitrary function $\tilde{J}_0(\bar{I}_1)$. The inverse transformation is acquired
through the composition of each of the $-W_i$ generated transformations. In
particular, the lowest order terms would be,
\begin{align*}
	J(\theta_1,I_1,\theta_2,I_2) &=  e^{-\varepsilon^2 X_{W_2}}e^{-\varepsilon X_{W_1}} \tilde{J}_0 \\
	 &= \tilde{J}_0(I_1) - \varepsilon \{W_1,\tilde{J}_0\} + \varepsilon^2
	\left( -\{W_2,\tilde{J}_0\} + \tfrac12 \{W_1,\{W_1,\tilde{J}_0\}\}  \right) + \dots. 
\end{align*}
Hence, in order $J_0$ be nonsingular, $\tilde{J}_0$ must be chosen to cancel
the singularities in both $W_2$ and $W_1$.

Using the outlined procedure a second order adiabatic invariant for the two-wave model can be computed as,
\begin{equation}
  \label{eq:invariant2}
  \begin{aligned}
    J &= \tfrac{1}{4} (p-1)^4 p^4 - \mu \left(p^2 (2 p-1) (p-1)^3 \cos (2 \pi q)+
      \nu p^3 (2 p-1) (p-1)^2 \cos (2 \pi(q-t))\right)\\
    &\qquad+\tfrac{1}{4} \mu ^2 \left[\nu ^2 p^2 (10 p^2 -12 p + 3)+
      \nu^2 p^2 (3 p-1) (4 p-3) \cos (4 \pi (q-t)) \right.\\
    &\qquad\quad +4 \nu p (p-1) (6p^2 - 6p + 1) \cos (2 \pi( 2q- t))
    +\nu p (5p^2 - 5p +1) \cos (2 \pi t)) \\
    &\qquad \quad\left. +(3 p-2) (4 p-1) (p-1)^2 \cos (4 \pi
      q)+(10p^2-8p+1) (p-1)^2\right] \end{aligned}
\end{equation}
	\bibliography{ConverseKAM}
	\bibliographystyle{alpha}	
\end{document}